\documentclass[12pt,x11names]{article}

\usepackage{polyglossia}
\setdefaultlanguage[variant=british]{english}

\usepackage{amsmath,amsthm,eucal,array}

\usepackage{appendix}

\usepackage[T1]{fontenc} 
\usepackage{enumerate}
\usepackage{textcomp}
\usepackage[obeyspaces,lowtilde]{url}

\usepackage{tikz}
\usetikzlibrary{calc}
\usepackage{tikz-cd}
\usetikzlibrary{arrows,decorations.markings}

\tikzset{commutative diagrams/.cd,arrow style=tikz,diagrams={>=stealth'}}

\DeclareFontFamily{U}{mathb}{\hyphenchar\font45}
\DeclareFontShape{U}{mathb}{m}{n}{
      <5> <6> <7> <8> <9> <10> gen * mathb
      <10.95> mathb10 <12> <14.4> <17.28> <20.74> <24.88> mathb12
      }{}
\DeclareSymbolFont{mathb}{U}{mathb}{m}{n}
\DeclareMathSymbol{\bigast}{1}{mathb}{"06}

\usepackage{newpxtext} 
\usepackage[type1]{cabin} 
\usepackage[bigdelims,vvarbb]{newpxmath} 
\usepackage[bb=boondox,scr=rsfso]{mathalfa} 
\usepackage[%
hyperfootnotes=false,
bookmarks=true,
breaklinks=true,
bookmarksnumbered = true,
colorlinks= true,
urlcolor= blue,
anchorcolor = yellow,
citecolor=blue,
pdfborder= {0 0 0}
]{hyperref}

\renewcommand{\epsilon}{\ensuremath{\varepsilon}}
\renewcommand{\to}{\ensuremath{\longrightarrow}}
\renewcommand{\mapsto}{\ensuremath{\longmapsto}}
\newcommand{\rp}[1][2]{\ensuremath{{\vvmathbb R}P^{#1}}}
\newcommand{\Et}{\ensuremath{\vvmathbb E}^{\,2}}
\newcommand{\N}{\ensuremath{\vvmathbb N}}
\newcommand{\Z}{\ensuremath{\vvmathbb Z}}
\newcommand{\dt}{\ensuremath{\vvmathbb D}^{\,2}}
\newcommand{\St}[1][2]{\ensuremath{\vvmathbb S}^{#1}}
\newcommand{\FF}{\ensuremath{\vvmathbb F}}
\newcommand{\F}[1][2]{\ensuremath{\FF_{{#1}}}}
\newcommand{\sn}[1][n]{\ensuremath{S_{{#1}}}}
\renewcommand{\ker}[1]{\ensuremath{\operatorname{\text{Ker}}({#1})}}

\newcommand{\pnm}[1][n]{\ensuremath{P_{{#1}}(M)}}

\newcommand{\Pa}{\ensuremath{\mathsf{Pa}}}

\makeatletter

\renewcommand{\theenumi}{\alph{enumi}}
\renewcommand{\theenumii}{\roman{enumii}}
\renewcommand{\labelenumi}{(\theenumi)}
\renewcommand{\labelenumii}{(\theenumii)}
\renewcommand{\p@enumii}{}

\def\@enum@{\list{\csname label\@enumctr\endcsname}%
           {\usecounter{\@enumctr}\def\makelabel##1{
\normalfont\ignorespaces\emph{{##1}~}}
\setlength{\labelsep}{3pt}
\setlength{\parsep}{0pt}
\setlength{\itemsep}{0pt}
\setlength{\leftmargin}{0pt}
\setlength{\labelwidth}{0pt}
\setlength{\listparindent}{\parindent}%
\setlength{\itemsep}{0pt}
\setlength{\itemindent}{0pt}
\setlength{\topsep}{3pt plus 1pt minus 1 pt}}}

\def\@map#1#2[#3]{\mbox{$#1 \colon\thinspace #2 \to #3$}}
\def\map#1#2{\@ifnextchar [{\@map{#1}{#2}}{\@map{#1}{#2}[#2]}}
\DeclareRobustCommand*\textsubscript[1]{\@textsubscript{\selectfont#1}}
\def\@textsubscript#1{{\m@th\ensuremath{_{\mbox{\fontsize\sf@size\z@#1}}}}}
\makeatother

\newcommand{\id}{\ensuremath{\operatorname{\text{Id}}}}
\newcommand{\aut}[1]{\ensuremath{\operatorname{\text{Aut}}(#1)}}
\newcommand{\im}[1]{\ensuremath{\operatorname{\text{Im}}(#1)}}
\newcommand{\lhra}{\mathrel{\lhook\joinrel\to}}
\newcommand{\diff}[1]{\ensuremath{\operatorname{\mathcal{H}^+}\!\!\left({#1}\right)}}
\newcommand{\mcg}[2][n]{\ensuremath{\operatorname{\text{MCG}}(#2,#1)}}
\renewcommand{\gcd}{\ensuremath{\operatorname{\text{gcd}}}}

\DeclareMathOperator{\ppcm}{\text{ppcm}}

\DeclareRobustCommand*{\up}[1]{\textsuperscript{#1}}
\renewcommand{\th}{\ensuremath{\up{th}}}

\newcommand{\garside}[1][n]{\ensuremath{\Delta_{#1}}}
\newcommand{\ft}[1][n]{\ensuremath{\Delta_{#1}^{2}}}
\newcommand{\ftalt}[2][n]{\ensuremath{\Delta_{#1}^{2#2}}}

\newcommand{\brak}[1]{\ensuremath{\left\{ #1 \right\}}}
\newcommand{\ang}[1]{\ensuremath{\left\langle #1\right\rangle}}
\newcommand{\set}[2]{\ensuremath{\brak{#1 \,\mid\, #2}}}
\newcommand{\ord}[1]{\ensuremath{\left\lvert #1\right\rvert}}
\newcommand{\floor}[1]{\ensuremath{\left\lfloor #1\right\rfloor}}

\newcommand{\setang}[2]{\ensuremath{\ang{#1 \,\mid\, #2}}}
\newcommand{\setangr}[2]{\ensuremath{\ang{#1 \,\left\lvert \, #2 \right.}}}
\newcommand{\setangl}[2]{\ensuremath{\ang{\left. #1 \,\right\rvert \, #2}}}

\newcommand{\setr}[2]{\ensuremath{\brak{#1 \,\left\lvert \, #2 \right.}}}
\newcommand{\setl}[2]{\ensuremath{\brak{\left. #1 \,\right\rvert \, #2}}}

\newcommand{\an}[1][n]{\ensuremath{A_{{#1}}}}

\newcommand{\dih}[1]{\ensuremath{\operatorname{\text{Dih}}_{#1}}}
\newcommand{\dic}[1]{\ensuremath{\operatorname{\text{Dic}}_{#1}}}
\newcommand{\quat}[1][8]{\ensuremath{\mathcal{Q}_{#1}}}

\newcommand{\tonestar}{\ensuremath{\operatorname{T}^{\ast}}}
\newcommand{\oonestar}{\ensuremath{\operatorname{O}^{\ast}}}
\newcommand{\istar}{\ensuremath{\operatorname{I}^{\ast}}}

\newtheoremstyle{theoremm}{}{}{\itshape}{}{\scshape}{.}{ }{}

\theoremstyle{theoremm}
\newtheorem{thm}{Theorem}
\newtheorem{lem}[thm]{Lemma}
\newtheorem{prop}[thm]{Proposition}
\newtheorem{cor}[thm]{Corollary}

\newtheoremstyle{remarkk}{}{}{}{}{\scshape}{.}{ }{}

\theoremstyle{remarkk}
\newtheorem{defn}[thm]{Definition}
\newtheorem{rem}[thm]{Remark}
\newtheorem{rems}[thm]{Remarks}

\newcommand{\reth}[1]{Theorem~\protect\ref{th:#1}}
\newcommand{\relem}[1]{Lemma~\protect\ref{lem:#1}}
\newcommand{\repr}[1]{Proposition~\protect\ref{prop:#1}}
\newcommand{\reco}[1]{Corollary~\protect\ref{cor:#1}}
\newcommand{\resec}[1]{Section~\protect\ref{sec:#1}}
\newcommand{\req}[1]{equation~(\protect\ref{eq:#1})}
\newcommand{\reqref}[1]{(\protect\ref{eq:#1})}
\newcommand{\rerem}[1]{Remark~\protect\ref{rem:#1}}
\newcommand{\rerems}[1]{Remarks~\protect\ref{rems:#1}}

\newcommand{\proclaim}[2]{\par\vspace{\partopsep}\vspace{\partopsep}\noindent\textsc{#1}
\emph{{#2}}\par\vspace{\partopsep}\vspace{\partopsep}}

\newcommand{\sii}[2][1]{\ensuremath{\sigma_{#2}^{-{#1}}}}

\newcommand{\aiji}[2]{\ensuremath{B_{#1,#2}^{-1}}}
\newcommand{\si}[2][{}]{\ensuremath{\sigma_{#2}^{#1}}}
\newcommand{\rh}[2][{}]{\ensuremath{\rho_{#2}^{#1}}}
\newcommand{\rhi}[2][1]{\ensuremath{\rho_{#2}^{-{#1}}}}
\newcommand{\rhs}[2]{\ensuremath{\rho_{#1}\cdots\rho_{#2}}}
\newcommand{\rhsi}[2]{\ensuremath{\rho_{#1}^{-1}\cdots\rho_{#2}^{-1}}}

\newcommand{\ssni}[2]{\ensuremath{\sigma_{#1}^{-1}\cdots\sigma_{#2}^{-1}}}

\newcommand{\ms}{\medskip}
\newcommand{\bs}{\bigskip}

\newcommand{\comj}[1]{\noindent\textcolor{DarkOrchid3}{\textbf{!!!J!!!~#1}}}

\def\bmod{\nonscript\mskip-\medmuskip\mkern5mu
  \mathbin{\text{mod}}\penalty900\mkern5mu\nonscript\mskip-\medmuskip}

\begin{document}

\title{The structure of the finite subgroups of the braid groups of the real projective plane}

\author{DACIBERG~LIMA~GON\c{C}ALVES\\
Departamento de Matem\'atica - IME-USP,\\
Rua~do~Mat\~ao~1010~CEP:~05508-090  - S\~ao Paulo - SP - Brazil.\\
e-mail:~\texttt{dlgoncal@ime.usp.br}\vspace*{4mm}\\
JOHN~GUASCHI\\
Université Caen Normandie, CNRS, Normandie Univ\\
LMNO UMR6139, F-14000 CAEN, FRANCE.\\
e-mail:~\texttt{john.guaschi@unicaen.fr}}

\date{\today}

\begingroup
\renewcommand{\thefootnote}{}
\endgroup


\maketitle

\begin{abstract} 
\noindent
In this paper we study the structure of the finite subgroups of the braid groups $B_n(\rp)$ of the real projective plane $\rp$. It is known that the isomorphism classes of these subgroups are comprised of cyclic, dicyclic and binary polyhedral groups. Using specific elements of $B_{n}(\rp)$, we construct embeddings of the cyclic and dicyclic subgroups, we show that there is a single conjugacy class of the maximal dicyclic groups, and for small values of $n$, we also construct embeddings of the binary tetrahedral and octahedral groups. One interesting phenomenon that occurs is the existence of isomorphism classes of finite subgroups that admit embeddings in $B_{n}(\rp)$ as both maximal finite and non-maximal subgroups. Up to conjugacy, we describe completely the lattice of the finite subgroups of $B_n(\rp)$ for $n=3,4,5$. Some of our results regarding the number of conjugacy classes of finite subgroups are obtained using the relation between $B_{n}(\rp)$ and the mapping class group of $\rp$ relative to $n$ points. We also provide a new presentation of $B_n(\rp)$ which has some interesting properties with respect to the classical presentation of Van Buskirk.
\end{abstract}

\section{Introduction}\label{sec:intro}

The study of the finite subgroups of a given group 
constitutes an interesting and natural problem. In this paper, we will study several aspects of the finite subgroups of  $B_n(\rp)$, the full $n$-string braid group of the real projective plane $\rp$. The following presentation of $B_{n}(\rp)$ is due to Van Buskirk.


\begin{prop}[Van Buskirk~\cite{vB}]\label{prop:present}
Let $n\geq 1$. The following constitutes a presentation of the group $B_n(\rp)$:
\begin{enumerate}
\item[\underline{\textbf{generators:}}]
$\si{1},\ldots,\si{n-1},\rh{1},\ldots,\rh{n}$.
\item[\underline{\textbf{relations:}}]
\begin{align}
\si{i}\si{j} &=\si{j}\si{i}\quad\text{if $\lvert i-j\rvert\geq 2$,} \label{eq:artin1}\\
\si{i}\si{i+1}\si{i}&=\si{i+1}\si{i}\si{i+1} \quad\text{for $1\leq i\leq n-2$,} \label{eq:artin2}\\
\si{i}\rh{j}&=\rh{j}\si{i}\quad\text{for $j\neq i,i+1$,}\label{eq:sirj}\\
\rh{i+1}&=\sii{i}\rh{i}\sii{i} \quad\text{for $1\leq i\leq n-1$,}\label{eq:sirisi}\\
\rh[-1]{i+1}\rh[-1]{i}\rh{i+1}\rh{i}&= \si[2]{i} \quad\text{for $1\leq i\leq n-1$,}\label{eq:rhocomm}\\
\rh[2]{1}&=\si{1}\si{2}\cdots\si{n-2}\si[2]{n-1} \si{n-2}\cdots\si{2}\si{1}.\label{eq:surfacerp2}
\end{align}
\end{enumerate}
\end{prop}
The generators $\sigma_{1},\ldots,\sigma_{n-1}$ are the usual Artin generators obtained via the embedding of a topological disc in $\rp$, while for $i=1,\ldots,n$, $\rho_{i}$ corresponds to a generator of the fundamental group of $\rp$ with basepoint at the $i$th point of the $n$ basepoints of the elements of $B_{n}(\rp)$. Surface braid groups may be defined topologically via configuration spaces, but for the purposes of this paper, we shall take the presentation of \repr{present} to be the definition of $B_{n}(\rp)$. We will refer to relations~\reqref{artin1} and~\reqref{artin2} as the \emph{Artin relations}, and relation~\reqref{surfacerp2} as the \emph{surface relation} of $B_n(\rp)$. We define the \emph{half-twist braid} $\garside$ of $B_{n}(\rp)$ by $\garside= \prod_{i=1}^{n-1} \sigma_{1} \cdots \sigma_{n-i}$. Viewed as an element of the Artin braid group $B_{n}$, $\garside$ is the well-known Garside element. One may check that the \emph{full-twist braid} $\ft$ satisfies $\ft=(\sigma_{1}\cdots \sigma_{n-1})^{n}$. 

It is clear that $B_1(\rp)$ is isomorphic to the cyclic group $\Z_2$ of order $2$. By~\cite[Theorem, page~94]{vB}, $B_2(\rp)$ is isomorphic to the generalised quaternion group of order $16$, and if $n\geq 3$ then $B_n(\rp)$ is infinite. Recall that if $m\geq 2$, the \emph{dicyclic group of order $4m$}, denoted by $\dic{4m}$, is defined by the following group presentation:
\begin{equation*}
\dic{4m}= \setangl{x,y}{\text{$x^m=y^2$ and $yxy^{-1}$}}.
\end{equation*}
If $m$ is a power of $2$, then $\dic{4m}$ is also known as the \emph{generalised quaternion group of order $4m$}, and will be denoted in what follows by $\quat[4m]$. If $m=2$ then $\quat[8]$ is the usual quaternion group of order $8$.


One of the questions that will interest us in this paper is that of the possible inclusions (up to conjugacy) among the finite subgroups of $B_n(\rp)$. To this end, we define the following notions of maximal and abstractly maximal finite subgroups.

\begin{defn}
Let $H$ be a finite group, and let $\iota\colon\thinspace H\lhra G$ be an embedding of $H$ in another group $G$. The embedding $\iota$ 
is said to be \emph{maximal} if for any finite subgroup $K$ of $G$ such that  $\iota(H)\subset K$, we have $\iota(H)=K$. If $H$ is itself a subgroup of $G$ then we use the phrase \emph{`$H$ is a finite maximal subgroup of $G$'} to mean that the inclusion $H\subset G$ is a maximal embedding.
\end{defn}


\begin{defn}
Let $H$ be a finite group that embeds in a group $G$. We shall say that \emph{$H$ is an abstractly maximal finite subgroup of $G$} 
if every embedding of $H$ in $G$ is maximal.
\end{defn}



Using the above terminology, we recall the following result that classifies the isomorphism classes of the finite subgroups of $B_n(\rp)$.

\begin{thm}\label{th:finitebnrp2}\cite[Theorem~5]{GG5}
Let $n\geq 2$. The abstractly maximal finite groups of $B_n(\rp)$ are as follows:
\begin{enumerate}
\item\label{it:finrpa} the dicyclic group $\dic{8n}$ of order $8n$.
\item\label{it:finrpb} 
the dicyclic group $\dic{8(n-1)}$ if $n\geq 4$.
\item\label{it:finrpc} the binary octahedral group, denoted by $\oonestar$, if $n\equiv 0,1\pmod 3$.
\item\label{it:finrpd} the binary icosahedral group, denoted by $\istar$, if $n\equiv
0,1,6,10\pmod{15}$.
\end{enumerate}
\end{thm}

More information about these groups may be found in~\cite{Wo} and in~\cite[Appendix]{GG8}. The groups $\oonestar$ and $\istar$ contain a subgroup that is isomorphic to the binary tetrahedral group $\tonestar$. 
So it follows from \reth{finitebnrp2} that if $\map{\iota}{\tonestar}[B_n(\rp)]$ is an embedding, then there exists an embedding $\iota_1$ of either $\oonestar$ or $\istar$ in $B_n(\rp)$ such that 
$\iota(\tonestar) \subset \im{\iota_1}$.
The groups $\tonestar,\oonestar$ and $\istar$ will be referred to collectively as the \emph{binary polyhedral groups}.

It follows from the above definitions that every abstractly maximal finite subgroup of a group $G$ is also maximal finite. One may ask whether the converse also holds. 
Theorem~1.3 of~\cite{GGagt2} characterises the maximal finite subgroups of $B_{n}(\St)$, and implies that if a finite group admits a maximal (resp.\ non-maximal) embedding in $B_{n}(\St)$ then any other embedding is also maximal (resp.\ non-maximal). In particular, the notions of being maximal finite subgroups and abstractly maximal finite subgroups are equivalent in $B_{n}(\St)$. Making use of the results of Sections~\ref{sec:subgroupsb3} and~\ref{sec:subgroupsb4} respectively, we will show that this is also the case in $B_3(\rp)$ and  $B_4(\rp)$.
However, as we shall see, if $n\geq 5$, $B_n(\rp)$ possesses maximal finite subgroups that are not abstractly maximal. One question that we will discuss in this paper is that of the characterisation of the finite subgroups of $B_n(\rp)$ that admit both maximal and non-maximal embeddings. As an example, we will show in \repr{nomax} that for all $n\geq 5$ there are dicyclic subgroups of $B_{n}(\rp)$ that are maximal finite, but not abstractly maximal finite. We also show in \resec{subgroupsB5} that there is an embedding of $\quat$ that is maximal but not abstractly maximal in $B_5(\rp)$.

We recall the following results
regarding the orders of the finite-order elements of $B_n(\rp)$.


\begin{prop}\label{prop:agt} \cite[Theorem~4 and Proposition~23]{GG2}  
Let $n\geq 2$. 
\begin{enumerate}
\item\label{it:agta} $B_n(\rp)$ has an element of order $\ell$ if and only if $\ell$ divides either $4n$ or $4(n-1)$.
\item\label{it:agtb} The full twist braid $\ft$ is the unique element of $B_n(\rp)$ of order~$2$, and it generates the centre of $B_n(\rp)$.
\end{enumerate} 
\end{prop}
It follows that any finite cyclic subgroup of $B_{n}(\rp)$ is isomorphic to a subgroup of $\Z_{4(n-i)}$ for some $i\in \brak{0,1}$. 

Apart from this introduction, this paper is comprised of seven sections and two appendices. In \resec{prelim}, we recall some relations between special elements of $B_{n}(\rp)$ that will allow us in what follows to construct explicit representatives of certain finite subgroups.

In \resec{conjclass}, we initiate our study of the conjugacy classes of finite subgroups of $B_n(\rp)$. It is well known that surface braid groups are closely related to the corresponding mapping class groups. In the case of $\rp$, if $n\geq 2$, we have the following short exact sequence~\cite[p.~617]{Sc}:
\begin{equation}\label{eq:bnrpmcg}
1\to \ang{\ft} \to B_n(\rp) \stackrel{\phi} \to \mcg{\rp} \to 1,
\end{equation}
where $\mcg{\rp}$ denotes the mapping class group of the $n$-punctured projective plane, and where the homomorphism $\phi$ may be viewed as the canonical projection ($\ang{\ft}$ is a normal subgroup of $B_{n}(\rp)$ of order $2$ by \repr{agt}(\ref{it:agtb})). Consequently, the short exact sequence~\reqref{bnrpmcg} may be used to compare the finite subgroups of $B_{n}(\rp)$ and $\mcg{\rp}$. Let $\phi^{\ast}$ denote the map induced by $\phi$ from the set of subgroups of $\mcg{\rp}$ to the set of subgroups of $B_n(\rp)$ given by $\phi^{\ast}(H)=\phi^{-1}(H)$, where $H$ is a subgroup of $\mcg{\rp}$. 

\begin{rems}\label{rems:finitemcgbnrp}\mbox{}
\begin{enumerate}
\item\label{it:finitemcgbnrpa} Under $\phi$, the finite subgroups of $B_{n}(\rp)$ isomorphic to $\dic{8n},\dic{8(n-1)}, \tonestar, \oonestar$ and $\istar$ are sent to finite subgroups of $\mcg{\rp}$ isomorphic to $\dih{4n}, \dih{4(n-1)}, A_{4}, \sn[4]$ and $A_{5}$ respectively, where $\dih{m}$ denotes the dihedral group of order $m$~\cite[Section~ 3.2]{GG5}, which using~\reqref{bnrpmcg}, implies that under $\phi^{\ast}$, the finite subgroups of $\mcg{\rp}$ isomorphic to $\dih{4n}, \dih{4(n-1)}, A_{4}, \sn[4]$ and $A_{5}$ are sent to the finite subgroups of $B_{n}(\rp)$ isomorphic to $\dic{8n},\dic{8(n-1)}, \tonestar, \oonestar$ and $\istar$ respectively. 

\item\label{it:finitemcgbnrpb} 
Using \cite[Theorem~2.1]{BCG} (resp.\ \reth{finitebnrp2}), one may check that the finite subgroups of $\mcg{\rp}$ (resp.\ of $B_{n}(\rp)$) of odd order (resp.\ of even order not divisible by $4$ or of odd order) are cyclic. 
\end{enumerate}
\end{rems}

In \resec{conjclass} we prove the following result.

\begin{prop}\label{prop:BraMap}
Let $n\geq 2$, and consider the map $\phi^{\ast}$ defined above.
\begin{enumerate}
\item\label{it:BraMapa} The restriction of $\phi^{\ast}$ to the set of finite subgroups of $\mcg{\rp}$ in the domain and to the set of finite subgroups of $B_n(\rp)$ of even order in the target is a bijection, and it induces a bijection between the corresponding sets of conjugacy classes of these finite subgroups. 

\item\label{it:BraMapb} Let
$K=\langle x \rangle$ be a cyclic subgroup 
of $\mcg{\rp}$ of odd order $d$. Then $\phi^{-1}(x)$ consists of two elements, one of which, $x_{1}$ say, is of order $d$, and the other, $x_{2}$ say, is of order $2d$. Further, $\phi^{\ast}(K)=\ang{x_{2}}$, and the correspondence $K \longmapsto \phi^{\ast}(K)$ (resp.\ $K \longmapsto \ang{x_{1}}$) is a bijection between the finite subgroups of $\mcg{\rp}$ of odd order and the finite subgroups of $B_n(\rp)$ of even order not divisible by $4$ (resp.\ the finite subgroups of $B_n(\rp)$ of odd order), and it preserves conjugacy classes.
\end{enumerate}
\end{prop}
 
 

In order to better understand the structure of the finite subgroups of $B_n(\rp)$, and to decide for example which subgroups are maximal finite, it is necessary to analyse certain properties of their conjugacy classes, such as the number and representatives of such classes. From a combinatorial point of view in $B_n(\rp)$, this seems to be a complicated problem. Using just the structure of $B_3(\rp)$, we have been able to carry out this analysis completely for $n=3$ in \resec{subgroupsb3} (see Figure~\ref{fig:omegai} for the subgroup lattices of $B_3(\rp)$). For $n\geq 4$, the precise computation of the number of conjugacy classes of $\quat$ is more difficult. But by \repr{BraMap}, the number of conjugacy classes of a given finite subgroup of $B_n(\rp)$ of even order is equal to that of the corresponding subgroup of $\mcg{\rp}$, and so may be computed using the results of~\cite{BCG} regarding the conjugacy classes of the finite subgroups of $\mcg{\rp}$. Some examples of such computations are given in \ref{app:quatbnrp2} for $n=3,4$ and $5$. This extra information allows us to determine completely the structure of the conjugacy classes of the finite subgroups of $B_n(\rp)$ for $n=4,5$ in Sections~\ref{sec:subgroupsb4} and~\ref{sec:subgroupsB5} respectively (see Figures~\ref{fig:omegai},~\ref{fig:omegaii} and~\ref{fig:omegaiii} for the corresponding subgroup lattices). 

For $n\geq 6$, the analysis of the structure of the finite subgroups of $B_n(\rp)$ becomes considerably more difficult. Even though we are able to compute the number of conjugacy classes of the isomorphism classes of the finite subgroups of $B_n(\rp)$ using~\cite{BCG}, the explicit realisation of their representatives in terms of the generators of $B_n(\rp)$ seems to be a very difficult problem in general. However, the structure of the conjugacy classes of the elements of finite order is understood.
Let $n\geq 2$. Applying Murasugi's classification~\cite{M} and its simplification~\cite[Proposition~8 and Corollary~12]{GG3} and making use of the presentation of \repr{present}, up to conjugacy and powers, every element of $B_{n}(\rp)$ of finite order is of the form: 
\begin{equation}\label{eq:muragen}
\delta_{r,j}=\bigl( \rho_{r}\sigma_{r-1}\cdots \sigma_{1}\bigr)^{2r/l} \bigl( \sigma_{r+1}\cdots \sigma_{n-1}\sigma_{r+1}^{j}\bigr)^{(n-j-r)/l},
\end{equation}
where $j\in \brak{0,1}$, $0\leq r \leq n-j$, $l=\gcd(n-j-r,2r)$, and $\delta_{r,j}$ is of order $2l$. In what follows, we will refer to the elements of the form~\reqref{muragen} and their powers as
\emph{Murasugi elements}.
Note that by \repr{present}, the two bracketed factors of \req{muragen} commute.

The following two results will also be proved in \resec{conjclass}.
 
\begin{thm}\label{th:ccrp2}
Let $n\geq 3$, and let $d$ be a divisor of $4(n-i)$, where $i\in \brak{0,1}$.
\begin{enumerate}[(a)]
\item\label{it:ccrp2a} 
Suppose that $d$ is divisible by $4$. If $x$ and $y$ are Murasugi elements of $B_n(\rp)$ of order $d$ then they are conjugate if and only if $y\in \brak{x,x^{-1}}$. 

\item\label{it:ccrp2b} 
Let $n\in \brak{3,4, 5}$. If $d$ is not divisible by $4$ then the cyclic subgroups of $B_n(\rp)$ of order $d$
are pairwise conjugate.
\end{enumerate}
\end{thm}



We have a proof of \reth{ccrp2}(\ref{it:ccrp2b}) for all $n\geq 3$, but it is currently extremely long. The result may also be proved using~\cite{BCG}. To illustrate our method, in this paper we will just prove the result for $n=3,4,5$.
With respect to part~(\ref{it:ccrp2a}), we may determine the exact number of conjugacy classes of cyclic subgroups of $B_{n}(\rp)$ of order $d$, where $d$ is a multiple of $4$, as follows.



\begin{cor}\label{cor:ccorder4} 
Let $n\geq 3$. 
\begin{enumerate}
\item\label{it:ccorder4a} The group $B_{n}(\rp)$ possesses $n+1$ conjugacy classes of subgroups isomorphic to $\Z_{4}$.
\item\label{it:ccorder4b} If $j\in \brak{0,1}$, let $d\geq 8$ be a multiple of $4$ that divides $4(n-j)$. Then the group $B_{n}(\rp)$ possesses $\lfloor 2(n-j)/d \rfloor+1$ conjugacy classes of subgroups isomorphic to $\Z_{d}$.
\end{enumerate}
\end{cor}

One may check that the results of \reco{ccorder4} are coherent with those of \cite{BCG}.



In \resec{inclcyclic}, we turn our attention to the exploration of the properties of the non-cyclic subgroups of $B_{n}(\rp)$, namely the dicyclic and binary polyhedral groups. We will be particularly interested in the abstractly maximal dicyclic subgroups $\dic{8(n-i)}$, where $i=0,1$ for which we show in \repr{conjdicmax} that there is a single conjugacy class of each. In \repr{centraliserx}, we will determine the centralisers of certain powers of two elements $a=\rho_{n} \sigma_{n-1}\cdots \sigma_{1}$ and $b=\rho_{n-1} \sigma_{n-2}\cdots \sigma_{1}$ that play an important rôle in the construction of finite subgroups of $B_{n}(\rp)$. The corresponding analysis of the binary polyhedral subgroups is much more difficult, but in \repr{binpolyI} we are able to obtain some information about their embedding in $B_{n}(\rp)$ for $n\geq 4$, notably with respect to their intersection with $P_{n}(\rp)$. As indicated above, in \repr{nomax}, we show that $B_{n}(\rp)$ contains dicyclic subgroups of order $4(n-i)$, where $i=0,1$, that are maximal but not abstractly maximal.

As we mentioned previously, in Sections~\ref{sec:subgroupsb3}--\ref{sec:subgroupsB5}, we determine the number and structure of the conjugacy classes of the finite subgroups of $B_{n}(\rp)$ for $n=3,4$ and $5$ respectively, and we give representatives in terms of the generators on the presentation of \repr{present}. In the case $n=3$, we are able to do this without reference to~\cite{BCG}. For $B_{3}(\rp)$ and $B_{4}(\rp)$, we also determine the centralisers of the elements of finite order, with the exception of two elements of $B_{4}(\rp)$. 

The final section of the paper, \resec{newpres}, is devoted to giving a new presentation of $B_n(\rp)$ in \reth{newpresrp2} that possesses some interesting features, highlighted in \repr{goodpropsnew}, that are distinct from those of \repr{present}. 

The structure of the conjugacy classes of the finite subgroups of $P_{n}(\rp)$ is completely understood. The group $P_{1}(\rp)$ is isomorphic to $\Z_{2}$, $P_{2}(\rp)$ is isomorphic to $\quat$~\cite[4th lemma, p.~87]{vB}. If $n\geq 3$, $P_{n}(\rp)$ is infinite, and as for $B_{n}(\rp)$ contains a unique subgroup $\ang{\ft}$ of order $2$. The remaining finite subgroups of $P_{n}(\rp)$ are isomorphic to $\Z_{4}$ and $\quat$ if $n=3$, and to $\Z_{4}$ if $n\geq 4$~\cite[Proposition~1]{GG3}. In $P_{3}(\rp)$, there is a single conjugacy class of subgroups isomorphic to $\quat$~\cite[Proposition 21(4)]{GG2}. In~\cite[Proposition 21(3)]{GG2}, we also determined the number of conjugacy classes in $P_{3}(\rp)$ of subgroups isomorphic to $\Z_{4}$, but as we pointed out in~\cite[Remark~10]{GGzeit}, there was an error in the statement and proof of~\cite[Proposition 21(3)]{GG2}. In \ref{app:conjclass}, we take the opportunity in \repr{5conjclasses} to correct this mistake as well as the proof of the fact that $P_{3}(\rp)$ does not possess any subgroup isomorphic to $\Z_{4} \times \Z$ (see \repr{p3noz4z}) whose proof in~\cite[Proposition~17]{GG3} made use of~\cite[Proposition 21(3)]{GG2}. 

In \ref{app:quatbnrp2}, making use of the results of~\cite{BCG} and applying combinatorial arguments, in \reth{GammanQ8} we give an explicit formula for the  number of conjugacy classes of subgroups of $B_n(\rp)$ isomorphic to $\quat$ for all $n\geq 2$, and in \repr{conjclass345} we also determine the number of conjugacy classes of the other non Abelian finite subgroups of $B_n(\rp)$ in the case that $n\in \brak{3,4,5}$. The number of such subgroups is coherent with our results in the cases $n=3,4$, and the results allow us to give a complete answer to the question of determining the number of conjugacy classes of the finite subgroups of $B_5(\rp)$.  

\subsection*{Acknowledgements}

The work on this paper took place during a number of visits, starting in 2015, of the first author to the Laboratoire de Math\'ematiques Nicolas Oresme, and of the second author to the Departamento de Matem\'atica do IME~--~Universidade de S\~ao Paulo. It was supported by the international Cooperation Capes-Cofecub project numbers Ma~733-12 (France) and Cofecub~1716/2012 (Brazil), the FAPESP Projeto Tem\'atico-FAPESP Topologia Alg\'ebrica, Geom\'etrica e Diferencial no.~2022/16455-6 (S\~ao Paulo-Brazil), and the CNRS-FAPESP project Sprint Fapesp no.~2024/19199-6. 
The first author thanks             
the Réseau Mathématique Franco-Brésilien for financial support for a   visit to the Laboratoire de Math\'ematiques Nicolas Oresme
in May 2023, and   
the second author thanks the Capes Print programme for financial support for two visits to the Universidade de S\~ao Paulo in March 2023 and March 2024.

\subsection*{AI declaration}

No sophisticated tools nor AI tools were used in the preparation or the writing of this paper.

\section{Preliminaries}\label{sec:prelim}

In this paper, we will frequently make use of certain special elements of $B_{n}(\rp)$ that we studied in previous papers. In terms of the generators of $B_n(\rp)$ of \repr{present}, in~\cite{GG2}, 
we obtained the following useful identities:
\begin{align}
\rh{j}&=\ssni{j-1}{1}\rh{1}\ssni{1}{j-1} \quad\text{for $j=1,\ldots, n$ (equation~(4),
p.~770)}\label{eq:rjr1}\\
\rho_{n}^{-2} &=\si{n-1}\cdots \si{2}\si[2]{1} \si{2}\cdots \si{n-1} \quad\text{(equations~(6)
and~(7), p.~771),}\label{eq:rn2}
\end{align}
and we exhibited elements of $B_n(\rp)$ of order $4n$ and of order $4(n-1)$ respectively:
\begin{equation}\label{eq:defab}
\left\{
\begin{aligned}
a &= \sii{n-1}\cdots\sii{1}\cdot\rh{1}=\rh{n} \si{n-1}\cdots \si{1}\\
b &= \sii{n-2}\cdots\sii{1}\cdot\rh{1}=\rh{n-1} \si{n-2}\cdots \si{1}
\end{aligned}\right.
\end{equation}
that satisfy:
\begin{equation}\label{eq:powerab}
\text{$a^{n}=\rho_{n}\cdots \rho_{1}$ and $b^{n-1}=\rho_{n-1}\cdots \rho_{1}$.}
\end{equation}
Together with the half-twist element $\garside$, the elements $a$ and $b$ play an important r\^{o}le in our combinatorial study of $B_{n}(\rp)$, notably via their action by conjugation on Van Buskirk's generators of $B_{n}(\rp)$. 
It is well known that:
\begin{equation}\label{eq:conjsigma}
\garside \sigma_{i} \garside^{-1}=\sigma_{n-i} \quad \text{for all $i=1,\ldots,n-1$.}
\end{equation}
In~\cite[pp.~864--65]{GG5}, just after Remark~14,  we showed that: 
\begin{align}
\garside^{-1} \rho_{i} \garside &=\rho_{n+1-i}^{-1} \quad \text{for all $i=1,\ldots,n$}\label{eq:conjrho}\\
\garside a \garside^{-1} &= a^{-1} \quad \text{and}\label{eq:conja}\\
\garside a^{-1} \ldotp b \ldotp a \garside^{-1} &= b^{-1}.\label{eq:conjb}
\end{align}

 
We recall the following results from~\cite[pp.~777--78]{GG2} and~\cite[Proposition~15]{GG5}, where  part~(\ref{it:embed}) also follows from equations~(\ref{eq:conja}) and~(\ref{eq:conjb}).

\begin{prop}\label{prop:permute} Let 
 $n\geq 2$. Then in $B_n(\rp)$, we have the following relations.
\begin{enumerate}[(a)]
\item\label{it:permsi} Conjugation by $a^{-1}$ permutes cyclically the following elements:
\begin{equation*}
\si{1},\ldots, \si{n-1}, a^{-1}\si{n-1}a, \sii{1},\ldots, \sii{n-1}, a^{-1}\sii{n-1}a.
\end{equation*}
\item\label{it:permrhi} Conjugation by $a^{-1}$ permutes cyclically the following elements:
\begin{equation*}
\rh{1},\ldots \rh{n}, \rhi{1},\ldots, \rhi{n}.
\end{equation*}
\item\label{it:conjbsigma} Conjugation by $b^{-1}$ permutes cyclically the following elements:
\begin{equation*}
\si{1},\ldots, \si{n-2}, b^{-1}\si{n-2}b, \sii{1},\ldots, \sii{n-2}, b^{-1}\sii{n-2}b
\end{equation*}
\item\label{it:embed} We have $\ang{a,\garside}\cong \dic{8n}$ and $\ang{b,\garside a^{-1}}\cong \dic{8(n-1)}$.
\end{enumerate}
\end{prop} 

The subgroups of $B_{n}(\rp)$ given in \repr{permute}(\ref{it:embed}) will be referred to as the \emph{standard copies} of $\dic{8n}$ and $\dic{8(n-1)}$ respectively. Since $\garside a^{-1}=a \garside$ and $\ang{a,\garside}= \ang{a,a\garside}$, they may also be expressed in the form $\ang{x,a\garside}$, where $x\in \brak{a,b}$. 

There is a surjective homomorphism $\map{\pi}{B_{n}(\rp)}[\sn]$, where $\sn$ is the symmetric group on $n$ letters, that to each element of $B_{n}(\rp)$ associates its permutation 
More precisely, for  $i=1,\ldots, n-1$ (resp.\ $j=1,\ldots,n$), $\pi$ sends $\sigma_{i}$ (resp.\ $\rho_{j}$) to the transposition $(i,i+1)$ (resp.\ to the identity). The kernel of $\pi$ is the $n$-string pure braid group of $\rp$, and is denoted by $P_{n}(\rp)$. We thus obtain the following short exact sequence:
\begin{equation}\label{eq:sessn}
1\to P_{n}(\rp) \to B_{n}(\rp) \stackrel{\pi}{\to} \sn \to 1.
\end{equation}
It follows from \repr{present} that if $n=1$ then $B_{1}(\rp)=P_{1}(\rp)\cong \Z_{2}$, and if $n\geq 2$, the Abelianisation of $B_{n}(\rp)$ is isomorphic to $\Z_{2} \oplus \Z_{2}$, where the Abelianisation homomorphism 
$\map{q}{B_{n}(\rp)}[\Z_{2} \oplus \Z_{2}]$ 
is given by:
\begin{equation}\label{eq:abmap}
\text{$q(\sigma_{i})=(\overline{1}, \overline{0})$ and $q(\rho_{j})=(\overline{0}, \overline{1})$ for all $1\leq i\leq n-1$ and $1\leq j\leq n$.}
\end{equation}

Recall from the proof of~\cite[Proposition~8]{GG6} that $P_{n}(\rp)$ is generated by the set $\setl{A_{i,j},\, \rho_{k}}{\text{$1\leq i<j\leq n$ and $1\leq k\leq n$}}$, where $A_{i,j}$ (which was denoted by $B_{i,j}$ in~\cite{GG6}) is defined in terms of the generators of \repr{present} by:
\begin{equation}\label{eq:defaij}
A_{i,j}=\sigma_{j-1}\cdots \sigma_{i+1}\sigma_{i}^{2} \sigma_{i+1}^{-1}\cdots \sigma_{j-1}^{-1}.
\end{equation}
Making use of the Artin relations, $A_{i,j}$ may also be expressed as follows:
\begin{equation}\label{eq:defaijalt}
A_{i,j}=\sigma_{i}^{-1}\cdots \sigma_{j-2}^{-1} \sigma_{j-1}^{2} \sigma_{j-2}\cdots \sigma_{i}.
\end{equation}

The following lemma will be used later in the paper.

\begin{lem}\label{lem:revrhocomm}
Let $n\geq 2$.
\begin{enumerate}
\item\label{it:revrhocomma} If $1\leq i\leq n-1$ then  $\rho_i \rho_{i+1} \rho_i^{-1} \rho_{i+1}^{-1}= \sigma_i^2$ in $B_n(\rp)$.
\item\label{it:revrhocommb} If $1\leq i\leq n-1$ then $\rho_{i+1} A_{i,i+1} \rho_{i+1}^{-1}= \rho_{i}^{-1} A_{i,i+1}^{-1} \rho_{i}$ and $\rho_{i+1}^{-1} A_{i,i+1} \rho_{i+1}= \rho_{i} A_{i,i+1}^{-1} \rho_{i}^{-1}$ in $B_n(\rp)$.
\item\label{it:revrhocommc} For all $1\leq i < k< j\leq n$, $\rho_{k} A_{i,j} \rho_{k}^{-1}= \rho_{i}^{-1}  C_{i,k} \rho_{i}  A_{i,j} \rho_{i}^{-1}  C_{i,k}^{-1} \rho_{i}$ in $B_n(\rp)$, where for $1\leq p < q \leq n$, $C_{p,q} = (\prod_{l=p+1}^{q} A_{p,l})(\prod_{l=p+1}^{q-1} A_{p,l})^{-1}$.
\end{enumerate}
\end{lem}

\begin{proof}
Let $n\geq 2$, and let $1\leq i\leq n-1$.
\begin{enumerate}
\item Conjugating~\reqref{rhocomm} by $\garside$ and using~\reqref{conjsigma} and~\reqref{conjrho}, we obtain $\rho_{n-i} \rho_{n-i+1} \rho_{n-i}^{-1} \rho_{n-i+1}^{-1}= \sigma_{n-i}^2$, and conjugating this relation by $a^{n-2i}$ and applying \repr{permute}(\ref{it:permsi}) and~(\ref{it:permrhi}), we obtain the required relation.
\item By part~(\ref{it:revrhocomma}), we see that $\rho_{i+1}\rho_{i}= \sigma_i^{-2} \rho_{i} \rho_{i+1}$, and using this, and equations~\reqref{sirisi} and~\reqref{rhocomm}, we have:
\begin{align*}
\rho_{i}\rho_{i+1} A_{i,i+1} \rho_{i+1}^{-1} \rho_{i}^{-1} &= \rho_{i}\sigma_{i}^{-1} \rho_{i} \sigma_{i} \rho_{i+1}^{-1} \rho_{i}^{-1} = \sigma_{i} \rho_{i+1} \rho_{i} \sigma_{i} \rho_{i+1}^{-1} \rho_{i}^{-1}=  \sigma_{i}^{-1} \rho_{i} \rho_{i+1} \sigma_{i} \rho_{i+1}^{-1} \rho_{i}^{-1}\\ 
&= \sigma_{i}^{-1} \rho_{i} \sigma_{i}^{-1} \rho_{i}  \rho_{i+1}^{-1} \rho_{i}^{-1}= \rho_{i+1} \rho_{i}  \rho_{i+1}^{-1} \rho_{i}^{-1}= \sigma_{i}^{-2} =A_{i,i+1}^{-1},
\end{align*}
and the first part of the statement follows. The second part may be obtained by conjugating the relation of the first part by $a^n$ and by applying \repr{permute}(\ref{it:permsi}) and~(\ref{it:permrhi}).

\item We first prove the result in the case $i=1$, $j=3$ and $k=2$. Note that the equality $\rho_{2} A_{1,3} \rho_{2}^{-1}= \rho_{1}^{-1}  A_{1,2} \rho_{1}  A_{1,3} \rho_{1}^{-1}  A_{1,2}^{-1} \rho_{1}$ is equivalent to showing that $\rho_{1}^{-1}  A_{1,2}^{-1} \rho_{1} \rho_{2}$ commutes with $A_{1,3}$, which is the case. Indeed, by part~(\ref{it:revrhocommb}), $\rho_{1}^{-1}  A_{1,2}^{-1} \rho_{1} \rho_{2}= \rho_{2}  A_{1,2}$, and using \repr{present},~\reqref{defaij} and~\reqref{defaijalt}, we have:
\begin{align*}
\rho_{1}^{-1}  A_{1,2}^{-1} \rho_{1} \rho_{2} A_{1,3} (\rho_{1}^{-1}  A_{1,2}^{-1} \rho_{1} \rho_{2})^{-1} &= \rho_{2}  A_{1,2} A_{1,3} A_{1,2}^{-1} \rho_{2}^{-1}=\rho_{2} \sigma_{1} \sigma_{2}^{2} \sigma_{1}^{-1} \rho_{2}^{-1}= \sigma_{1}^{-1} \rho_{1} \sigma_{2}^{2} \rho_{1}^{-1} \sigma_{1}\\
& = \sigma_{1}^{-1} \sigma_{2}^{2} \sigma_{1}= A_{1,3}.
\end{align*}
This proves the result in this case. Now let $1\leq i < k< j\leq n$.
Conjugating the equality $\rho_{2} A_{1,3} \rho_{2}^{-1}= \rho_{1}^{-1}  A_{1,2} \rho_{1}  A_{1,3} \rho_{1}^{-1}  A_{1,2}^{-1} \rho_{1}$ by $a^{-(k-2)}$ and using \repr{permute}(\ref{it:permsi}) and~(\ref{it:permrhi}) yields $\rho_{k} A_{k-1,k+1} \rho_{k}^{-1}= \rho_{k-1}^{-1}  A_{k-1,k} \rho_{k-1}  A_{k-1,k+1} \rho_{k-1}^{-1}  A_{k-1,k}^{-1} \rho_{k-1}$. Conjugating this by $\sigma_{j-1} \sigma_{j-2} \cdots \sigma_{k+1}$ and using \repr{present} and~\reqref{defaij}, we see that:
\begin{equation*}
\rho_{k} A_{k-1,j} \rho_{k}^{-1}= \rho_{k-1}^{-1}  A_{k-1,k} \rho_{k-1}  A_{k-1,j} \rho_{k-1}^{-1}  A_{k-1,k}^{-1} \rho_{k-1}.
\end{equation*}
Setting $\lambda=\sigma_{i}^{-1} \sigma_{i+1}^{-1}  \cdots \sigma_{k-2}^{-1}$ and $\widehat{\lambda}= \sigma_{k-2}^{-1} \sigma_{k-3}^{-1} \cdots \sigma_{i}^{-1}$, it follows from this and~\reqref{defaijalt} that: 
\begin{align*}
\rho_{k} A_{i,j} \rho_{k}^{-1}& \!=\! \lambda \rho_{k} A_{k-1,j} \rho_{k}^{-1}\lambda^{-1}= \lambda \rho_{k-1}^{-1}  A_{k-1,k} \rho_{k-1}  A_{k-1,j} \rho_{k-1}^{-1}  A_{k-1,k}^{-1} \rho_{k-1} \lambda^{-1}\\
& \!=\! \lambda \rho_{k-1}^{-1} \widehat{\lambda} \ldotp \widehat{\lambda}^{-1} A_{k-1,k} \widehat{\lambda} \ldotp \widehat{\lambda}^{-1} \rho_{k-1}  \lambda^{-1} \!\ldotp \lambda A_{k-1,j} \lambda^{-1} \! \ldotp \lambda \rho_{k-1} ^{-1} \widehat{\lambda} \ldotp \widehat{\lambda}^{-1} A_{k-1,k}^{-1} \widehat{\lambda} \ldotp \widehat{\lambda}^{-1} \rho_{k-1} \lambda^{-1}\\
&\!=\! \rho_{i}^{-1} \widehat{\lambda}^{-1} A_{k-1,k} \widehat{\lambda} \rho_{i} A_{i,j} \rho_{i}^{-1} (\widehat{\lambda}^{-1} A_{k-1,k} \widehat{\lambda})^{-1} \rho_{i}.
\end{align*}
The result then follows from the fact that:
\begin{equation*}
\widehat{\lambda}^{-1} A_{k-1,k} \widehat{\lambda}= \sigma_{i} \sigma_{i+1} \cdots \sigma_{k-2} \sigma_{k-1}^{2} \sigma_{k-2}^{-1} \cdots \sigma_{i+1}^{-1} \sigma_{i}^{-1}= \left(\prod_{l=i+1}^{k} A_{i,l} \right)\left( \prod_{l=i+1}^{k-1} A_{i,l} \right)^{-1}=C_{i,k}.\qedhere
\end{equation*}
\end{enumerate}
\end{proof}

\section{Conjugacy classes of finite subgroups of $B_n(\rp)$ and $\mcg{\rp}$}\label{sec:conjclass}


%
%
%


In this section, we prove \repr{BraMap}, \reth{ccrp2} and \reco{ccorder4}. One of the main aims is to determine (as far as possible) the number of conjugacy classes of finite cyclic subgroups of $B_{n}(\rp)$. We start by proving \repr{BraMap}. Given a finite group $G$, we denote its order by $\ord{G}$.




\begin{proof}[Proof of \repr{BraMap}] 
Let $n\geq 2$, and let $\map{\phi}{B_n(\rp)}[\mcg{\rp}]$ be the homomorphism defined in \req{bnrpmcg}. Recall that $\phi$ may be viewed as the canonical projection of $B_{n}(\rp)$ onto the quotient by its centre of order $2$ generated by $\ft$, which is the unique element of $B_{n}(\rp)$ of order $2$. So if $L$ is a finite subgroup of $B_{n}(\rp)$ then $L\cap \ang{\ft}$ is the trivial subgroup if $\ord{H}$ is odd, and is the subgroup $\ang{\ft}$ of $B_{n}$ if $\ord{H}$ is even.

\begin{enumerate}[(a)]
\item From the above considerations, if $H$ is a finite subgroup of $\mcg{\rp}$ then $\phi^{\ast}(H)$ is a subgroup of $B_{n}(\rp)$ that contains $\ft$, so $\ord{\phi^{\ast}(H)}= 2\ord{H}$, and $\phi(\phi^{\ast}(H))=H$. Conversely, if $K$ is a subgroup of even order then $\phi(K)$ is a subgroup of $\mcg{\rp}$ of order $\ord{K}/2$, and $\phi^{\ast}(\phi(K))$ is a subgroup of $B_{n}(\rp)$ of order $\ord{K}$ that contains $K$, so $\phi^{\ast}(\phi(K))=K$. It follows that the restriction described in the statement of part~(\ref{it:BraMapa}) is a bijection, and it is straightforward to see that it induces a bijection of the corresponding conjugacy classes.
\item Let $K=\langle x \rangle$ be a cyclic subgroup of $\mcg{\rp}$ of odd order $d$. Since $\ker{\phi}$ is of order $2$, $\phi^{\ast}(K)$ is of order $2d$, and $\phi^{-1}(x)$ has two elements, one of which, $x_{1}$, is of order $d$, and the other, $x_{2}$, is of the form $\ft x_{1}$. Since $d$ is odd and $\ft$ is central and of order $2$, it follows that $\ang{x_{2}}$ is cyclic and of order $2d$, and hence $\phi^{\ast}(K)=\ang{x_{2}}$. Now $\ang{x_{1}}\cap \ang{\ft}=\brak{1}$ (resp.\ $\ang{x_{2}}\cap \ang{\ft}= \ang{\ft}$), and the correspondence $K \longmapsto \ang{x_{1}}$ (resp.\ $K \longmapsto \phi^{\ast}(K)$) is a bijection between the finite cyclic subgroups of $\mcg{\rp}$ of odd order and the finite cyclic subgroups of $B_n(\rp)$ of odd order (resp.\ the finite cyclic subgroups of $B_n(\rp)$ of even order not divisible by $4$), and it preserves conjugacy classes. The result then follows from \rerems{finitemcgbnrp}(\ref{it:finitemcgbnrpb}).\qedhere
\end{enumerate}  
\end{proof}

\begin{rem}
Let $m\in \N$, and let $H$ (resp.\ $K$) be a cyclic subgroup of $B_{n}(\rp)$ (resp.\ of $\mcg{\rp}$) of order $m$. In a vein similar to that of \rerems{finitemcgbnrp}(\ref{it:finitemcgbnrpb}), using the proof of \repr{BraMap}, we have $\phi^{-1}(K) \cong \Z_{2m}$, $\phi(H) \cong \Z_m$ if $m$ is odd, and $\phi(H) \cong \Z_{m/2}$ if $m$ is even.
\end{rem}


Before coming to the proof of \reth{ccrp2}, we prove the following two lemmas.


\begin{lem}\label{lem:conjinv}
Let $n\geq 2$, let $x\in B_{n}(\rp)$ be an element of finite order, and let $r,s\in \Z$. Then $x^r$ and $x^s$ are conjugate in $B_{n}(\rp)$ if and only if $x^r\in \brak{x^{s}, x^{-s}}$. 
\end{lem}

\begin{proof} 
We first prove that every element of $B_{n}(\rp)$ of finite order is conjugate to its inverse. Applying Murasugi's classification and \req{muragen}, it suffices to prove the result for the elements $\delta_{r,j}$, where $j\in \brak{0,1}$, $0\leq r \leq n-j$ and  $l=\gcd(n-j-r,2r)$. Suppose first that $j=0$, so $0\leq j\leq n$, and let $y=a^{n-r}\garside$. Then:
\begin{align}
y (\rho_{r}\sigma_{r-1}\cdots \sigma_{1}) y^{-1}&= a^{n-r}\garside (\rho_{r}\sigma_{r-1}\cdots \sigma_{1}) \garside^{-1}a^{-(n-r)}\notag\\
&= a^{n-r} \rho_{n+1-r}^{-1} \sigma_{n-r+1}\cdots \sigma_{n-1} a^{-(n-r)}\; \text{by~\reqref{conjsigma} and~\reqref{conjrho}}\notag\\
&= a^{-r} \rho_{n+1-r} \sigma_{n-r+1}^{-1}\cdots \sigma_{n-1}^{-1} a^{r}=\rho_{1}^{-1} \sigma_{1}\cdots \sigma_{r-1}\;\text{by \repr{permute}}\notag\\
&= \sigma_{1}^{-1}\cdots \sigma_{r-1}^{-1}\rho_{r}^{-1}= (\rho_{r}\sigma_{r-1}\cdots \sigma_{1})^{-1} \;\text{by \reqref{defab},}\label{eq:ytype1}
\end{align}
and 
\begin{align}
y (\sigma_{r+1}\cdots \sigma_{n-1}) y^{-1}&= a^{n-r}\garside \sigma_{r+1}\cdots \sigma_{n-1} \garside^{-1}a^{-(n-r)}\notag\\
&= a^{n-r} \sigma_{n-r-1}\cdots \sigma_{1} a^{-(n-r)}\; \text{by~\reqref{conjsigma}}\notag\\
&= a^{-r} \sigma_{n-r-1}^{-1}\cdots \sigma_{1}^{-1} a^{r}= \sigma_{n-1}^{-1}\cdots \sigma_{r+1}^{-1} \quad\text{by \repr{permute}}\\
&= (\sigma_{r+1}\cdots \sigma_{n-1})^{-1}.\label{eq:ytype1a}
\end{align}
Since the two elements $\rho_{r}\sigma_{r-1}\cdots \sigma_{1}$ and $\sigma_{r+1}\cdots \sigma_{n-1}$ commute, it follows from~\reqref{muragen},~\reqref{ytype1} and~\reqref{ytype1a} that $y\delta_{r,0} y^{-1}=\delta_{r,0}^{-1}$. 

Now suppose now that $j=1$, so $0\leq r\leq n-1$. If $r=n-1$ then $\delta_{n-1,1}=b$, and taking $y=\garside a^{-1}$ and using \repr{permute}(\ref{it:embed}), we have $\ang{\delta_{n-1,1},y}\cong \dic{4(n-1)}$, in particular, $y\delta_{n-1,1} y^{-1}=\delta_{n-1,1}^{-1}$. If $r=n-2$ then $\delta_{n-2,1}$ is of order $2$, so $\delta_{n-2,1}=\ft$, and the result is clear. So suppose that $0\leq r\leq n-3$, and let $y=\sigma_{r+2}\sigma_{r+3}\cdots \sigma_{n-2}a^{n-r}\garside$. Using~\reqref{ytype1} and the fact that $\rho_{r}\sigma_{r-1}\cdots \sigma_{1}$ and $\sigma_{r+2}\cdots \sigma_{n-2}$ commute, we see that:
\begin{equation}\label{eq:ycomm2}
y (\rho_{r}\sigma_{r-1}\cdots \sigma_{1}) y^{-1}= (\rho_{r}\sigma_{r-1}\cdots \sigma_{1})^{-1}.
\end{equation}
As in \req{ytype1a}, we obtain:
\begin{equation}\label{eq:yconjsigma}
y z y^{-1}= (\sigma_{r+2}\cdots \sigma_{n-2}) \sigma_{n-1}^{-1}\cdots \sigma_{r+1}^{-1} \sigma_{n-1}^{-1}(\sigma_{r+2}\cdots \sigma_{n-2})^{-1}.
\end{equation}
where $z=\sigma_{r+1}\cdots \sigma_{n-1}\sigma_{r+1}$. If $r=n-3$ then~\reqref{yconjsigma} becomes:
\begin{equation*}
y (\sigma_{n-2} \sigma_{n-1}\sigma_{n-2}) y^{-1}= \sigma_{n-1}^{-1} \sigma_{n-2}^{-1} \sigma_{n-1}^{-1}=(\sigma_{n-2} \sigma_{n-1} \sigma_{n-2})^{-1},
\end{equation*}
while if $r\leq n-2$ then:
\begin{align*}
yzy^{-1}&=\bigl( (\sigma_{r+2}\cdots \sigma_{n-2}) \sigma_{n-1}\sigma_{r+1} \cdots \sigma_{n-1}(\sigma_{r+2}\cdots \sigma_{n-2})^{-1}\bigr)^{-1}\\
&=  \bigl( (\sigma_{r+2}\cdots \sigma_{n-2}) \sigma_{n-1}\sigma_{r+1} \cdots \sigma_{n-1}(\sigma_{n-2}^{-1}\cdots \sigma_{r+2}^{-1})\sigma_{n-1}^{-1}\cdots \sigma_{r+1}^{-1} \ldotp \sigma_{r+1} \cdots \sigma_{n-1}\bigr)^{-1}\\
&= \bigl( \sigma_{r+2}\cdots \sigma_{n-2} \sigma_{n-1} \sigma_{n-1}^{-1}\cdots \sigma_{r+3}^{-1} \sigma_{r+1} \cdots \sigma_{n-1}\bigr)^{-1}= \bigl( \sigma_{r+2}\sigma_{r+1} \sigma_{r+2} \sigma_{r+3}\cdots \sigma_{n-1}\bigr)^{-1}\\
&= \bigl(\sigma_{r+1}\sigma_{r+2} \sigma_{r+1} \sigma_{r+3}\cdots \sigma_{n-1}\bigr)^{-1}= z^{-1}.
\end{align*}
So for all $0\leq r\leq n-3$, $yz y^{-1}=z^{-1}$ by \req{yconjsigma}. Using the fact that $\rho_{r}\sigma_{r-1}\cdots \sigma_{1}$ and $\sigma_{r+1}\cdots \sigma_{n-1}\sigma_{r+1}$ commute as well as equations~\reqref{muragen} and~\reqref{ycomm2}, we obtain $y \delta_{r,1}y^{-1}=\delta_{r,1}^{-1}$. We conclude that every finite order element of $B_n(\rp)$ is conjugate to its inverse.

Conversely, let $x\in B_{n}(\rp)$ be an element of finite order, and let $r,s\in \Z$ be such that $x^r$ and $x^s$ are conjugate in $B_{n}(\rp)$. Let $u\in B_{2n}(\St)$ be the image of $x$ under the embedding of $B_{n}(\rp)$ in $B_{2n}(\St)$ given by~\cite[Corollary~2]{GG5}. Then $u^r$ and $u^s$ are conjugate in $B_{2n}(\St)$. Using Murasugi's classification of the finite order elements of the braid groups of the sphere~\cite[Theorem~A]{M}, it follows from~\cite[Proposition~9]{GG8} that $u^{r}\in \brak{u^{s}, u^{-s}}$, and hence $x^{r} \in \brak{x^{s}, x^{-s}}$ as required. 
\end{proof}

\begin{lem}\label{lem:rhorrho1}
Let $1\leq r\leq n$.
In $B_{n}(\rp)$, we have the following equality 
$\bigl( \rho_{r}\sigma_{r-1}\cdots \sigma_{1}\bigr)^{r}=\rho_{r}\cdots \rho_{1}$, so $\bigl( \rho_{r}\sigma_{r-1}\cdots \sigma_{1}\bigr)^{r}\in P_{n}(\rp)$.
\end{lem}

\begin{proof}
The proof is by induction on $r$. If $r=1$, there is nothing to prove. So suppose that the result holds for some $1\leq r\leq n-1$. Thus $(\rho_{r}\cdots \sigma_{1})^{r}=\rho_{r}\cdots \rho_{1}$, and conjugating this relation by $a^{-1}$ and using \repr{permute}(\ref{it:permsi}) and~(\ref{it:permrhi}), we obtain $(\rho_{r+1} \cdots \sigma_{2})^{r}=\rho_{r+1}\cdots \rho_{2}$. Multiplying each side of this equation on the right by the corresponding side of the relation $\sigma_1 \cdots \sigma_r \rho_{r+1} \sigma_r \cdots \sigma_1= \rho_1$, which may be obtained using the fact that $\rho_i= \sigma_i \rho_{i+1} \sigma_i$ for all $i=1,\ldots, n-1$ by~\reqref{sirisi}, we see that:
\begin{equation}\label{eq:rhr1r1a}
\rho_{r+1}\cdots \rho_{1}= (\rho_{r+1} \cdots \sigma_{2})^{r} \sigma_1 \cdots \sigma_r \rho_{r+1} \sigma_r \cdots \sigma_1
\end{equation}
Let us prove by induction on $0\leq j\leq r$ that: \begin{equation}\label{eq:rhr1r1b}
\rho_{r+1}\cdots \rho_{1}= (\rho_{r+1} \cdots \sigma_{2})^{r-j} (\rho_{r+1} \cdots \sigma_{1})^{j} \sigma_{j+1} \cdots \sigma_r (\rho_{r+1} \sigma_r \cdots \sigma_1).
\end{equation}
If $j=0$ then the equality~\reqref{rhr1r1b} is just~\reqref{rhr1r1a}. So suppose that~\reqref{rhr1r1b} is valid for some $0\leq j\leq r-1$. Then by  relations~\reqref{artin1}--\reqref{sirj}, for all $2\leq k\leq r$, we have $\rho_{r+1} \sigma_r \cdots \sigma_1 \ldotp \sigma_k= \sigma_{k-1} \ldotp \rho_{r+1} \sigma_r \cdots \sigma_1$, and thus by~\reqref{rhr1r1b}, we have:
\begin{align*}
\rho_{r+1}\cdots \rho_{1} &= (\rho_{r+1} \cdots \sigma_{2})^{r-j} (\rho_{r+1} \cdots \sigma_{1})^{j} \sigma_{j+1} \cdots \sigma_r (\rho_{r+1} \sigma_r \cdots \sigma_1)\\
&= (\rho_{r+1} \cdots \sigma_{2})^{r-j} \sigma_{1}(\rho_{r+1} \cdots \sigma_{1})^{j} \sigma_{j+2} \cdots \sigma_r (\rho_{r+1} \sigma_r \cdots \sigma_1)\\
&= (\rho_{r+1} \cdots \sigma_{2})^{r-j-1} (\rho_{r+1} \cdots \sigma_{1})^{j+1} \sigma_{j+2} \cdots \sigma_r (\rho_{r+1} \sigma_r \cdots \sigma_1),
\end{align*}
which yields~\reqref{rhr1r1b} for $j+1$, and the result is obtained by taking $j=r$.
\end{proof}
 
\begin{proof}[Proof of \reth{ccrp2}.]\mbox{}
%
Let $n\geq 2$.
\begin{enumerate}[(a)]
\item Suppose that $d$ is divisible by $4$, and let $x,y$ be elements of $B_n(\rp)$ of order $d$ that are powers of elements of the form~\reqref{muragen}, and that are conjugate in $B_n(\rp)$. We claim that they are powers of the same $\delta_{r,j}$, where $j\in \brak{0,1}$ and $0\leq r \leq n-j$, from which it follows that $y\in \brak{x,x^{-1}}$ by \relem{conjinv}, which proves the result. To prove the claim, suppose on the contrary that $x$ and $y$ are powers of distinct $\delta_{r,j}$. Since they are conjugate, their $d/4$th powers are of order $4$ and are conjugate.
 
We first identify the Murasugi elements of order $4$. If $j\in \brak{0,1}$, and $0\leq r \leq n-j$, recall that $\delta_{r,j}$ is of order $2l$, where $l=\gcd(n-j-r,2r)$. Note that:
\begin{equation}\label{eq:4divides}
\text{$4$ divides the order of $\delta_{r,j}$ $\Longleftrightarrow$ $2 \mid l$ $\Longleftrightarrow$ $2 \mid n-j-r$ $\Longleftrightarrow r+j\equiv n\bmod{2}$.}
\end{equation}
This being the case, we have $l=2l'$, where $l'=\gcd((n-j-r)/2,r)$, and $\delta_{r,j}$ is of order $4l'$. Setting:
\begin{equation}\label{eq:muragen2}
\eta_{r,j}=\delta_{r,j}^{l'}=\bigl( \rho_{r}\sigma_{r-1}\cdots \sigma_{1}\bigr)^{r} \bigl( \sigma_{r+1}\cdots \sigma_{n-1}\sigma_{r+1}^{j}\bigr)^{(n-j-r)/2},
\end{equation}
it follows that $\eta_{r,j}$ and $\eta_{r,j}^{3}$ are the powers of $\delta_{r,j}$ that are of order $4$. Since these two elements are conjugate by \relem{conjinv}, it suffices to 
show that the $\eta_{r,j}$ are pairwise non conjugate. More precisely, if $0\leq r\leq n-j$,  $0\leq r'\leq n-j'$, where $j,j'\in \brak{0,1}$, and $r+j \equiv r'+j' \equiv n \bmod{2}$, then $\eta_{r,j}$ and $\eta_{r',j'}$ are conjugate if and only if $(r,j)=(r',j')$. Clearly the condition is sufficient. Conversely, assume that $(r,j)\neq (r',j')$. First suppose that $j\neq j'$. Exchanging the rôles of $(r,j)$ and $(r',j')$, we may suppose that $j=0$ and $j'=1$. 
Then $\eta_{r,0}$ is not conjugate to $\eta_{r',1}$, since under the Abelianisation homomorphism~\reqref{abmap}, the first factor of $q(\eta_{r,0})$ (that corresponding to the $\rho_j$) is equal to $\overline{r}=\overline{n} \bmod{2}$ by~\reqref{4divides}, while that of $q(\eta_{r',1})$ is equal to $\overline{r'}= \overline{n-1} \bmod{2}$, from which we conclude that $q(\eta_{r,0})\neq q(\eta_{r',1})$, so $\eta_{r,j}$ and $\eta_{r',j'}$ are non-conjugate. Thus it remains to show that if $j\in \brak{0,1}$ and $0\leq r,r'\leq n-j$, with $r\neq r'$ and $r,r'\equiv n-j\bmod{2}$, then $\eta_{r,j}$ is not conjugate to $\eta_{r',j}$. We shall do this by determining their cycle types. By \relem{rhorrho1}, we have $\bigl( \rho_{r}\sigma_{r-1}\cdots \sigma_{1} \bigr)^{r}= \rho_{r}\cdots \rho_{1}\in P_{n}(\rp)$, and the first factor on the right-hand side term~\reqref{muragen2} is a pure braid. As for the second factor of this term, its permutation is $\bigl(n, \frac{n+r}{2}\bigr) \bigl(n-1, \frac{n+r-2}{2}\bigr)\cdots \bigl(\frac{n+r+2}{2}, r+1\bigr)$ (resp.\ $\bigl(n, \frac{n+r+1}{2}\bigr) \bigl(n-1, \frac{n+r-1}{2}\bigr)\cdots \bigl(\frac{n+r+5}{2}, r+3\bigr)\bigl(\frac{n+r+3}{2}, r+1\bigr)$) if $j=0$ (resp.\ if $j=1$), which is the product of $\frac{n-r}{2}$ (resp.\ $\frac{n-1-r}{2}$) disjoint transpositions. Since $r\neq r'$, it follows that the permutations of $\eta_{r,j}$ and $\eta_{r',j}$ have distinct cycle types, and hence $\eta_{r,j}$ and $\eta_{r',j}$ are non conjugate, which proves the result.


\item 
Suppose that $d$ is not divisible by $4$, and assume for the moment that $n\geq 2$. If $d\in \brak{1,2}$ then $B_{n}(\rp)$ has a unique central element of order $d$ and we are done. So assume that $d\geq 3$. Then $d$ is either odd, or is equal to twice an odd number. We claim that it suffices to prove the result for the second possibility. To see this, let $d$ be odd, and assume that any two elements of order $2d$ appearing in Murasugi's classification are conjugate. Let $x,x'\in B_{n}(\rp)$ be two elements of order $d$ of the form~\reqref{muragen}. Since $d$ is odd and $\ft$ is central and of order $2$ by \repr{agt}(\ref{it:agtb}), $x\ft$ and $x'\ft$ are of order $2d$, so are conjugate, and since $\ft$ is central, this implies that $x$ and $x'$ are conjugate, which proves the claim.

Now let $n\in \brak{3,4}$ (resp.\ $n=5$). By \repr{agt}(\ref{it:agta}) the only non-trivial elements of odd order in $B_n(\rp)$ are of order $3$ (resp.\ of order $5$). By the previous paragraph, to show that all such elements are pairwise conjugate in $B_n(\rp)$, it suffices to prove that the Murasugi elements of order $6$ (resp.\ of order $10$) are pairwise conjugate. Using~\reqref{muragen} and~\relem{conjinv}, up to conjugacy, there are two such elements for each value of $n$: $(\rho_3 \sigma_2 \sigma_1)^2$ and $\sigma_1 \sigma_2$ if $n=3$; $(\rho_3 \sigma_2 \sigma_1)^2$ and $\sigma_1 \sigma_2 \sigma_3 \sigma_1$ if $n=4$; $\sigma_1 \sigma_2 \sigma_3 \sigma_4$ and $(\rho_5 \sigma_4 \sigma_3 \sigma_2 \sigma_1)^2$ if $n=5$. First suppose that $n=3,4$. Let $\lambda=\rho_{3} \rho_{1}$. Then by \repr{permute}(\ref{it:permrhi}) and~\reqref{sirisi}, we have:
\begin{align*}
\lambda^{-1} (\rho_3 \sigma_2 \sigma_1)^2 \lambda&= \rho_{1}^{-1}  \sigma_2 \sigma_1 \rho_3 \sigma_2 \sigma_1 \rho_{3} \rho_{1}= \rho_{1}^{-1}  \sigma_2 \sigma_1 \rho_{2} \rho_3 \sigma_2 \sigma_1  \rho_{1}= \rho_{1}^{-1}  \sigma_2 \sigma_1 \rho_{2} \sigma_1 \sigma_1^{-1}\rho_3 \sigma_2 \sigma_1  \rho_{1}\\
& = \rho_{1}^{-1}  \sigma_2 \rho_{1} \sigma_1^{-1}\rho_3 \sigma_2 \sigma_1  \rho_{1}=  \sigma_2 \sigma_1^{-1} \sigma_2^{-1} \sigma_1^{-1} \sigma_1 \sigma_2 \rho_3 \sigma_2 \sigma_1  \rho_{1}= \sigma_2 \sigma_1^{-1} \sigma_2^{-1} \sigma_1^{-1} \rho_{1}^{2}. 
\end{align*}
Let $\sim$ mean `is conjugate to'. Then by~\reqref{surfacerp2} and the Artin relations, if $n=3$ we see that:
\begin{align*}
\lambda^{-1} (\rho_3 \sigma_2 \sigma_1)^2 \lambda&= \sigma_2 \sigma_1^{-1} \sigma_2^{-1} \sigma_1^{-1} \sigma_{1} \sigma_{2}^{2} \sigma_{1} =\sigma_2 \sigma_1^{-1} \sigma_{2} \sigma_{1} \sim \sigma_{1}  \sigma_2 \sigma_1^{-1} \sigma_{2} = \sigma_1^{-1} \sigma_{1}  \sigma_2^{2}\sim \sigma_{1}  \sigma_2,
\end{align*}
while if $n=4$, 
\begin{align*}
\lambda^{-1} (\rho_3 \sigma_2 \sigma_1)^2 \lambda&= \sigma_2 \sigma_1^{-1} \sigma_2^{-1} \sigma_1^{-1} \sigma_{1} \sigma_{2} \sigma_{3}^{2} \sigma_{2} \sigma_{1}= \sigma_2 \sigma_1^{-1} \sigma_{3}^{2} \sigma_{2} \sigma_{1}= \sigma_2  \sigma_{3}^{2} \sigma_{2} \sigma_{1} \sigma_2^{-1} \sim \sigma_{3}^{2} \sigma_{2} \sigma_{1}\\
& \sim \sigma_{3}\sigma_{1} \sigma_{2} \sigma_{3}= \sigma_{1} \sigma_{2} \sigma_{3} \sigma_{2} \sim \sigma_{2}  \sigma_{1} \sigma_{2} \sigma_{3} = \sigma_{1} \sigma_{2} \sigma_{3} \sigma_{1}.  
\end{align*}
This proves the result in the cases $n=3,4$. Finally suppose that $n=5$, and let $\lambda=\rho_{5} \rho_{3} \rho_{1}$. Then by \repr{permute}(\ref{it:permrhi}), equations~\reqref{sirisi} and~\reqref{surfacerp2} and the Artin relations, we have:
\begin{align*}
\lambda^{-1} a^{2} \lambda&= \rho_{1}^{-1} \rho_{3}^{-1} \sigma_4 \sigma_3 \sigma_2 \sigma_1 a \rho_{5} \rho_{3} \rho_{1}= \rho_{1}^{-1} \rho_{3}^{-1} \sigma_4 \sigma_3 \sigma_2 \sigma_1 \rho_{4} \rho_{2} a  \rho_{1}\\
&= \rho_{1}^{-1} \rho_{3}^{-1} \sigma_4 \sigma_3 \rho_{4} \sigma_3 \sigma_3^{-1} \sigma_2 \sigma_1  \rho_{2} \sigma_1 \sigma_1^{-1} a  \rho_{1}= \rho_{1}^{-1} \rho_{3}^{-1} \sigma_4 \rho_{3} \sigma_3^{-1} \sigma_2 \rho_{1}  \sigma_1^{-1} a  \rho_{1}\\
& = \sigma_4  \sigma_3^{-1} \sigma_2  \sigma_1^{-1} a  \rho_{1}= \sigma_4  \sigma_3^{-1} \sigma_2  \sigma_1^{-1} \sigma_4^{-1} \sigma_3^{-1} \sigma_2^{-1} \sigma_1^{-1} \rho_{1}^{2}\\
&= \sigma_4  \sigma_3^{-1} \sigma_2  \sigma_1^{-1} \sigma_4^{-1} \sigma_3^{-1} \sigma_2^{-1} \sigma_1^{-1} \sigma_1 \sigma_2 \sigma_3 \sigma_4^{2}\sigma_3 \sigma_2 \sigma_1 = \sigma_4  \sigma_3^{-1} \sigma_2  \sigma_1^{-1} \sigma_4\sigma_3 \sigma_2 \sigma_1\\
&= \sigma_4  \sigma_3^{-1} \sigma_2  \sigma_4\sigma_3 \sigma_2 \sigma_1 \sigma_2^{-1} \sim  \sigma_4  \sigma_3 \sigma_2^{-1} \sigma_3^{-1}  \sigma_4\sigma_3 \sigma_2 \sigma_1 = \sigma_4  \sigma_3 \sigma_2^{-1} \sigma_4\sigma_3 \sigma_2 \sigma_1 \sigma_4^{-1}\\
& \sim \sigma_3 \sigma_2^{-1} \sigma_4\sigma_3 \sigma_2 \sigma_1= \sigma_3 \sigma_4\sigma_3 \sigma_2 \sigma_1 \sigma_3^{-1} \sim \sigma_1\sigma_2 \sigma_3 \sigma_4,
\end{align*}
which proves the result in the case $n=5$.\qedhere
\end{enumerate}
\end{proof}

We are now able to prove \reco{ccorder4}.

\begin{proof}[Proof of \reco{ccorder4}]
Let $n\geq 3$. First suppose that $d=4$. Making use of the proof of \reth{ccrp2}(\ref{it:ccrp2a}), \relem{conjinv} and~\req{4divides}, in $B_n(\rp)$, that there are $2(n+1)$ Murasugi elements of order $4$, and so there are $n+1$ conjugacy classes of cyclic subgroups of order $4$. More precisely, if $n$ is even (resp.\ odd), then there are $n+2$ such elements if $j=0$ and $n$ such elements if $j=1$ (resp.\ $n+1$ such elements if $j=0$ or if $j=1$). This proves the result in this case.

Now suppose that $d>4$, where $d$ is divisible by $4$, so $d\geq 8$. By \repr{agt}(\ref{it:agta}), there exists a unique $j\in \brak{0,1}$ such that $d$ divides $4(n-j)$, and so up to conjugacy, the Murasugi elements of order $d$ are powers of certain $\delta_{r,j}$, where $0\leq r\leq n-j$. Let us determine the values of $r$ for which this is the case. Since $\delta_{r,j}$ is of order $2\gcd(n-j-r,2r)$, $d/2$ must divide $n-j-r$ and $2r$. Since $d/2$ divides $2(n-j)$, for $d/2$ to divide $2r$, we must have $r=n-j-\frac{kd}{4}$, where $0\leq k\leq 4(n-j)/d$, so $n-j-r=kd/4$. But $d/2$ must also divide $n-j-r$, which implies that $k$ is even. One may check that the number of even values of $k$ for which $0\leq k\leq 4(n-j)/d$ is equal to $\lfloor 2(n-j)/d \rfloor$. It follows from the proof of \reth{ccrp2}(\ref{it:ccrp2a}) that this number is equal to the number of conjugacy classes of cyclic subgroups of $B_{n}(\rp)$ of order $d$.
\end{proof}


Recall from~\cite[Corollary~19]{GG2} that with the exception of the two elements of $\ang{\ft}$, the torsion elements of $P_{n}(\rp)$ are of order $4$. We have the following result regarding their conjugacy classes in $B_{n}(\rp)$.

\begin{prop}\label{prop:ord4Pn}
Let $n\geq 2$. Then every element of order $4$ of $P_{n}(\rp)$ is conjugate in $B_{n}(\rp)$ to one of $a^{n}$ or $b^{n-1}$. 
\end{prop}

\begin{proof}
The result follows from~\cite[Proposition~10]{GG3} and \relem{conjinv}.
\end{proof}

If $M$ is a compact, connected surface without boundary different from $\St$ and $\rp$, the inclusion of a topological disc of $M$ in the surface $M$ induces an embedding of the Artin braid group $B_{n}$ in the $n$-string braid group $B_{n}(M)$ (see~\cite[Lemma~3]{Gold} for the restriction to the pure braid subgroups and~\cite[Proposition~12]{GGcambridge}). 
In the case $M=\St$, the kernel of this embedding is isomorphic to the direct product of a free group of rank $n-1$ with the infinite cyclic subgroup generated by the square $\garside^{4}$ of the full-twist braid~\cite[Theorem~4]{DL}.

\begin{prop}\label{prop:inclus}
Let $n \geq 2$, and consider the homomorphism $\map{\iota_{n}}{B_{n}}[B_{n}(\rp)]$ induced by the inclusion $\iota$ of a topological disc in $\rp$. Then $\ker{\iota_n}=\ang{\garside^{4}}$.
\end{prop}

\begin{proof}
Clearly $\ang{\garside^{4}} \subset \ker{\iota_{n}}$ by \repr{agt}(\ref{it:agtb}). Conversely, let $x \in \ker{\iota_{n}}$. Since the permutation of any element $y\in B_{n}$ is equal to that of $\iota_{n}(y)$, it follows that $x\in P_{n}$. It thus suffices to consider the restriction $\map{\iota_{n}\left\lvert_{P_{n}}\right.}{P_{n}}[P_{n}(\rp)]$ of $\iota_{n}$ to $P_{n}$. If $n=2$ then $P_{2}=\ang{\ft}$, so $\im{\iota_{n}}=\ang{\ft}$, and the result follows. So suppose by induction that the result holds for some $n\geq 2$. We have the following commutative diagram of short exact sequences:
\begin{equation*}
\begin{tikzcd}
& & 1 \arrow{d}& 1\arrow{d} &\\
& & \ker{\iota_{n+1}} \arrow{d} \arrow[dashrightarrow]{r}{p\left\lvert_{\ker{\iota_{n+1}}}\right.} & \ker{\iota_{n}} \arrow{d} & \\
1 \arrow{r} & \pi_{1}(\dt \setminus \brak{x_{1}, \ldots, x_{n}}, x_{n+1}) \arrow{r} \arrow{d}{f} & P_{n+1} \arrow{d}{\iota_{n+1}} \arrow{r}{p} & P_{n} \arrow{d}{\iota_{n}} \arrow{r} & 1\\
1 \arrow{r} & \pi_{1}(\rp \setminus \brak{x_{1}, \ldots, x_{n}}, x_{n+1}) \arrow{r}  & P_{n+1}(\rp) \arrow{r}{\widetilde{p}} & P_{n}(\rp) \arrow{r}  & 1,
\end{tikzcd}
\end{equation*}
where $p, \widetilde{p}$ denote the surjective homomorphisms defined by forgetting the last string, and $f$ is the homomorphism on the level of fundamental groups induced by $\iota$ 
The fact that the above diagram is commutative implies that $p$ restricts to a homomorphism $\map{p\left\lvert_{\ker{\iota_{n+1}}}\right.}{\ker{\iota_{n+1}}}[\ker{\iota_{n}}]$. Since $\ker{\iota_{n}}=\ang{\garside^{4}}$ by induction, $\ang{\garside[n+1]^{4}} \subset \ker{\iota_{n+1}}$ by \repr{agt}(\ref{it:agtb}) and $p(\ft[n+1])=\ft$, we see that  $p\left\lvert_{\ker{\iota_{n+1}}}\right.$ is surjective. Now $\ker{p}$ (resp.\ $\ker{\widetilde{p}}$) is a free group of rank $n$, and up to identification with the corresponding subgroup of $P_{n+1}$ (resp.\ $P_{n+1}(\rp)$), the set $\brak{A_{i,n+1}}_{i=1,\ldots,n}$ (resp.\ the set $\brak{A_{i,n+1}} _{i=1,\ldots,n-1} \cup \brak{\rho_{n+1}}$) is a basis of $\ker{p}$ (resp.\ $\ker{\widetilde{p}}$), where for $1\leq i<j\leq n+1$, $A_{i,j}$ is defined by~\reqref{defaij}. Further, $\iota_{n+1}(A_{i,n+1})=A_{i,n+1}$ for $i=1,\ldots n-1$, and $\iota_{n+1}(A_{n,n+1})=A_{n,n+1}=\rho_{n+1}^{2} A_{1,n+1}, \ldots, A_{n-1,n+1}$ by~\reqref{rn2} and~\reqref{defaij}, which implies that $\iota_{n+1}(\ker{p})$ is the free group $\langle \rho_{n+1}^2, A_{1,n+1},\ldots ,A_{n-1,n+1}\rangle$ of rank $n$. Since free groups of finite rank are Hopfian, it follows that the restriction of $\iota_{n+1}$ to $\ker{p}$ is injective.
A straightforward diagram-chasing argument then shows that $p\left\lvert_{\ker{\iota_{n+1}}}\right.$ is injective, and we conclude that $\ker{\iota_{n+1}}=\bigl\langle\garside[n+1]^{4}\bigr\rangle$ as required.
\end{proof}

\section{Embeddings of certain finite subgroups of $B_{n}(\rp)$}\label{sec:inclcyclic}   

Having studied the conjugacy classes of the finite cyclic subgroups of $B_{n}(\rp)$ in \resec{conjclass}, we now turn to the analogous problem for dicyclic subgroups. One of the aims of this section is to understand better the inclusions between the two classes of subgroups, and to determine the centralisers of these subgroups. We start by showing that no finite cyclic subgroup of $B_{n}(\rp)$ is maximal as a finite subgroup.

\begin{prop}\label{prop:nomaxcyclic}
Let $n\geq 2$, and let $x$ be an element of $B_{n}(\rp)$ of finite order $m$.  If $m\in \brak{1,2}$ then $x$ belongs to a dicyclic group of order $8n$. If $m\geq 3$ then $x$ belongs to a dicyclic group of order $2m$ (resp.\ of order $4m$) if $m$ is even (resp.\ if $m$ is odd). 
\end{prop}

\begin{proof}
Let $x\in B_{n}(\rp)$ be of order $m$. If $m\in \brak{1,2}$ then $x\in \ang{\ft}$ by~\repr{agt}(\ref{it:agtb}), and so $x$ belongs to the standard copy $\ang{a,\garside}$ of $\dic{8n}$. So suppose that $m\geq 3$. If $m$ is odd then $x\ft$ is of order $2m$ because $\ft$ is central and of order $2$. Thus $\ang{x\ft}$ contains $\ft$, and hence contains $x$. So considering $x\ft$ if necessary, we may suppose that $m$ is even, and thus $m\geq 4$. Conjugating $x$ if necessary, we may suppose by Murasugi's classification that $x$ is a power of some $\delta_{r,j}$ given by~\reqref{muragen}, where $j\in\brak{0,1}$ and $0\leq r\leq n-j$. To prove the result, we consider the elements $y\in B_{n}(\rp)$ that appear in the proof of \relem{conjinv}. Since $y\delta_{r,j}y^{-1}= \delta_{r,j}^{-1}$, we have $yxy^{-1}=x^{-1}$. We analyse the cases of the proof of \relem{conjinv} in turn.
\begin{enumerate}
\item If $j=0$ then $y=a^{n-r}\garside$, and using \req{conja}, we obtain:
\begin{equation}\label{eq:ankdelta}
y^{2}=a^{n-r}\garside a^{n-r}\garside=\ft.
\end{equation}
Since $\ft$ is the unique element of $B_{n}(\rp)$ of order $2$ by \repr{agt}(\ref{it:agtb}), we have $x^{m/2}=y^{2}$. Hence $\ang{x,y}$ is isomorphic to a quotient of $\dic{2m}$. But $x$ is of order $m\geq 4$, so $x\neq x^{-1}$, which implies that $y\notin \ang{x}$, and hence $\ang{x,y}\cong \dic{2m}$ (if the original $m$ is odd then we replace $x$ by $x\ft$, and thus $m$ is replaced by $2m$).

\item If $j=1$ and $r=n-1$ then as we saw, $x$ belongs to the subgroup $\ang{\delta_{n-1,1}, y}$ which is isomorphic to $\dic{4(n-1)}$.

\item If $j=1$ and $r=n-2$ then $\delta_{n-2,1}$ is of order $2$, which contradicts the fact that $m\geq 4$.

\item If $j=1$ and $0\leq r\leq n-3$ then $y=\sigma_{r+2}\sigma_{r+3}\cdots \sigma_{n-2}a^{n-r}\garside$, and: 
\begin{align*}
y^{2}&=\sigma_{r+2}\cdots \sigma_{n-2}a^{n-r}\garside \sigma_{r+2}\cdots \sigma_{n-2}a^{n-r}\garside\\
& =\sigma_{r+2}\cdots \sigma_{n-2}(a^{n-r}\garside \sigma_{r+2}\cdots \sigma_{n-2} \garside^{-1}a^{-(n-r)})(a^{n-r}\garside)^2 \\
&= \sigma_{r+2}\cdots \sigma_{n-2}(a^{n-r} \sigma_{n-r-2}\cdots \sigma_{2} a^{-(n-r)}) \ft\; \text{by~\reqref{conjsigma} and~\reqref{ankdelta}}\\
&= \sigma_{r+2}\cdots \sigma_{n-2}(a^{-r} \sigma_{n-r-2}^{-1}\cdots \sigma_{2}^{-1} a^{r}) \ft=\ft \;\text{by \repr{permute},}
\end{align*}
so once more $x^{m/2}=y^{2}$, and as in the case $j=0$, we see that $\ang{x,y}\cong \dic{2m}$ as required.\qedhere
\end{enumerate}
\end{proof}

From \repr{nomaxcyclic}, we obtain the following corollary.

\begin{cor}\label{cor:cyclicmax}
No finite cyclic subgroup of $B_{n}(\rp)$ is finite maximal.
\end{cor}

If $x\in B_{n}(\rp)$ is of finite order, it follows from \reco{cyclicmax} and the proof of \repr{nomaxcyclic} that $x$ and $x^{-1}$ are conjugate by an element of finite order. One may ask whether they are conjugate by an element of infinite order. This question is closely related to the existence of certain virtually-cyclic subgroups of $B_{n}(\rp)$, and we will return to it after \repr{centraliserx}.


The following result is the analogue of~\cite[Proposition~28]{GG8} for the braid groups of $\rp$.

\begin{prop}\label{prop:bnrpsubgp}
Let $n\geq 3$. Any infinite subgroup of $B_{n}(\rp)$ contains an element of infinite order. In particular, any subgroup of $B_{n}(\rp)$ consisting entirely of elements of finite order is itself finite.
\end{prop}

\begin{proof}
If $H$ is a subgroup of $B_{n}(\rp)$ consisting entirely of elements of finite order, then it embeds abstractly in $B_{2n}(\St)$ by~\cite[Corollary~2]{GG5}, and the result follows from~\cite[Proposition~28]{GG8}.
\end{proof}

If $\gamma$ is an element of a group $\Gamma$, let $Z_{\Gamma}(\gamma)$ denote the centraliser of $\gamma$ in $\Gamma$, and if $H$ is a subgroup of $\Gamma$, let $Z_{\Gamma}(H)$ (resp.\ $N_{\Gamma}(H)$) denote the centraliser (resp.\ normaliser) of $H$ in $\Gamma$. The following result generalises~\cite[Proposition~9]{GGzeit} and~\cite[Proposition~8]{GG8}. 

\begin{prop} \label{prop:centraliserx} 
Assume that $n\geq 2$. Let $x\in \brak{a,b}$, let $N=n$ (resp.\ $N=n-1$) if $x=a$ (resp.\ if $x=b$). If $k\in \N$ and $k/\gcd(k,N)$ is odd then:
\begin{enumerate}
\item\label{it:centraliser} $Z_{B_{n}(\rp)}(x^{k})=\ang{x}$.
\item\label{it:normaliser} The normaliser $N_{B_{n}(\rp)}(\ang{x^{k}})$ is equal to the standard copy $\ang{x,a\garside}$ of $\dic{8N}$ in $B_{n}(\rp)$, unless $x=b$ and $n=2$, in which case:
\begin{equation*}
\ang{b,a\garside[2]} \subsetneqq N_{B_{2}(\rp)}(\ang{x^{k}})=B_{2}(\rp).
\end{equation*}
In particular, there is no element $\gamma\in B_{n}(\rp)$ of infinite order for which $\gamma x^k \gamma^{-1}=x^{-k}$.
\end{enumerate}
\end{prop}

\begin{rem}
With the notation of the statement of \repr{centraliserx}, note that if $k$ divides $N$ then the condition that $k/\gcd(k,N)$ is odd is satisfied. 
\end{rem}

\begin{proof}[Proof of \repr{centraliserx}]\mbox{}
\begin{enumerate}
\item Let $k=2^{c}k'$ and $N=2^{d}N'$, where $k',N'\in \N$ are odd, and $c,d\geq 0$, and let $l=\gcd(k,N)$. Since $k/\gcd(k,N)$ is odd, we have $c\leq d$, and so $l= \gcd(k,N)=\gcd(k,4N)$ 
Thus $\ang{x^{l}}= \ang{x^{k}}$ because $x$ is of order $4N$, and $\ang{x}\subset Z_{B_{n}(\rp)}(x^{k})=Z_{B_{n}(\rp)}(x^{l})$. Conversely, let $y\in Z_{B_{n}(\rp)}(x^{l})$. Now $l$ divides $N$, so $y$ belongs to $Z_{B_{n}(\rp)}(x^{N})$, but $Z_{B_{n}(\rp)}(x^{N})=\ang{x}$ by~\cite[Proposition~9]{GGzeit}, and thus $Z_{B_{n}(\rp)}(x^{l})\subset \ang{x}$, from which we conclude that $Z_{B_{n}(\rp)}(x^{k})=\ang{x}$.

%

\item With the above notation, we have $N_{B_{n}(\rp)}(\ang{x^{l}})=N_{B_{n}(\rp)}(\ang{x^{k}})$. So it suffices to prove the result in the case where $k=l$. The inclusion $\ang{x,a\garside}\subset N_{B_{n}(\rp)}(\ang{x^{l}})$ follows from equations~\reqref{conja} and~\reqref{conjb}. So suppose that $y\in N_{B_{n}(\rp)}(\ang{x^{l}})$. Then $y^{2}\in Z_{B_{n}(\rp)}(x^{l})$ by \relem{conjinv}, and so $y$ is of finite order by part~(\ref{it:centraliser}). It follows from \repr{bnrpsubgp} that $N_{B_{n}(\rp)}(\ang{x^{l}})$ is finite. If $N=2$ and $x=a$, or if $N\geq 3$ then $\ang{x,a\garside}$ is maximal finite by \reth{finitebnrp2}, which proves that $N_{B_{n}(\rp)}(\ang{x^{l}})=\ang{x,a\garside}$ in these cases. So assume that either $N=1$, or that $N=2$ and $x=b$.
\begin{enumerate}[(i)]
\item If $N=1$ then $x=b$, $n=2$ and $l=1$. Then $x^{l}=b=\rho_{1}$, and so $x^{l} \in P_{2}(\rp)$. 
Now $\ang{x,a\garside[2]}= \ang{\rho_{1}, \rho_{2}\sigma_{1}^{2}}= \ang{\rho_{1}, \rho_{2}}$ by~\reqref{surfacerp2}, and thus $\ang{x,a\garside[2]}=P_{2}(\rp) \cong \quat$~\cite[p.~87]{vB}. Now $a\notin N_{B_{2}(\rp)}(\ang{\rho_{1}})$ by~\repr{permute}(\ref{it:permrhi}), so $P_{2}(\rp) \subset N_{B_{2}(\rp)}(\ang{x^{l}}) \subsetneqq B_{2}(\rp)$, from which it follows that $N_{B_{2}(\rp)}(\ang{x^{l}})=P_{2}(\rp)= \ang{x,a\garside[2]}$. 

\item If $N=2$ and $x=b$ then $l\in \brak{1,2}$. Then $n=3$, $x$ is of order $8$, and $\ang{x,a\garside}\cong \quat[16]$. If $N_{B_{3}(\rp)}(\ang{x^{l}})\cong \quat[16]$ then $\ang{x,a\garside}=N_{B_{3}(\rp)}(\ang{x^{l}})$ as required. So assume that $N_{B_{3}(\rp)}(\ang{x^{l}})\ncong \quat[16]$. By \reth{finitebnrp2}, $N_{B_{3}(\rp)}(\ang{x^{l}})\cong \oonestar$. However, $\oonestar$ has no normal cyclic subgroups of order greater than two~\cite[Proposition~85]{GG8}, which contradicts the fact that $\ang{x^{l}}$ is normal in $N_{B_{3}(\rp)}(\ang{x^{l}})$. The last part of the statement then follows. \qedhere
\end{enumerate}
\end{enumerate}
\end{proof}

The following result shows that each of the abstractly maximal dicyclic groups of \reth{finitebnrp2} possesses a single conjugacy class in $B_n(\rp)$.

\begin{prop}\label{prop:conjdicmax}
Let $n\geq 3$, and let $i\in \brak{0,1}$. Then any two subgroups of $B_{n}(\rp)$ that are abstractly isomorphic to $\dic{8(n-i)}$ are conjugate.
\end{prop}

\begin{proof}
Let $i\in \brak{0,1}$, and set $x=a$ and $y=\garside$ if $i=0$, and $x=b$ and $y=\garside a^{-1}$ if $i=1$. Let $G$ be a subgroup of $B_{n}(\rp)$ isomorphic to $\dic{8(n-i)}$. It suffices to show that $G$ is conjugate to the standard copy $\ang{x,y}$ of $\dic{8(n-i)}$ (see \repr{permute}(\ref{it:embed})). By~\cite{M} (see also \req{muragen} and \reth{ccrp2}(\ref{it:ccrp2a})), any two subgroups of $B_{n}(\rp)$ that are isomorphic to $\Z_{4(n-i)}$ are conjugate. Hence the unique cyclic subgroup of $G$ isomorphic to $\Z_{4(n-i)}$ is conjugate to $\ang{x}$. Conjugating $G$ if necessary, we may suppose that $G=\ang{x,z}$, where $z^{2}=x^{n-i}=\ft$, and $zxz^{-1}=x^{-1}$. Now $yxy^{-1}=x^{-1}$, so $y^{-1}z x z^{-1}y=x$, and hence $y^{-1}z\in Z_{B_{n}(\rp)}(x)$. So by \repr{centraliserx}, there exists $j\in \brak{0,1,\ldots,4(n-i)-1}$ such that $z=yx^{j}$, and thus $G=\ang{x,yx^{j}}=\ang{x,y}$ as required.
\end{proof}

It follows from the proof of \repr{conjdicmax} that if $G$ is a subgroup of $B_{n}(\rp)$ isomorphic to $\dic{8(n-i)}$ that contains the subgroup $\ang{x}$ then $G=\ang{x,y}$. Using \repr{centraliserx}, this may be generalised to obtain the following rigidity result concerning dicyclic groups of $B_{n}(\rp)$ that contain $a$ or $b$.

\begin{prop}\label{prop:rigiddyc}
Let $n\geq 3$, let $x\in \brak{a,b}$, let $N=n$ if $x=a$, and let $N=n-1$ if $x=b$. 
Suppose that 
$H$ is a dicyclic subgroup of $B_{n}(\rp)$ that contains $x^{k}$ for some $k\in \N$ for which $k/\gcd(k,N)$ is odd. If $l=\gcd(k,4N)$,
there exists
$0\leq c < l$ such that $H=\ang{x^{l},x^{c} a\garside}$.
\end{prop}

\begin{proof}
With the given hypotheses, $l$ divides $k$, and from the proof of \repr{centraliserx}(\ref{it:centraliser}), we have $\langle x^k\rangle=\langle x^l\rangle$. Since $H$ is dicyclic, there exists $y\in B_{n}(\rp)$ such that $H= \ang{x^{l},y}$ and $yx^{l} y^{-1}=x^{-l}$. Thus $y\in N_{B_{n}(\rp)}(\ang{x^{l}})\setminus \ang{x^{l}}$, and hence $y\in \ang{x}a\garside$ by \repr{centraliserx}(\ref{it:normaliser}). 
So there exists $\alpha\in \Z$ such that $y=x^{\alpha} a \garside$. Writing $\alpha=\beta l+c$, where $\beta\in \Z$ and $c\in \brak{0,1, \ldots, l-1}$, we obtain $H=\ang{x^{l},x^{\alpha} a \garside}= \ang{x^{l},x^{c} a\garside}$ 
as required.
\end{proof}

It is considerably more difficult to exhibit explicit embeddings in $B_{n}(\rp)$ of the binary polyhedral groups.
Nevertheless, we may obtain some information about such embeddings. For example, part~(\ref{it:binpolyb}) of the following proposition shows that for all $n\geq 4$, the intersection of any such subgroup with $P_{n}(\rp)$ is reduced to $\ang{\ft}$, and so is small as it can possibly be. In \repr{binpolyII}, we will see that for $n=3$, a similar result holds, but the intersection is slightly larger, isomorphic to $\quat$. Note that we do not need to consider the case $n=2$ because  $B_2(\rp)$ does not contain an embedded binary polyhedral group.

\begin{prop}\label{prop:binpolyI}\mbox{}
\begin{enumerate}
\item\label{it:binpolya} Let $n\geq 4$, 
let $x\in \brak{a,b}$, let $N=n$ if $x=a$ and let 
$N=n-1$ if $x=b$. Suppose that $k\in \N$ is such that $k/\gcd(k,N)$ is odd. 
Then no conjugate of $x^{k}$ is contained in any subgroup of $B_{n}(\rp)$ that is isomorphic to a binary polyhedral group.

\item\label{it:binpolyb}  Let $n\geq 4$, 
and let $G$ be a subgroup of $B_{n}(\rp)$ isomorphic to a binary polyhedral group. Then $G\cap P_{n}(\rp)= \ang{\ft}$. 
 In particular, $\pi(G)\cong \an[4]$ if $G\cong \tonestar$,  $\pi(G)\cong \sn[4]$ if $G\cong \oonestar$  and $\pi(G)\cong \an[5]$ if $G\cong \istar$. If $n=4$ then the first two isomorphisms are equalities. 
\end{enumerate}
\end{prop}

\begin{proof}\mbox{}
\begin{enumerate}
\item Let $n\geq 4$. Suppose that $k\in \N$ is such that $k/\gcd(k,N)$ is odd, and assume on the contrary that some conjugate of $x^{k}$ is belongs to a subgroup $H$ of $B_{n}(\rp)$ that is isomorphic to a binary polyhedral group. Replacing $k$ by $\gcd(k,4N)$ and making use of the proof of \repr{centraliserx}(\ref{it:centraliser}) if necessary, up to conjugation and taking an appropriate power of $x^{k}$, we may suppose that $k=N$ and so $x^{N}\in H\cap P_{n}(\rp)$. 
Consider the restriction $\map{\pi\left\lvert_{H}\right.}{H}[\sn]$ to $H$ of the permutation homomorphism $\pi$ of the short exact sequence~\reqref{sessn}. Since $x^{N}\in P_{n}(\rp)$, $x^{N}$ belongs to $\ker{\pi\left\lvert_{H}\right.}$. By~\cite[Proposition~85(b)]{GG8}, the normal subgroups of $\tonestar$ different from $\brak{e}$ are isomorphic to $\Z_{2}$, $\quat$ or $\tonestar$, those of $\oonestar$ are isomorphic to $\Z_{2},\quat,\tonestar$ or $\oonestar$, and those of $\istar$ are isomorphic to $\Z_{2}$ or $\istar$. Now $x^{N}$ is of order $4$, so $\ker{\pi\left\lvert_{H}\right.}\ncong \Z_{2}$, and since $\ker{\pi\left\lvert_{H}\right.}$ is a subgroup of $P_{n}(\rp)$, it follows using~\cite[Proposition~1]{GG3} that it cannot be isomorphic to $\quat,\tonestar, \oonestar$ or $\istar$ either, which yields a contradiction.

\item Let $G$ be a subgroup of $B_{n}(\rp)$ isomorphic to a binary polyhedral group. By~\cite[Proposition~1]{GG3}, $G\cap P_{n}(\rp)$ is isomorphic to a subgroup of $\Z_{4}$. Since $\ord{G}$ is even and $\ft\in P_{n}(\rp)$, it follows from \repr{agt}(\ref{it:agtb}) that $\ang{\ft}\subset G\cap P_{n}(\rp)$. With the above notation, using \repr{ord4Pn}, any element of $G\cap P_{n}(\rp)$ of order $4$ is conjugate in $B_{n}(\rp)$ to $x^{N}$, but by part~(\ref{it:binpolya}), $G$ does not contain any such conjugate. Thus $\ang{\ft}= G\cap P_{n}(\rp)$, which proves the first part of the statement. 
The second part of the statement follows from \reqref{sessn}, the fact that the quotient of $\tonestar$ (resp.\ $\oonestar, \istar$) by its unique subgroup of order $2$ is isomorphic to $\an[4]$ (resp.\ to $\sn[4], \an[5]$), and the fact that if $n=4$ then the target of $\pi$ is $\sn[4]$.\qedhere
\end{enumerate} 
\end{proof}

The result of the following proposition contrasts with that in the case of $B_{n}(\St)$~\cite[Theorem~1.3]{GGagt2}. 


\begin{prop}\label{prop:nomax}
Let $n\geq 5$, and let $i=0$ (resp.\ $i=1$) if $n$ is even (if $n$ is odd). Then $B_{n}(\rp)$ possesses  dicyclic subgroups of order $4(n-i)$ that are maximal but not abstractly maximal.
\end{prop}

\begin{rem}
As we shall see in Sections~\ref{sec:subgroupsb3} and~\ref{sec:subgroupsb4}, every finite maximal subgroup of $B_{3}(\rp)$ and $B_{4}(\rp)$ is also abstractly maximal. One may ask whether $\quat$ can be realised as a maximal  subgroup in $B_n(\rp)$ for some $n$. In \resec{subgroupsB5}, we will show that this is the case if $n=5$.
\end{rem}

\begin{proof}[Proof of \repr{nomax}] 
Let $n\geq 5$. If $n$ is even (resp.\ odd) let $i=0$ (resp. $i=1$), let $w=a$ (resp.\ $w=b$), and with the notation of~\reqref{muragen}, let $u=\delta_{0,i}$, so $u=\sigma_{1}\cdots \sigma_{n-1}$ (resp.\ $u=\sigma_{1}\cdots \sigma_{n-1}\sigma_{1}$), and let $z=\garside$ (resp.\ $z=\garside a^{-1}$). By \repr{conjdicmax}, $B_{n}(\rp)$ contains a single conjugacy class of subgroups isomorphic to $\dic{8(n-i)}$ of which a representative $G$ is given by $\ang{w,z}$. Using the structure of dicyclic groups, $G$ possesses exactly two subgroups $H_{1}=\ang{w^{2},z}$ and $H_{2}=\ang{w^{2},wz}$ isomorphic to $\dic{4(n-i)}$, and the unique cyclic subgroup of these subgroups of order $2(n-i)$ is generated by $w^{2}$. Now $u$ is also of order $2(n-i)$, and by \repr{nomaxcyclic}, it generates the unique cyclic subgroup of order $2(n-i)$ of a dicyclic subgroup $H_{3}$ of order $4(n-i)$. Since $4\mid 2(n-i)$, it follows from \reth{ccrp2}(\ref{it:ccrp2a}) that $\ang{w^{2}}$ and $\ang{u}$ generate non conjugate subgroups, which implies that $H_{3}$ is neither conjugate to $H_{1}$ nor to $H_{2}$.

Clearly $H_{3}$ is not abstractly maximal. Suppose that it is not maximal. Then it is strictly contained in a finite subgroup $L$ of $B_{n}(\rp)$. First suppose that $n=5$. By \reth{finitebnrp2}, the abstractly maximal 
finite subgroups of $B_{5}(\rp)$ are isomorphic to $\dic{40}$ or to $\quat[32]$. Since $H_{3}\cong \quat[16]$, we must have $L\cong \quat[32]$. From \repr{conjdicmax}, $L$ and $G$ are conjugate, and hence $H_{3}$ is conjugate to $H_{1}$ or $H_{2}$, which contradicts the conclusion of the previous paragraph, and thus $H_{3}$ is finite maximal, but not abstractly maximal. Now suppose that $n\geq 6$. Then $\ord{H_{3}}=4(n-i)\geq 24$, and since $H_{3}$ is dicyclic, it follows from~\cite[Proposition~85]{GG8} that it cannot be contained in a subgroup of $B_{n}(\rp)$ isomorphic to a binary polyhedral group. By Lagrange's theorem and \reth{finitebnrp2}, $L$ must therefore be isomorphic to a subgroup of $\dic{8(n-i)}$, but since $L$ contains strictly $H_{3}$, $L$ is of order $8(n-i)$, and so is conjugate to $G$ by \repr{conjdicmax}. As in the case $n=5$, we obtain a contradiction, and once more, $H_{3}$ is finite maximal, but not abstractly maximal.
\end{proof}

\section{The finite subgroups of $B_{3}(\rp)$} \label{sec:subgroupsb3}

In this section we will determine the structure of the finite subgroups of $B_{3}(\rp)$, the aim being to determine the subgroup lattice up to conjugacy, which we do in \resec{latticeB3}. In \resec{ccb3rp}, we first determine the number of such conjugacy classes, and then in \resec{centralB3} we compute the centralisers of the elements of order $4$. 

\subsection{Conjugacy classes of the finite subgroups of $B_{3}(\rp)$}\label{sec:ccb3rp}

Using \reth{finitebnrp2} and the subgroup structure of binary polyhedral and dicyclic groups, every  non-trivial finite subgroup of $B_{3}(\rp)$ is isomorphic to one of $\oonestar, \tonestar$, $\dic{24}, \quat[16], \quat$, $\dic{12}$ or $\Z_{m}$, where $m\in \brak{2,3,4,6,8,12}$. From \req{muragen},~\cite{M} and~\cite[Proposition~8 and Corollary~12]{GG3}, every finite-order element of $B_{3}(\rp)$ is conjugate to a power of one of the following elements: $\sigma_1\sigma_2$ (order $6$), $\rho_1\sigma_2$ (order $4$), $a=\rho_3\sigma_2\sigma_1$ (order $12$),  $\garside[3]=\sigma_1\sigma_2\sigma_1$ (order $4$), and $b=\rho_2\sigma_1$ (order $8$).
By \repr{conjdicmax}, there is a single conjugacy class of subgroups isomorphic to $\dic{24}$ or to $\quat[16]$, and by \reth{finitebnrp2} and \reco{ccorder4}, $B_3(\rp)$ possesses four conjugacy classes of subgroups isomorphic to $\Z_{4}$, and a single conjugacy class of subgroups isomorphic to $\Z_{m}$ for $m\in \brak{2,3,6,8,12}$. The number of conjugacy classes of the remaining subgroups $\oonestar, \tonestar, \dic{12}$ and $\quat$ may be determined using~\cite{BCG} and the computations of \ref{app:quatbnrp2}, but as we shall see, we may also calculate them directly using the structure of $B_3(\rp)$. 

We first study the elements of order $4$. By \reth{ccrp2} and \reco{ccorder4} and the above list of finite-order elements, every element of $B_{3}(\rp)$ of order $4$ is conjugate to one of $\rho_1\sigma_2$, $a^3$, $\garside[3]$ or $b^2$. These four elements are non conjugate since their images under Abelianisation are $(\overline{1},\overline{1})$, $(\overline{0},\overline{1})$, $(\overline{1},\overline{0})$ and $(\overline{0},\overline{0})$ respectively, and the two elements $a^3$ and $b^2$ belong to $P_3(\rp)$. From~\reqref{conja}, $a^{3}\garside[3]$ is of order $4$, and since $q(a^{3}\garside[3])=(\overline{1},\overline{1})$, $a^{3}\garside[3]$ is conjugate to $\rho_1\sigma_2$ (and is a more convenient representative of its conjugacy class). We now make use of this information to analyse the subgroups of $B_3(\rp)$ that are isomorphic to $\quat$. 


\begin{prop}\label{prop:conjclassesq8}\mbox{}
\begin{enumerate}
\item\label{it:conjclassesq8a} There are precisely three conjugacy classes of subgroups of $B_{3}(\rp)$ isomorphic to $\quat$ for which representatives are given by $\ang{b^2, \garside[3]a^{-1}}$, $\ang{b^2, b\garside[3]a^{-1}}$ and $\ang{a^3, \garside[3]}$.
\item\label{it:conjclassesq8b}
If $H$ is a subgroup of $B_{3}(\rp)$ isomorphic to $\quat$ then the centraliser of $H$ in $B_{3}(\rp)$ is equal to $\ang{\ft[3]}$. In particular, $B_3(\rp)$ does not possess any subgroup isomorphic to $\quat\times \Z$. 
\end{enumerate}
\end{prop}

\begin{rems}\mbox{}\label{rems:quatb3}
\begin{enumerate}[(a)] 
\item By \repr{BraMap}(\ref{it:BraMapa}), the number of conjugacy classes of subgroups of $B_{3}(\rp)$ isomorphic to $\quat$ is equal to that of subgroups of $\mcg[3]{\rp}$ isomorphic to $\Z_{2}\oplus \Z_{2}$, which by \reco{clasq8} of \ref{app:quatbnrp2} is indeed equal to $3$. 

\item\label{it:quatb3b} Since $\ang{b,\garside[3]a^{-1}}$ is the standard copy of $\quat[16]$ in $B_{3}(\rp)$, it follows from the general subgroup structure of dicyclic groups that $\ang{b^2,b\garside[3]a^{-1}}$ is isomorphic to $\quat$. As we shall see during the course of the proof of \repr{conjclassesq8}, up to conjugacy, $\ang{b^2,b\garside[3]a^{-1}}$, is the only subgroup of $B_{3}(\rp)$ isomorphic to $\quat$ whose elements of order $4$ are pairwise conjugate in $B_{3}(\rp)$. Further, this subgroup is contained in $P_{3}(\rp)$. Since $B_3(\rp)$ contains a copy of $\tonestar \cong \quat \rtimes \Z_{3}$, where the action permutes cyclically the elements $i,j$ and $k$ of $\quat$, it follows that
the elements of order $4$ of the copy of $\quat$ lying in $\tonestar$ are pairwise conjugate, and so up to conjugacy, if $G$ is a subgroup of $B_{3}(\rp)$ isomorphic to $\tonestar$ then the subgroup of $G$ isomorphic to $\quat$ is conjugate to $\ang{b^2,b\garside[3]a^{-1}}$, and the subgroups of $G$ isomorphic to $\Z_{4}$ are pairwise conjugate, which is indeed the case using the above analysis of the elements of order $4$.
\item\label{it:quatb3c} Any subgroup of $B_{3}(\rp)$ isomorphic to $\quat$ is non maximal. The first two copies given in \repr{conjclassesq8}(\ref{it:conjclassesq8a}) are subgroups of the standard copy $\ang{b, \garside[3]a^{-1}}$ of $\quat[16]$, and the third copy is a subgroup of the standard copy $\ang{a,\garside[3]}$ of $\dic{24}$.

\item Since the abstractly maximal subgroups of $B_{3}(\rp)$ are isomorphic to $\oonestar$ or $\dic{24}$, both of which contain subgroups isomorphic to $\quat$, it follows from \repr{conjclassesq8}(\ref{it:conjclassesq8b}) that $B_{3}(\rp)$ possesses very few infinite virtually cyclic subgroups. By~\cite[Theorem~10 and Proposition~11]{GJ}, $B_3(\rp)$ is isomorphic to $\oonestar \bigast_{\dic{12}} \dic{24}$, and up to isomorphism, its infinite virtually-cyclic subgroups are $\Z, \Z/2\times\Z, \Z/4\times\Z, \Z/4\bigast_{\Z/2}\Z/4$ and $\quat\bigast_{\Z/4}\quat$. Some of the computations in a preliminary version of this paper were used to prove certain results of~\cite{GJ}.
\end{enumerate}
\end{rems}

In order to prove \repr{conjclassesq8}, we shall make use of the following lemma.

\begin{lem}\label{lem:order4}
Let $G$ be a group possessing a unique element of order $2$, let $z\in G$ be of order $4$, and let $\lambda\in G$. Then $\lambda z$ is of order $4$ if and only if $z\lambda z^{-1}=\lambda^{-1}$. 
\end{lem}

\begin{proof}[Proof of \relem{order4}.]
Let $w$ denote the order $2$ element of $G$. Then $z^{2}=w$ by hypothesis, and we have:
\begin{align*}
\text{$\lambda z$ is of order $4$} \Longleftrightarrow \text{$(\lambda z)^2$ is of order $2$} \Longleftrightarrow \lambda z \lambda z=w=z^2 \Longleftrightarrow z \lambda z^{-1}= \lambda^{-1},
\end{align*}
which proves the result.
\end{proof}

\begin{proof}[Proof of \repr{conjclassesq8}.]\mbox{}
\begin{enumerate}
\item Let $H$ be a subgroup of $B_{3}(\rp)$ that is isomorphic to $\quat$. By \rerems{quatb3}(\ref{it:quatb3c}), the three subgroups given in the statement are subgroups of one of the standard copies of the maximal dicyclic groups $\dic{24}$ and $\quat[16]$ of $B_{3}(\rp)$, and using standard properties of dicyclic groups, they are each abstractly isomorphic to $\quat$. Further, they are non conjugate since under Abelianisation they project to non-conjugate subgroups in $\Z_{2}\oplus \Z_{2}$. Indeed, $q(\ang{b^2,\garside[3]a^{-1}})=\ang{(\overline{1}, \overline{1})}$, $q(\ang{b^2,b\garside[3]a^{-1}})=\brak{(\overline{0}, \overline{0})}$ and $q(\ang{a^3,\garside[3]})=\Z_{2} \oplus \Z_{2}$. So to prove the result, it suffices to show that $H$ is conjugate to one of these three copies of $\quat$. 

As we mentioned at the beginning of the section, every element of $B_{3}(\rp)$ of order $4$ is conjugate to one of $a^3 \garside[3],a^{3},\garside[3]$ or $b^{2}$, and its conjugacy class is determined by its Abelianisation. Now $a^3, b^2\in P_{3}(\rp)$ and $a^3 \garside[3], \garside[3]\in B_{3}(\rp)\setminus P_{3}(\rp)$. The subgroup $H$ contains three subgroups isomorphic to $\Z_{4}$, at least one of which must lie in $P_{3}(\rp)$, for if two of these three subgroups are both conjugate to $\ang{a^3 \garside[3]}$ or both conjugate to $\ang{\garside[3]}$ then projecting into $\Z_{2}\oplus \Z_{2}$ shows that the third subgroup is conjugate to $\ang{b^2}$, while if one of these subgroups is conjugate to $\ang{a^3 \garside[3]}$ and another is conjugate to $\ang{\garside[3]}$ then the third subgroup is conjugate to $\ang{a^3}$.

So taking $x\in \brak{a,b}$, $N=3$ and $z=\garside[3]$ if $x=a$, $N=2$ and  $z=\garside[3]a^{-1}$ if $x=b$, and setting $v=x^N$, up to conjugacy, we may suppose that $H=\ang{v,u}$, where $u\in H\setminus \ang{v}$ is of order $4$. 
Set $\xi=uz^{-1}$. Taking successively $\lambda=\xi$ and $\lambda=v\xi$, we obtain respectively $\lambda z=u$ and $\lambda z= vu$, which are of order $4$, and since $z$ is also of order $4$, \relem{order4} implies that $z\xi z^{-1}=\xi^{-1}$ and $z(v\xi)z^{-1}=(v\xi)^{-1}$ using \repr{agt}(\ref{it:agtb}). From these two equations and the fact that $zvz^{-1}=v^{-1}$ by equations~\reqref{conja} and~\reqref{conjb}, we obtain:
\begin{equation*}
v\xi= z(v\xi)^{-1}z^{-1}=z \xi^{-1} z^{-1} \ldotp zv^{-1}z^{-1} = \xi v,
\end{equation*}
hence $\xi \in Z_{B_{3}(\rp)}(v)$, and so $\xi\in \ang{x}$ by \repr{centraliserx}(\ref{it:centraliser}).
Now $u=\xi z$, and multiplying $u$ on the left by powers of $v$ if necessary, we obtain $H=\ang{v, x^{k} z}$, where $0\leq k\leq N-1$.
If $x=a$ then the three possibilities for $H$ are conjugate since by~\reqref{conja}, we have:
\begin{equation*}
\text{$\ang{a^3, a \garside[3]}=\ang{a^3, a^4 \garside[3]}=a^2 \ang{a^3, \garside[3]} a^{-2}$ and  
$\ang{a^3, a^2 \garside[3]}=a \ang{a^3, \garside[3]} a^{-1}$,}
\end{equation*}
It follows that any subgroup of $B_{3}(\rp)$ isomorphic to $\quat$ and of the form $\ang{a^3,u}$ is conjugate to $\ang{a^3, \garside[3]}$. If $x=b$ then $H=\ang{b^2, \garside[3]a^{-1}}$ if $k=0$, and $H=\ang{b^2, b\garside[3]a^{-1}}$ if $k=1$. We have thus shown that any subgroup of $B_{3}(\rp)$ isomorphic to $\quat$ is conjugate to one of the three given subgroups, and this completes the proof of part~(\ref{it:conjclassesq8a}).

\item
Let $H$ be a subgroup of $B_{3}(\rp)$ isomorphic to $\quat$. Then by part~(\ref{it:conjclassesq8a}), $H$ is conjugate in $B_{3}(\rp)$ to $H'$, where $H'$ is one of the three subgroups given there, and with the same notation as above, there exists $x\in \brak{a,b}$ such that $x^N\in H'$. It follows from \repr{centraliserx}(\ref{it:centraliser})
that $\ang{x}=Z_{B_{3}(\rp)}(x^N) \supset Z_{B_{3}(\rp)}(H')$. Since $x\notin Z_{B_{3}(\rp)}(H')$, we see that $Z_{B_{3}(\rp)}(H')=\ang{x^{2}}=\ang{\ft[3]}$, and we see using \repr{agt}(\ref{it:agtb}) that $Z_{B_{3}(\rp)}(H)=\ang{x^{2}}=\ang{\ft[3]}$. In particular, $Z_{B_{3}(\rp)}(H)$ contains no element of infinite order.\qedhere
\end{enumerate}
\end{proof}

We now turn our attention to the subgroups of $B_3(\rp)$ that are isomorphic to one of the binary polyhedral groups $\oonestar$ or $\tonestar$. Recall that $\oonestar$ contains $\tonestar$ as a subgroup of index $2$.
The following proposition exhibits explicit representatives of these two subgroups, and shows that there is a single conjugacy class of each, which is coherent with the results of 
\repr{conjclass345} of
\ref{app:quatbnrp2} 
regarding the number of conjugacy classes of subgroups of $\mcg[3]{\rp}$ isomorphic to $\sn[4]$ and $\an[4]$ respectively.

\begin{prop}\label{prop:ostar}
In $B_{3}(\rp)$:
\begin{enumerate}[(a)]
\item\label{it:ostar2} $\ang{\rho_{3}\rho_{2},\rho_{1}\rho_{2}, a^{4}}\cong \tonestar$ and $\ang{\rho_{3}\rho_{2},\rho_{1}\rho_{2}, a^{4}, a^{3}\garside[3]}\cong \oonestar$.
\item\label{it:ostar3} there is a single conjugacy class of subgroups isomorphic to $\tonestar$.
\item\label{it:ostar4} there is a single conjugacy class of subgroups isomorphic to $\oonestar$.
\end{enumerate}
\end{prop}

\begin{proof}
First recall that $\oonestar$ is generated by elements $X,P,Q$ and $R$ that are subject to the following relations~\cite[page~198]{Wo}:
\begin{equation}\label{eq:presostar}
\left\{ 
\begin{aligned}
& X^{3}=1,\,P^2=Q^2=R^2,\,PQP^{-1}=Q^{-1},\\
& XPX^{-1}=Q,\,XQX^{-1}=PQ,\\
& RXR^{-1}=X^{-1},\, RPR^{-1}=QP,\,RQR^{-1}=Q^{-1},
\end{aligned}
\right.
\end{equation}
where its index~$2$ subgroup isomorphic to $\tonestar$ is generated by $P,Q$ and $X$ subject to those relations of~\reqref{presostar} that do not involve $R$.

\begin{enumerate}[(a)]
\item Let $G$ be a subgroup of $B_{3}(\rp)$ isomorphic to $\tonestar$, and let $F$ be its unique subgroup that is isomorphic to $\quat$. 
Conjugating $G$ if necessary, we may suppose by \rerems{quatb3}(\ref{it:quatb3b}) that $F= \ang{b^2, b\garside[3]a^{-1}}$.
We now construct such a subgroup $G$ of $B_3(\rp)$. Using \repr{present}, we have:
\begin{align*}
b\garside[3]a^{-1}&=\rho_{2}\sigma_{1} \sigma_{1}\sigma_{2}\sigma_{1} \sigma_{1}^{-1}\sigma_{2}^{-1} \rho_{3}^{-1} =\rho_{2}\sigma_{1} \sigma_{1} \rho_{3}^{-1}= \sigma_{1}^{-1} (\sigma_{1}\rho_{2}\sigma_{1}) \rho_{3}^{-1} \sigma_{1}= \sigma_{1}^{-1} \rho_{1} \rho_{3}^{-1} \sigma_{1}\\
\sigma_{1} b^{2} \sigma_{1}^{-1}&= \sigma_{1} \rho_{2}\sigma_{1}\ldotp \sigma_{1}^{-1}\rho_{1} \sigma_{1}^{-1}=\rho_{1} \rho_{2}.
\end{align*}
Taking $F'=\sigma_{1} F \sigma_{1}^{-1}$, it follows that $F'=\ang{P,Q}\cong \quat$, where $P=\rho_{1} \rho_{2}$ and $Q=\rho_{3} \rho_{1}^{-1}$. 
Let $X=a^4$. Then $X$ is of order $3$, and using \repr{permute}(\ref{it:permrhi}), we may check that $XPX^{-1}=Q$ and $XQX^{-1}=\rho_{2}^{-1} \rho_{3}^{-1}= P^{-1} Q^{-1}= PQ$ because $F'\cong \quat$. Since $a^4\notin F'$, we conclude from~\reqref{presostar} that $\tonestar \cong \ang{P,Q,X}=  \ang{QP, P, X}= \ang{\rho_{3}\rho_{2},\rho_{1}\rho_{2}, a^{4}}$. 



Now let $R=a^{3}\garside[3]$, and let $K=\ang{P,Q,X,R}$.
Then $R^{2}=\ft[3]$ by \req{conja}, so $R$ is of order $4$. Further, using the structure of $\tonestar$, the elements of order $4$ of $\ang{P,Q,X}$ belong to $F'$, so belong to $P_3(\rp)$, but $\pi(R)= (1,3)$, whence $R\notin \ang{P,Q,X}$. Now:
\begin{align*}
RXR^{-1}&= a^{3}\garside[3] a^{4} \garside[3]^{-1} a^{-3} = a^{-4}= X^{-1}\\
RPR^{-1}&= a^{3}\garside[3] \rho_{1} \rho_{2} \garside[3]^{-1} a^{-3} =\rho_{3} \rho_{2}= QP\\
RQR^{-1}&= a^{3}\garside[3] \rho_{3} \rho_{1}^{-1} \garside[3]^{-1} a^{-3}= (\rho_{3}\rho_{1}^{-1})^{-1}=Q^{-1}
\end{align*}
by \reqref{conjrho}, \reqref{conja} and \repr{permute}(\ref{it:permrhi}), and it follows from~\reqref{presostar} that $K= \ang{\rho_{3}\rho_{2},\rho_{1}\rho_{2}, a^{4}, a^{3}\garside[3]}$ is isomorphic to $\oonestar$. 


\item Let $\Gamma_{1},\Gamma_{2}$ be two subgroups of $B_{3}(\rp)$ isomorphic to $\tonestar\cong \quat\rtimes \Z_{3}$. 
Let $i=1,2$. By \rerems{quatb3}(\ref{it:quatb3b}), the (unique) subgroup of $\Gamma_{i}$ that is isomorphic to $\quat$ is conjugate to $F$. So there exists a conjugate $\Gamma_{i}'$ of $\Gamma_{i}$ such that $F\subset \Gamma_{i}'$. Let $z_{i}\in \Gamma_{i}'$ be an element of order $3$. Then $\Gamma_{i}'= \ang{F\cup\brak{z_{i}}}$ using the structure of $\tonestar$. 
Replacing $z_{i}$ by $z_{i}^{-1}$ if necessary, we may suppose that $z_{i} b^{2} z_{i}^{-1}\in \brak{b\garside[3]a^{-1}, (b\garside[3]a^{-1})^{-1}}$.  
Thus $z_{2}^{-1}z_{1} b^{2} z_{1}^{-1}z_{2}\in \brak{b^{2}, b^{-2}}$. Hence $z_{2}^{-1}z_{1} \in N_{B_{3}(\rp)}(b^{2})$, and so by 
\repr{centraliserx}(\ref{it:normaliser}), there exist $0\leq k\leq 7$ and $l\in \brak{0,1}$ such that $z_{1}= z_{2}b^{k}(a\garside[3])^l$ using the fact that $(a\garside[3])^2=\ft[3]=b^4$. By \reth{ccrp2}(\ref{it:ccrp2b}) and \relem{conjinv}, $z_{1}$ and $z_{2}$ are conjugate, thus $q(b^k(a\garside[3])^l)= (\overline{0},\overline{0})$, 
and hence $k+l$ is even. If $l=0$ (resp.\ $l=1$) then $k$ is even (resp.\ odd), and $z_{1}=z_{2} b^{k}$ (resp.\ $z_{1}=z_{2} b^{k-1} \ldotp b a\garside[3]= z_{2} b^{k-1} \ldotp b \garside[3] a^{-1}$). Since $b^2, b \garside[3] a^{-1}\in F$, we conclude in both cases that $\Gamma_{1}'=\Gamma_{2}'$, and hence $\Gamma_{1}$ and $\Gamma_{2}$ are conjugate in $B_{3}(\rp)$.


\item Let $G_{1}$ and $G_{2}$ be subgroups of $B_{3}(\rp)$ that are isomorphic to $\oonestar$. From parts~(\ref{it:ostar2}) and~(\ref{it:ostar3}), for $i=1,2$ there exists a conjugate $G_{i}'$ of $G_{i}$ such that $\Gamma\subset G_{1}'\cap G_{2}'$, where $\Gamma$ is the copy of $\tonestar$ constructed in part~(\ref{it:ostar3}). Suppose that 
$G_{1}'\neq G_{2}'$, so $G_{1}'\cap G_{2}'=\Gamma$, and let $L=\ang{G_{1}'\cup G_{2}'}$. By \reth{finitebnrp2}, $L$ is infinite, so $L\cong \oonestar \bigast_{\tonestar} \oonestar$  by~\cite[Lemma~15]{GG3} and~\cite[Remark~24]{GG8}. Thus $L$ possesses an index $2$ subgroup that is isomorphic to a semi-direct product of the form $\tonestar\rtimes \Z$ (\emph{cf.} the proof of~\cite[Lemma~15]{GG3}), and since the subgroup of $\tonestar$ isomorphic to $\quat$ is characteristic in $\tonestar$, $L$ possesses a subgroup that is isomorphic to a semi-direct product of the form $\quat\rtimes \Z$, which contradicts \repr{conjclassesq8}(\ref{it:conjclassesq8b}). It follows that $G_1'=G_2'$, and hence $G_{1}$ and $G_{2}$ are conjugate in $B_{3}(\rp)$.\qedhere
\end{enumerate}
\end{proof}

\begin{rems}\mbox{}
\begin{enumerate}
\item It follows from the proof of \repr{ostar}(\ref{it:ostar3}) and~(\ref{it:ostar4}) that if $\Gamma_{1}$ and $\Gamma_{2}$ are subgroups of $B_{3}(\rp)$ isomorphic to $\tonestar$ (resp.\ to $\oonestar$) such that $\Gamma_{1}\cap \Gamma_{2}$ contains their subgroup isomorphic to $\quat$ (resp.\ to $\tonestar$) then $\Gamma_{1}=\Gamma_{2}$.
 
\item Using~\cite[Corollary~2]{GG5}, the subgroup $\ang{\rho_{3}\rho_{2},\rho_{1}\rho_{2}, a^{4}, a^{3}\garside[3]}$ of the statement of \repr{ostar}(\ref{it:ostar2}) gives rise to an explicit subgroup of $B_{6}(\St)$ that is isomorphic to $\oonestar$.

\item
It remains to determine the number of conjugacy classes in $B_3(\rp)$ for $\dic{12}$. This will be accomplished in \reco{conjdic12} after having computed the centraliser of $a^{2}$ in \repr{centregarside}.

\end{enumerate}
\end{rems}

\begin{prop}\label{prop:binpolyII}\mbox{} 
\begin{enumerate}
\item There exists a subgroup of $B_{3}(\rp)$ that is isomorphic to $\oonestar$ and that contains the subgroup $\ang{b}$.
\item If $k\in \N$, then $a^k$ belongs to subgroup of $B_{3}(\rp)$ that is isomorphic to a binary polyhedral group if and only if $k$ is even.

\item\label{it:binpolybII} 
Let $G$ be a subgroup of $B_{3}(\rp)$ isomorphic to a binary polyhedral group. 
Then $G\cap P_{3}(\rp)$ is isomorphic to $\quat$.
 In particular,  $\pi(G)\cong \mathbb{Z}_3$ if $G\cong \tonestar$ and   $\pi(G)= \sn[3]$ if $G\cong \oonestar$.
\end{enumerate}
\end{prop}

\begin{proof} 
Let $H= \ang{\rho_{3}\rho_{2},\rho_{1}\rho_{2}, a^{4}}$ and      $K=\ang{\rho_{3}\rho_{2},\rho_{1}\rho_{2}, a^{4}, a^{3}\garside[3]}$ be the copies of $\tonestar$ and $\oonestar$ respectively in $B_{3}(\rp)$ given in \repr{ostar}(\ref{it:ostar2}), and let $L=\ang{\rho_{3}\rho_{2},\rho_{1}\rho_{2}}$, so $L\subset H \subset K$. 
\begin{enumerate}
\item
\begin{enumerate}
\item As we saw at the beginning of this section, $B_{3}(\rp)$ contains a single conjugacy class of cyclic subgroups of order $8$, represented by $\ang{b}$. If $G$ is a subgroup of $B_{3}(\rp)$ isomorphic to $\oonestar$, it contains elements of order $8$, and so contains a conjugate of $\ang{b}$. Hence there exists a conjugate of $G$ that contains $\ang{b}$. 

\item First suppose that $k$ is even. Recall that $a$ is of order $12$ and that $a^{6}=\ft[3]$. Using the fact that $a^{2} =\ft[3] a^{-4}$, it follows from \repr{ostar}(\ref{it:ostar2}) that $a^{k}\in H$. Conversely, let $k\in \N$, and suppose that some conjugate $z$ of $a^{k}$ belongs to a subgroup of $B_{3}(\rp)$ that is isomorphic to $\tonestar$ or $\oonestar$. Using \repr{ostar}(\ref{it:ostar3}) and~(\ref{it:ostar4}) and the fact that $H\subset K$, and conjugating $z$ if necessary, we may suppose that $z\in K$. Since $K$ contains no element of order $12$, we see that $k\not\equiv 1,5,7,11 \bmod{12}$. Further, $q(K)=\ang{(\overline{1}, \overline{1})}$, and since $q(a^{3})= q(a^{9})= (\overline{0}, \overline{1})$, if follows also that $k\not\equiv 3,9 \bmod{12}$, hence $k$ is even. 
\end{enumerate}

\item Let $G$ be a subgroup of $B_{3}(\rp)$ that is isomorphic to a binary polyhedral group. Applying \repr{ostar}(\ref{it:ostar3}) and~(\ref{it:ostar4}), using the fact that $P_{3}(\rp)$ is normal in $B_{3}(\rp)$, and conjugating $G$ if necessary, we may suppose that $G=H$ or $G=K$. Now $L\cong \quat$, $[H:L]=3$ and $[K:H]=2$, and making use of the constructions of $H$ and $K$ in the proof of \repr{ostar}(\ref{it:ostar2}), we see that $H=\coprod_{i=0}^{2} a^{4i} L$ and $K= (\coprod_{i=0}^{2} a^{4i} L) \coprod (\coprod_{i=0}^{2} a^{3}\garside[3] a^{4i} L)$. Since $\pi(a)=(1,2,3)$, $\pi(a^{3}\garside[3])=(1,2)$ and $\pi(L)=\brak{\id}$, we have $H \cap P_{3}(\rp)= K \cap P_{3}(\rp)=L$, and thus $G\cap P_{3}(\rp)=L$ as required.\qedhere
\end{enumerate}
\end{proof}

\subsection{Centralisers of the finite-order elements of $B_3(\rp)$}\label{sec:centralB3}

The following two results show that the centraliser of most of the finite-order elements of $B_{3}(\rp)$ given at the beginning of \resec{subgroupsb3} are themselves finite. The question of such centralisers is interesting in its own right, and the results have already been applied in the classification of the virtually cyclic subgroups of $B_{3}(\rp)$~\cite{GJ}. 
As mentioned at the beginning of \resec{ccb3rp}, every finite-order element of $B_{3}(\rp)$ is conjugate to a power of one of $\sigma_1\sigma_2, \rho_1\sigma_2, a, \garside[3]$ and $b$. 
Making use of \repr{centraliserx}(\ref{it:centraliser}) and the notation defined there, if $x\in \brak{a,b}$ and $1\leq k \leq 4N$, we are able to determine $Z_{B_{3}(\rp)}(\ang{x^{k}})$, with the exception of the case $x=a$ and $k\in \brak{2,4}$. Further:
\begin{enumerate}
\item $a^{2}\ft[3]=a^{-4}$ and it follows from \repr{agt}(\ref{it:agtb}) that $a^{2}$ and $a^{4}$ have the same centraliser.
\item by \reth{ccrp2}(\ref{it:ccrp2b}), $\sigma_1\sigma_2$ (resp.\ $(\sigma_1\sigma_2)^{2}$) is conjugate to $a^{2}$ (resp.\ to $a^{4}$), so they have conjugate centralisers.
\end{enumerate}
For the case $x=a$ and $k=4$, it thus suffices to study $Z_{B_{3}(\rp)}(\ang{a^{2}})$, which we now do, together with the centraliser of $\garside[3]$ of order $4$, and that of $\sigma_{1} \sigma_{2}$. The centraliser of $\rho_1\sigma_2$ will be determined in \repr{centralr1s2b3}.

\begin{prop}\label{prop:centregarside}\mbox{}
\begin{enumerate}
\item\label{it:centrea2} $Z_{B_{3}(\rp)}(a^{2})=\ang{a}$.
\item\label{it:centregar3} $Z_{B_{3}(\rp)}(\garside[3])=\ang{\garside[3]}$.
\item\label{it:centresig12} $Z_{B_{3}(\rp)}(\sigma_{1}\sigma_{2})= \rho_{1}\sigma_{2}\sigma_{1}\ang{a} \sigma_{1}^{-1}\sigma_{2}^{-1} \rho_{1}^{-1}$.
\end{enumerate}
\end{prop}
\begin{proof}
Recall that $P_{3}(\rp)$ is generated by the set $\setl{A_{i,j},\, \rho_{k}}{\text{$1\leq i<j\leq 3$ and $1\leq k\leq 3$}}$, where $A_{i,j}$ is defined in terms of the generators of \repr{present} by~\reqref{defaij}. The elements $\rho_{3}$ and $A_{2,3}$ generate a free subgroup $\F[2]=\F[2](\rho_{3},A_{2,3})$ of $P_{3}(\rp)$ of rank $2$ that may be identified with the kernel of the projection $\map{p}{P_{3}(\rp)}[P_{2}(\rp)]$ defined geometrically by forgetting the third string. Recall that $P_2(\rp)\cong \quat$, and that using this projection, $P_3(\rp) \cong \F[2] \rtimes \quat$~\cite[pp.~87 and~97]{vB}. Further, $\F[2]$ has an index $2$ subgroup $\F[3]= \F[3](\rho_{3}^{2}, A_{2,3},\rho_{3} A_{2,3} \rho_{3}^{-1})$ of rank $3$ that is the kernel of the homomorphism $\F[2]\to \Z_{2}$ defined by sending $\rho_{3}$ to the generator $\overline{1}$ of $\Z_{2}$, and $A_{2,3}$ to $\overline{0}$. We first prove parts~(\ref{it:centrea2}) and~(\ref{it:centregar3}). Let $w\in \brak{\garside[3],a^{2}}$. To prove the result in these cases, we claim that it suffices to show that:
\begin{equation}\label{eq:centralw}
Z_{B_{3}(\rp)}(w)\cap \F[3]=\brak{e}.
\end{equation}
For the moment, we suppose that this is indeed the case. Let $z\in Z_{B_{3}(\rp)}(w)$. Now $P_{3}(\rp) \mathrel{\lhd} B_{3}(\rp)$, $[B_{3}(\rp):P_{3}(\rp)]=\ord{\sn[3]}= 6$, $\F[2]\mathrel{\lhd} P_{3}(\rp)$ and $P_{3}(\rp)/\F[2] \cong \quat$, and by construction, $\F[3]\mathrel{\lhd}\F[2]$ and $[\F[2]:\F[3]]=2$. Then $z^{6}\in P_{3}(\rp)$, $z^{24} \in \F[2]$ and $z^{48} \in \F[3]$, and each of these elements commutes with $w$. Then $Z_{B_{3}(\rp)}(w)$ consists entirely of elements of finite order by~\reqref{centralw}, and so is itself finite by \repr{bnrpsubgp}.
%
It follows from~\reth{finitebnrp2} and the structure of dicyclic and binary polyhedral groups (see~\cite[Proposition~85]{GG8} for example) that $Z_{B_{3}(\rp)}(w)$ is either finite cyclic, dicyclic or binary polyhedral. By construction, $\ang{w}$ is contained in $Z(Z_{B_{3}(\rp)}(w))$, and thus $Z_{B_{3}(\rp)}(w)$ cannot be dicyclic or binary polyhedral since the centre of such a group is isomorphic to $\Z_{2}$ while $w$ is of order $4$ or $6$, hence $Z_{B_{3}(\rp)}(w)$ is cyclic. If $w=a^{2}$ then $Z_{B_{3}(\rp)}(a^{2})$ contains the group $\ang{a}$, which is maximal as a finite cyclic group, and so $Z_{B_{3}(\rp)}(a^{2})=\ang{a}$, which proves part~(\ref{it:centrea2}).
%
%
If $w=\garside[3]$ then it belongs to finite cyclic group $Z_{B_{3}(\rp)}(\garside[3])$. If $Z_{B_{3}(\rp)}(\garside[3])$ contained strictly $\ang{\garside[3]}$, by Murasugi's classification~\cite{M} and \reqref{muragen}, 
it would be conjugate to the cyclic subgroups $\ang{b}$ or $\ang{a}$ of order $8$ and $12$ respectively, and since it contains $\garside[3]$ which is of order $4$, $\garside[3]$ would be conjugate to $a^{\pm 3}$ or $b^{\pm 2}$, which is impossible since these elements are pure while $\garside[3]$ is not. This implies that $Z_{B_{3}(\rp)}(\garside[3]) =\ang{\garside[3]}$, and proves part~(\ref{it:centregar3}).


Finally, to prove part~(\ref{it:centresig12}), let $w=\sigma_{1} \sigma_{2}$ and $\lambda=\rho_{1}\sigma_{2}\sigma_{1}$. By equations~\reqref{sirisi}--\reqref{surfacerp2} and~\reqref{defab}--\reqref{powerab}, we have:
\begin{align*}
\lambda a^{2} \lambda^{-1} &= \rho_{1}\sigma_{2}\sigma_{1} a^{2} (\rho_{1}\sigma_{2}\sigma_{1})^{-1} = \rho_{1}\rho_{3}^{-1} a^{3} \sigma_{1}^{-1}\sigma_{2}^{-1}  \rho_{1}^{-1}=
\rho_{1} \rho_{2}\rho_{1} \sigma_{1}^{-1}\sigma_{2}^{-1}  \rho_{1}^{-1}\\
&= \rho_{1} \ldotp \rho_{1} \rho_{2} \sigma_{1}^{2}\ldotp \sigma_{1}^{-1} \sigma_{2}^{-1}\rho_{1}^{-1}= \rho_{1}^{2} \sigma_{1}^{-1}\ldotp \sigma_{1} \rho_{2} \sigma_{1} \rho_{1}^{-1}\sigma_{2}^{-1}= \rho_{1}^{2} \sigma_{1}^{-1} \rho_{1}  \ldotp\rho_{1}^{-1}\sigma_{2}^{-1}\\
&= \rho_{1}^{2} \sigma_{1}^{-1} \sigma_{2}^{-1}= \sigma_{1} \sigma_{2}^{2} \sigma_{1} \ldotp \sigma_{1}^{-1} \sigma_{2}^{-1}=\sigma_{1} \sigma_{2}=w.
\end{align*}
This implies that the centralisers of $a^{2}$ and $w$ are conjugated by $\lambda$, and hence part~(\ref{it:centresig12}) follows from part~(\ref{it:centrea2}).

It thus remains to prove the validity of \req{centralw} for $w\in \brak{a^{2}, \garside[3]}$. We achieve this by computing the action by conjugation of $w$ on the basis elements of $\F[3]$ (in fact on the subgroup $\F[3]\oplus \ang{\ft[3]}$ of $B_{3}(\rp)$, which as we shall see, is invariant under this action). To do so, let $x=b^{2}= \rho_{2}\rho_{1}$ and $y=\rho_{2}A_{1,2}\rho_{3}^{-1}$ be the generators of the subgroup of $P_{3}(\rp)$ isomorphic to $\quat$ described in~\cite[Section~2]{GG6}.
So $x^2=y^2=\ft[3]$, $x^{-1}yx=y^{-1}$, and by~\cite[equations~(15), (17) and (18)]{GG6}, the action by conjugation of $x$ and $y$ on $\F[2]$ is given by:
\begin{equation}\label{eq:conjxy}
\text{$x\rho_{3}x^{-1}= \rho_{3}^{-1}$, $y\rho_{3}y^{-1}=A_{2,3}\rho_{3}$, 
$x A_{2,3}x^{-1}=\rho_{3}^{-1}A_{2,3}^{-1} \rho_{3}$ and $y A_{2,3}y^{-1}=A_{2,3}^{-1}$.}
\end{equation}
Applying these relations and \req{rhocomm}, we obtain:
\begin{align}
\rho_{3}yx=x\ldotp x^{-1} \rho_{3}y x &=\rho_{2}\rho_{1}\ldotp \rho_{3}^{-1} \ldotp \rho_{3} A_{1,2}^{-1} \rho_{2}^{-1}=\rho_{1}\label{eq:rho1}\\
\rho_{3}^{-1} A_{2,3}^{-1}y^{-1}&= \rho_{3}^{-1} A_{2,3}^{-1} \ldotp \rho_{3} A_{1,2}^{-1} \rho_{2}^{-1}= x A_{2,3}x^{-1} \ldotp A_{1,2}^{-1} \rho_{2}^{-1}\notag\\
&= \rho_{2}\rho_{1} A_{2,3} \rho_{1}^{-1}\rho_{2}^{-1} A_{1,2}^{-1} \rho_{2}^{-1}= \rho_{2} A_{2,3} \rho_{2}^{-1} A_{1,2}^{-1} \rho_{2}^{-1}=\rho_{2},\label{eq:rho2}
\end{align}
where we have used the following calculation:
\begin{equation*}
\text{$\rho_{2} A_{1,2} \rho_{2}= \sigma_{1}^{-1}(\sigma_{1}\rho_{2}\sigma_{1})\sigma_{1}^{-1}= \sigma_{1}^{-1} \rho_{1}^2 \sigma_{1}^{-1} = A_{2,3}$ by equations~\reqref{sirisi} and~\reqref{surfacerp2}.}
\end{equation*}
Equations~\reqref{rho1} and~\reqref{rho2} describe $\rho_{1}$ and $\rho_{2}$ as elements of $P_{3}(\rp)$ in terms of the semi-direct product structure $\F[2]\rtimes \quat$. We now compute the action by conjugation of $\garside[3]$ on the basis elements of $\F[3]$. By equations~\reqref{rn2} and~\reqref{defaij}, as well as the fact that $\ft[3]=A_{1,2} A_{1,3} A_{2,3}= A_{1,3} A_{2,3} A_{1,2}$,
we also have: 
\begin{equation}\label{eq:sig1sq}
A_{1,2} = \sigma_{2}^{-1} \sigma_{1}^{-2} \sigma_{2}^{-1} \ldotp \sigma_{2} \sigma_{1}^{2} \sigma_{2} \sigma_{1}^{2}= \rho_{3}^2 \ldotp A_{1,3} A_{2,3} A_{1,2}= \rho_{3}^2 \ft[3],
\end{equation}
so by \req{conjsigma}, we obtain:
\begin{equation}\label{eq:conjA23}
\garside[3] A_{2,3}\garside[3]^{-1}=A_{1,2}=\rho_{3}^2 \ft[3].
\end{equation}
Then by equations~\reqref{conjrho},~\reqref{conjxy} and~\reqref{rho1}, we have:
\begin{align}
\garside[3] \rho_{3}\garside[3]^{-1}&=\rho_{1}^{-1}= x^{-1}y^{-1} \rho_{3}^{-1}=xy \rho_{3}^{-1} y^{-1}x^{-1} \ldotp xy=x \rho_{3}^{-1}A_{2,3}^{-1}x^{-1} \ldotp xy\notag\\
&= \rho_{3}\rho_{3}^{-1}A_{2,3}\rho_{3}xy= A_{2,3}\rho_{3}xy.\label{eq:conjrho3}
\end{align}
Equations~\reqref{conjxy} and~\reqref{conjrho3} imply that:
\begin{align}
\garside[3] \rho_{3}^2\garside[3]^{-1}&= A_{2,3}\rho_{3}xy A_{2,3}\rho_{3}xy=A_{2,3}\rho_{3}xy A_{2,3}\rho_{3}y^{-1}x^{-1} (xy)^2\notag\\
&= A_{2,3}\rho_{3}x \rho_{3} x^{-1} \ft[3]=A_{2,3} \ft[3],\label{eq:conjrho3sq}
\end{align}
and from equations~\reqref{conjxy}, \reqref{conjA23} and~\reqref{conjrho3}, we obtain:
\begin{align}
\garside[3] \rho_{3} A_{2,3} \rho_{3}^{-1} \garside[3]^{-1}&= A_{2,3} \rho_{3} xy \rho_{3}^{2} \ft[3] y^{-1} x^{-1} \rho_{3}^{-1} A_{2,3}^{-1}= A_{2,3} \rho_{3} x (A_{2,3} \rho_{3}) x^{-1} \rho_{3}^{-1} A_{2,3}^{-1} \ft[3]\notag\\
&= A_{2,3} \rho_{3} (\rho_{3}^{-1} A_{2,3}^{-1}\rho_{3} \ldotp \rho_{3}^{-1}) \rho_{3}^{-1} A_{2,3}^{-1} \ft[3]=
(\rho_{3}^{-1} A_{2,3}^{-1})^2 \ft[3]\notag\\
&= \rho_{3}^{-2} \ldotp \rho_{3} A_{2,3}^{-1} \rho_{3}^{-1} \ldotp A_{2,3}^{-1} \ft[3],\label{eq:conjmixed}
\end{align}
using also the fact that $\ft[3]$ is central in $B_3(\rp)$. Equations~\reqref{conjA23},~\reqref{conjrho3sq} and~\reqref{conjmixed} thus show that $\garside[3]$ acts by conjugation on $\F[3] \oplus \ang{\ft[3]}$.
Setting $\alpha=\rho_{3}^{2}$, $\beta=A_{2,3}$ and $\gamma=\rho_{3}A_{2,3}\rho_{3}$, we see that $(\alpha,\beta,\gamma)$ is a basis of $\F[3]$, and that:
\begin{equation}\label{eq:actDelta3}
\garside[3] \alpha \garside[3]^{-1}= \beta\ft[3], \quad \garside[3] \beta \garside[3]^{-1}=\alpha \ft[3] \quad\text{and}\quad \garside[3] \gamma \garside[3]^{-1}=\gamma^{-1}.
\end{equation}
It follows from~\reqref{actDelta3} that no non-trivial element of $\F[3]$ commutes with $\garside[3]$,
which proves that \req{centralw} holds in the case $w=\garside[3]$. So let $w=a^{2}$. We determine the action by conjugation of $a^{-2}$ on the basis $(\alpha,\beta,\gamma)$ of $\F[3]$. By \repr{permute}(\ref{it:permsi}) and~(\ref{it:permrhi}) and equations~\reqref{conjxy},~\reqref{rho2} and~\reqref{sig1sq}, we see that:
\begin{align*}
a^{-2} \alpha a^{2} &= \rho_{2}^{-2}= (yA_{2,3}\rho_{3})^{2}=\rho_{3}A_{2,3}\rho_{3} y^{2}= \gamma \ft[3]\\
a^{-2} \beta a^{2} &=A_{1,2}^{-1}=\alpha^{-1}\ft[3]\\
a^{-2} \gamma a^{2} &=\rho_{2}^{-1}A_{1,2}^{-1} \rho_{2}^{-1}= \rho_{2}^{-1}\rho_{3}^{-2} \rho_{2}^{-1} \ft[3]= yA_{2,3}\rho_{3}^{-1} y A_{2,3}\rho_{3}\ft[3]=  A_{2,3}^{-1} =\beta^{-1}.
\end{align*}
Using these expressions, we see that no non-trivial element of $\F[3]$ commutes with $a^{2}$, which proves that \req{centralw} holds in the case $w=a^{2}$, and this concludes the proof of the proposition. 
%
%
\end{proof}

As a corollary, we may thus determine the number of conjugacy classes of $\dic{12}$ in $B_{3}(\rp)$. The result is coherent with the results of~\cite{BCG}.

\begin{cor}\label{cor:conjdic12}
In $B_{3}(\rp)$, there are two conjugacy classes of subgroups isomorphic to $\dic{12}$.
\end{cor}

\begin{proof}
Let $G$ be a subgroup of $B_{3}(\rp)$ isomorphic to $\dic{12}$. Recall that $G$ is generated by two elements $x$ and $y$ of order $6$ and $4$ respectively, subject to the relations $x^{3}=y^{2}$ and $yxy^{-1}=x^{-1}$, so $y\in N_{B_{3}(\rp)}(\ang{x})$. By Murasugi's classification and \repr{nomaxcyclic} if necessary, up to conjugacy, we may suppose that $x=a^{2}$. Using the fact that $Z_{B_{3}(\rp)}(\ang{a^{2}})=\ang{a}$ by \repr{centregarside}(\ref{it:centrea2}), and arguing as in the proof of \repr{centraliserx}(\ref{it:normaliser}), we see that $N_{B_{3}(\rp)}(\ang{a^{2}})=\ang{a, \garside[3]}$. It follows from the structure of dicyclic groups that $y=a^{k} \garside[3]$, where $0\leq k\leq 5$, so $G=\ang{a^{2}, a^{l} \garside[3]}$, where $l\in \brak{0,1}$, which are indeed dicyclic groups of order $12$ (they are subgroups of the standard copy of $\dic{24}$), and they are non conjugate since under Abelianisation, $G$ is sent to $\ang{(\overline{0},\overline{1})}$ (resp.\ to $\ang{(\overline{1},\overline{1})}$) if $l=0$ (resp.\ $l=1$).
\end{proof}

\begin{rem}
If $G$ is a subgroup of $B_{3}(\rp)$ that is isomorphic to $\dic{12}$ and that contains $\ang{a^{2}}$, the proof of \reco{conjdic12} implies that $G=\ang{a^{2}, \garside[3]}$ or $G=\ang{a^{2}, a\garside[3]}$.
\end{rem}

It follows from \repr{centregarside} and the list of finite-order elements given at the beginning of this section that with the exception of $\rho_{1} \sigma_{2}$, the centraliser of every element of $B_{3}(\rp)$ of order $4$ is finite. Since $\rho_{1}$ and $\sigma_{2}$ commute, they both belong to the centraliser of $\rho_{1} \sigma_{2}$, so $Z_{B_{3}(\rp)}(\rho_{1} \sigma_{2})$ is infinite, and as we shall now show, they generate this centraliser. This enables us to determine the centralisers of all of the finite-order elements of $B_{3}(\rp)$.

\begin{prop}\label{prop:centralr1s2b3}
We have $Z_{B_{3}(\rp)}(\rho_{1} \sigma_{2})=\ang{\rho_{1}, \sigma_{2}} \cong \Z\oplus \Z_{4}$.
\end{prop}

\begin{proof}
The proof is similar in spirit to that of the second part of \repr{centregarside}. We start by determining $Z_{B_{3}(\rp)}(\rho_{1} \sigma_{2}) \cap P_{3}(\rp)$. To do so, we compute the conjugates by $\rho_{1} \sigma_{2}$ of the generators of $\F[3]$ and $\quat$. We claim that the subgroup $\F[3]\oplus \ang{\ft[3]}$ of $B_{3}(\rp)$ is invariant under conjugation by $\rho_{1}\sigma_{2}$. Indeed:
\begin{equation}\label{eq:conjrho312}
\left\{\begin{aligned}
\rho_{1}\sigma_{2} \ldotp \rho_{3} \ldotp \sigma_{2}^{-1}\rho_{1}^{-1}&= \rho_{1}\rho_{2} A_{2,3}^{-1}\rho_{1}^{-1}= \rho_{3} yx (\rho_{3}^{-1} A_{2,3}^{-1} y^{-1} A_{2,3}^{-1}) x^{-1} y^{-1} \rho_{3}^{-1}\\
&=\rho_{3}^{2} y \rho_{3}^{-1} A_{2,3}^{-1} y^{-1} \ldotp y= \rho_{3} y\\
\rho_{1}\sigma_{2} \ldotp A_{2,3} \ldotp \sigma_{2}^{-1}\rho_{1}^{-1}&= A_{2,3}
\end{aligned}\right.
\end{equation}
by equations~\reqref{sirisi} and \reqref{conjxy}--\reqref{rho2}, so:
\begin{equation}\label{eq:actr1s2}
\left\{\begin{aligned}
\rho_{1}\sigma_{2} \alpha \sigma_{2}^{-1}\rho_{1}^{-1}&=(\rho_{3} y)^{2}=\rho_{3} y \rho_{3}y^{-1}\ldotp y^{2}=\rho_{3} A_{2,3}\rho_{3}\ft[3]=\gamma \ft[3]\\
\rho_{1}\sigma_{2} \beta \sigma_{2}^{-1}\rho_{1}^{-1}&=\beta\\
\rho_{1}\sigma_{2} \gamma \sigma_{2}^{-1}\rho_{1}^{-1}&= \rho_{3} y A_{2,3} \rho_{3} y= \rho_{3} \ldotp y A_{2,3} \rho_{3} y^{-1}\ldotp y^{2}=\rho_{3}^{2} \ft[3]=\alpha \ft[3],
\end{aligned}\right.
\end{equation}
by~\reqref{conjxy} and~\reqref{conjrho312}, which proves the claim.
On the other hand, using \repr{present} and equations~\reqref{conjxy}--\reqref{rho1}, we have:
\begin{align}
\rho_{1}\sigma_{2} x \sigma_{2}^{-1}\rho_{1}^{-1} &=  \rho_{1}\sigma_{2}^{2} \sigma_{2}^{-1}\rho_{2}\sigma_{2}^{-1}=\rho_{1} A_{2,3} \rho_{3}= \rho_{3} (yx A_{2,3} \rho_{3} x^{-1}y^{-1}) yx= yx\label{eq:r1s2x}\\
\rho_{1}\sigma_{2} y \sigma_{2}^{-1}\rho_{1}^{-1} &= \rho_{1}\sigma_{2}^{2} \sigma_{2}^{-1}\rho_{2}\sigma_{2}^{-1} \sigma_{2} A_{1,2}\sigma_{2}^{-1}\sigma_{2}^{2}\sigma_{2}^{-1}\rho_{3}^{-1}  \sigma_{2}^{-1}\rho_{1}^{-1}\notag\\
&= \rho_{1}\sigma_{2}^{2} \rho_{3} A_{1,3}A_{2,3}\rho_{2}^{-1} \rho_{1}^{-1}= \rho_{1} A_{2,3} \rho_{3}^{-1}\rho_{2}^{-1} \rho_{1}^{-1}\notag\\
&= \rho_{3} (yx A_{2,3} \rho_{3}^{-1} y A_{2,3} \rho_{3} x^{-1} y^{-1}) \rho_{3}^{-1}\notag\\
&= A_{2,3} \rho_{3} A_{2,3} \rho_{3} y^{-1} \rho_{3}^{-1} A_{2,3} \rho_{3} \rho_{3}^{-1} A_{2,3}^{-1} \rho_{3}^{-1} y \ldotp y^{-1}=y^{-1}.\label{eq:r1s2y}
\end{align}
We conclude that $Z_{B_{3}(\rp)}(\rho_{1} \sigma_{2}) \cap \quat=\ang{\ft[3]}$.

Now suppose that $u \in Z_{B_{3}(\rp)}(\rho_{1} \sigma_{2})\cap P_{3}(\rp)$. Making use of the decomposition $P_{3}(\rp)=\F[2](\rho_{3},A_{2,3})\rtimes \quat(x,y)$ and the fact that $\F[3]$ is the kernel of the homomorphism $\F[2] \to \Z_{2}$ defined at the beginning of the proof of \repr{centregarside}, there exist $z\in \F[3]$, $\epsilon\in \brak{0,1}$ and $\xi\in \quat$ such that $u=z\rho_{3}^{\epsilon} \xi$. Then:
\begin{equation}\label{eq:conjzrho}
z\rho_{3}^{\epsilon} \xi=u=\rho_{1} \sigma_{2} u\sigma_{2}^{-1}\rho_{1}^{-1}=\underbrace{\rho_{1} \sigma_{2} z\sigma_{2}^{-1}\rho_{1}^{-1} \rho_{3}^{\epsilon}}_{\in \F[2]\oplus \ang{\ft[3]}} \underbrace{y^{\epsilon} \rho_{1} \sigma_{2} \xi\sigma_{2}^{-1}\rho_{1}^{-1}}_{\in \quat},
\end{equation}
by equations~\reqref{conjrho312} and~\reqref{actr1s2}. From~\reqref{actr1s2}, if $z=z(\alpha, \beta, \gamma)$, then $\rho_{1} \sigma_{2} z\sigma_{2}^{-1}\rho_{1}^{-1}= z(\gamma,\beta,\alpha) \ftalt[3]{m}$, where $m\in \Z$ is the sum of the exponents of $\alpha$ and $\gamma$ in $z$. Comparing the factors of $P_{3}(\rp)$ in~\reqref{conjzrho} that belong to $\F[2]$ and then to $\quat$, 
we see that $z(\alpha, \beta, \gamma) \rho_{3}^{\epsilon}= z(\gamma, \beta, \alpha) \rho_{3}^{\epsilon}$, so $z(\alpha, \beta, \gamma)= z(\gamma, \beta, \alpha)$, and thus $z\in \ang{\beta}$ and $m=0$ by~\reqref{actr1s2}, and that:
\begin{equation}\label{eq:xiyxi}
\xi=y^{\epsilon}(\rho_{1} \sigma_{2} \xi\sigma_{2}^{-1}\rho_{1}^{-1}).
\end{equation}
A straightforward computation using equations~\reqref{r1s2x} and~\reqref{r1s2y} gives rise to the following table:
\begin{center}\renewcommand{\arraystretch}{1.3}
\begin{tabular}{|>{$}c<{$}||*{8}{>{$}c<{$}|}}\hline
\xi & 1 & \ft[3] & x & x^{-1} & y & y^{-1} & xy & (xy)^{-1}\\ \hline\hline
\rho_{1} \sigma_{2} \xi\sigma_{2}^{-1}\rho_{1}^{-1} & 1 & \ft[3] & yx & xy & y^{-1} & y & x^{-1} & x\\ \hline
y(\rho_{1} \sigma_{2} \xi\sigma_{2}^{-1}\rho_{1}^{-1}) & y & y^{-1} & x^{-1} & x & 1 & \ft[3] & xy & (xy)^{-1}\\ \hline
\end{tabular}
\end{center}
If $\epsilon=0$ then it follows from~\reqref{xiyxi} and the table that $\xi\in \ang{\ft[3]}$, and so $u=\beta^{k} \ftalt[3]{\delta}$, where $k\in \Z$ and $\delta\in \brak{0,1}$. If $\epsilon=1$ then using the table, \req{xiyxi} possesses two solutions in $\quat$, namely $\xi=xy$ and $\xi=(xy)^{-1}$. Now $\rho_{3} (xy)^{-1}=\rho_{1}$ by~\reqref{rho1}, and $\rho_{3} xy=\rho_{3} (xy)^{-1} (xy)^{2}=\rho_{1}\ft[3]$, so $z=\beta^{k} \rho_{1} \ftalt[3]{\delta}$, where $k \in \Z$ and $\delta \in \brak{0,1}$. Setting $X=\setl{\sigma_{2}^{2k} \rho_{1}^{\epsilon}}{\text{$k\in \Z$ and $0\leq \epsilon\leq 3$}}$, we conclude from these two cases that $Z_{B_{3}(\rp)}(\rho_{1} \sigma_{2})\cap P_{3}(\rp) \subset X$,
where 
we have used the fact that $\rho_{1}\sigma_{2}$ is of order $4$, which implies that $\ft[3]= (\rho_{1}\sigma_{2})^{2}= \rho_{1}^{2}\sigma_{2}^{2}$ and $\rho_{1}^{4}=\sigma_{2}^{-4}$. It is clear that $X\subset Z_{B_{3}(\rp)}(\rho_{1} \sigma_{2})\cap P_{3}(\rp)$, so $X= Z_{B_{3}(\rp)}(\rho_{1} \sigma_{2})\cap P_{3}(\rp)$. 

Now let $\tau\in Z_{B_{3}(\rp)}(\rho_{1} \sigma_{2})$. Then $\pi(\tau)\in Z_{\sn[3]}((2,3))$, so $\pi(\tau)\in \ang{(2,3)}$. Since $\sigma_{2}$ commutes with $\rho_{1}\sigma_{2}$, 
it follows from above that there exist $u\in X$ and $\epsilon \in \brak{0,1}$ such that $\tau=\sigma_{2}^{\epsilon}u$, and then that $Z_{B_{3}(\rp)}(\rho_{1} \sigma_{2})=\setl{\sigma_{2}^{l} \rho_{1}^{\epsilon}}{\text{$l\in \Z$ and $0\leq \epsilon\leq 3$}}=\ang{\rho_{1},\sigma_{2}}$, which proves the equality in the statement.
To obtain the isomorphism, observe that $\ang{\rho_{1},\sigma_{2}}$ is Abelian, and is subject to the relation $\rho_{1}^{4}=\sigma_{2}^{-4}$, which implies that it is isomorphic to a quotient of $\Z\oplus \Z_{4}$. Any other (independent) relation would be of the form $\rho_{1}^{r}=\sigma_{2}^{s}$, where $(r,s)\in \Z^{2} \setminus \setl{(4t,-4t)}{t\in \Z}$. Since $\rho_{1}$ is a non-trivial element of $P_{3}(\rp)$ of infinite order, this relation is equivalent to one of the form $\rho_{1}^{r'}=\sigma_{2}^{2}$, where $r'\in \Z$. Conjugating~\reqref{sig1sq} by $\garside[3]$ and making use of~\reqref{conjsigma}--\reqref{conjrho}, we obtain $\sigma_{2}^{2}=\rho_{1}^{-2}\ft[3]$, so $\rho_{1}^{r'+2}=\ft[3]$, which is impossible since $\rho_{1}$ is of infinite order. Hence $Z_{B_{3}(\rp)}(\rho_{1} \sigma_{2})\cong \Z\oplus \Z_{4}$.
\end{proof}


\subsection{The subgroup lattice of $B_{3}(\rp)$}\label{sec:latticeB3}

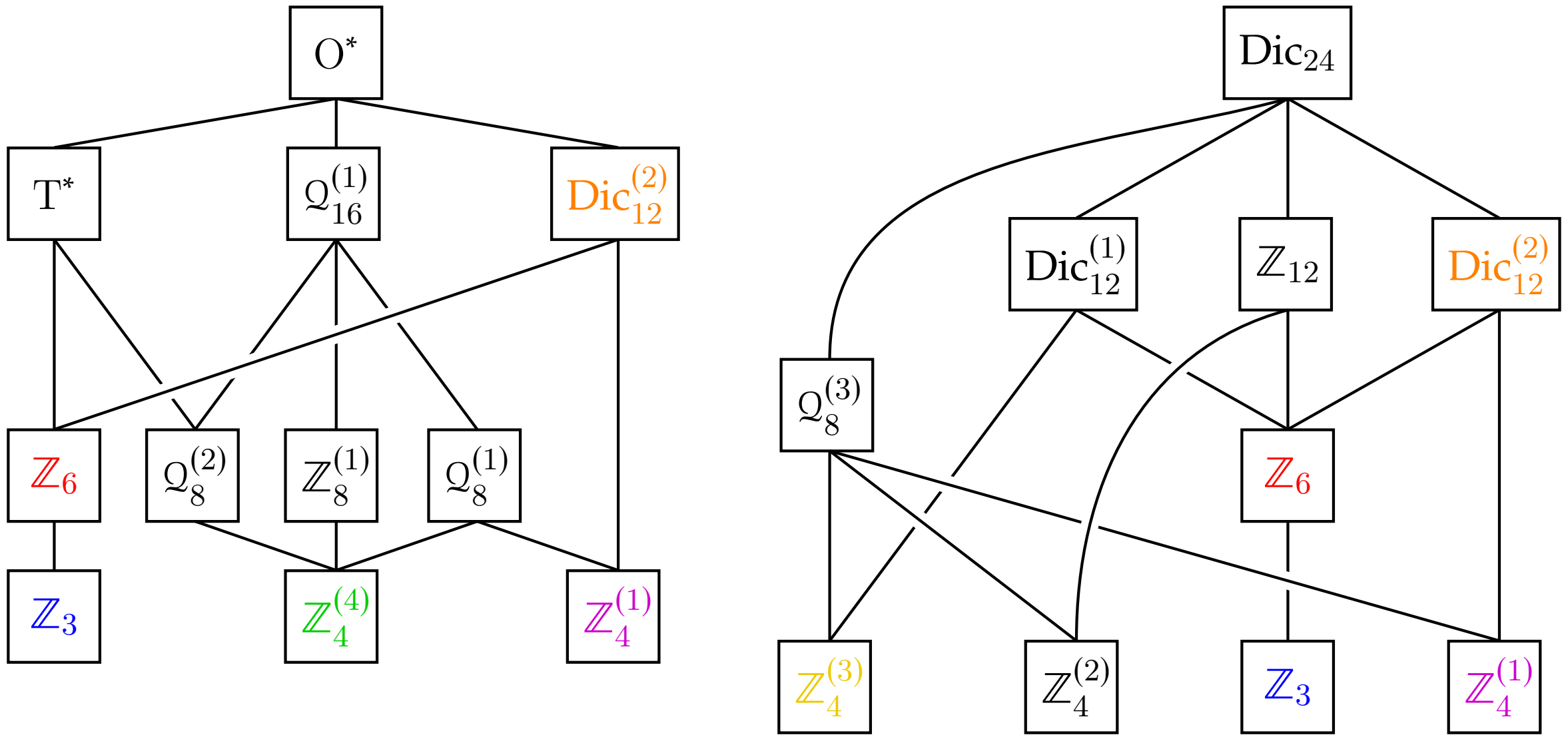
\begin{figure}[!h]
\hfill
\begin{tikzpicture}[scale=0.7]
\tikzstyle{every node} =[node contents={}, minimum height=0.9cm,minimum width=0.9cm]
\node (a) at (5,11) {$\oonestar$};
\node (b) at (1,9){$\tonestar$};
\node (c) at (5,9){$\quat[16]^{(1)}$};
\node (d) at (9,9){\textcolor{orange}{$\dic{12}^{(2)}$}};
\foreach \toa in {b,c,d} \draw[thick] (a.south)--(\toa.north);
\node (n) at (1,5){\textcolor{red}{$\Z_{6}$}};
\node (f) at (3,5){$\quat^{(2)}$};
\node (g) at (5,5){$\Z_{8}^{(1)}$};
\node (m) at (7,5){$\quat^{(1)}$};
\node (j) at (5,3){\textcolor{Green3}{$\Z_{4}^{(4)}$}};
\node (k) at (9,3){\textcolor{Magenta3}{$\Z_{4}^{(1)}$}};
\node (e) at (1,3){\textcolor{blue}{$\Z_{3}$}};
\foreach \tob in {f} \draw[thick] (b.south)--(\tob.north);
\foreach \toc in {f,g,m} \draw[thick] (c.south)--(\toc.north);
\foreach \tof in {j} \draw[thick] (f.south)--(\tof.north);
\foreach \tof in {j} \draw[thick] (g.south)--(\tof.north);
\foreach \tof in {e} \draw[thick] (n.south)--(\tof.north);
\foreach \tof in {j,k} \draw[thick] (m.south)--(\tof.north);
\foreach \tob in {n} \draw[line width=5pt,white] (d.south)--(\tob.north);
\foreach \tob in {n,k} \draw[thick] (d.south)--(\tob.north);

\node (aa) at (18.5,11) {$\dic{24}$};
\node (ab) at (15.5,8) {$\dic{12}^{(1)}$};
\node (ac) at (12,6) {$\quat^{(3)}$};
\node (ad) at (18.5,8) {$\Z_{12}$};
\node (ae) at (21.5,8) {\textcolor{orange}{$\dic{12}^{(2)}$}};
\node (af) at (18.5,5) {\textcolor{red}{$\Z_{6}$}};
\node (ag) at (18.5,2) {\textcolor{blue}{$\Z_{3}$}};
\node (ah) at (21.5,2) {\textcolor{Magenta3}{$\Z_{4}^{(1)}$}};
\node (ai) at (15.5,2) {$\Z_{4}^{(2)}$};
\node (aj) at (12,2) {\textcolor{Gold2}{$\Z_{4}^{(3)}$}};
\foreach \tob in {ab,ad,ae} \draw[thick] (aa.south)--(\tob.north);
\foreach \tob in {aj,af} \draw[thick] (ab.south)--(\tob.north);
\foreach \tob in {ag} \draw[thick] (af.south)--(\tob.north);

\path[thick] (aa.south) edge [out=-165, in=90] (ac.north);

\draw[line width=5pt,white] ($0.9*(ac.south) + 0.1*(ah.north)$)--($0.1*(ac.south) + 0.9*(ah.north)$);

\foreach \tob in {ah} \draw[thick] (ac.south)--(\tob.north);
\foreach \tob in {ai} \path[line width=5pt,white] (ad.south) edge [out=-165, in=90] (\tob.north);
\foreach \tob in {ah,af} \draw[thick] (ae.south)--(\tob.north);

\draw[line width=5pt,white] ($0.9*(ac.south) + 0.1*(ai.north)$)--($0.1*(ac.south) + 0.9*(ai.north)$);

\foreach \tob in {ai} \draw[thick] (ac.south)--(\tob.north);

\foreach \tob in {ai} \draw[thick] (ad.south) edge [out=-165, in=90] (\tob.north);
\foreach \tob in {aj} \draw[thick] (ac.south)--(\tob.north);

\foreach \tob in {af} \draw[thick] (ad.south)--(\tob.north);
\foreach \tob in {n} \draw[thick] (b.south)--(\tob.north);
\foreach \boxx in {d,aa,ab,ae} \draw[thick] (\boxx.south west) rectangle ++(1.8,1.3); 
\foreach \boxx in {a,b,c,n,f,g,m,j,k,e,ac,ad,af,ag,ah,ai,aj} \draw[thick] (\boxx.south west) rectangle ++(1.3,1.3); 
\end{tikzpicture}\hspace*{\fill}
\caption{The lattice of the conjugacy classes of the finite subgroups of $B_{3}(\rp)$.}\label{fig:omegai}
\end{figure}

We are now in a position to make explicit the lattice of the conjugacy classes of the finite subgroups of $B_{3}(\rp)$ (see Figure~\ref{fig:omegai} on page~\pageref{fig:omegai}). To simplify the discussion and Figure~\ref{fig:omegai}, we omit the isomorphism classes of $\Z_{2}$ and $\brak{e}$ in what follows.
Table~\ref{tab:ccB3} summarises the number of such conjugacy classes for the remaining subgroups. The third column refers to the result where the number of conjugacy classes is computed. For those isomorphism classes where there are at least two conjugacy classes, and also for $\quat[16]$ and  and $\Z_{8}$, 
we number the copies as shown for future reference. 
\begin{table}[ht]
\begin{center}\renewcommand{\arraystretch}{1.3}\setlength{\tabcolsep}{2pt}
\begin{tabular}{|>{$}c<{$}||>{$}c<{$}|c|}\hline
G &  \text{representatives of the conjugacy classes of $G$} & result\\ \hline\hline
\oonestar &  \ang{\rho_{3}\rho_{2},\rho_{1}\rho_{2}, a^{4}, a^{3}\garside[3]} & \repr{ostar}(\ref{it:ostar4})\\ \hline
\tonestar & \ang{\rho_{3}\rho_{2},\rho_{1}\rho_{2}, a^{4}} & \repr{ostar}(\ref{it:ostar3})\\ \hline
\dic{24} &  \ang{a, \garside[3]}& \repr{conjdicmax}\\  \hline
\quat[16] &  \quat[16]^{(1)}= \ang{b, \garside[3] a^{-1}}& \repr{conjdicmax}\\  \hline
\dic{12} &  \text{$\dic{12}^{(1)}= \ang{a^{2}, \garside[3]}$, $\dic{12}^{(2)}= \ang{a^{2}, a\garside[3]}$}& \reco{conjdic12}\\  \hline
\quat[8]  &  \text{$\quat^{(1)}= \ang{b^2, \garside[3]a^{-1}}$, $\quat^{(2)}= \ang{b^2, b\garside[3]a^{-1}}$, $\quat^{(3)}= \ang{a^3, \garside[3]}$} & \repr{conjclassesq8}(\ref{it:conjclassesq8a}) \\  \hline
\Z_{12}  &  \ang{a} & \cite{M} \& \reco{ccorder4}(\ref{it:ccorder4b})\\  \hline
\Z_{8} &  \Z_{8}^{(1)}= \ang{b} & \cite{M} \& \reco{ccorder4}(\ref{it:ccorder4b})\\  \hline
\Z_{6} &  \ang{a^{2}}& \cite{M} \& \reth{ccrp2}(\ref{it:ccrp2b})\\  \hline
\Z_{4}  &  \text{$\Z_{4}^{(1)}= \ang{\rho_1\sigma_2}$, $\Z_{4}^{(2)}= \ang{a^{3}}$, $\Z_{4}^{(3)}= \ang{\garside[3]}$, $\Z_{4}^{(4)}= \ang{b^{2}}$}& \cite{M} \& \reth{ccrp2}(\ref{it:ccrp2a})\\  \hline
\Z_{3} &  \ang{a^{4}}& \cite{M} \& \reth{ccrp2}(\ref{it:ccrp2b})\\  \hline
\end{tabular}
\end{center}
\caption{The conjugacy classes of the finite subgroups of $B_{3}(\rp)$}
\label{tab:ccB3}
\end{table}
Observe that there are a total of seventeen conjugacy classes of finite subgroups of order at least three in $B_{3}(\rp)$ that correspond to the vertices of the lattice of Figure~\ref{fig:omegai}. It now remains to determine the inclusions (up to conjugacy) between the given representatives of these conjugacy classes. Many of them may be obtained directly from the choice of representatives given in Table~\ref{tab:ccB3}. We analyse the subgroups of the two abstractly maximal subgroups $\oonestar$ and $\dic{24}$ of $B_{3}(\rp)$.
\begin{enumerate}[(a)]
\item\label{it:OstarB3} Abstractly, $\tonestar,\quat[16]$ and $\dic{12}$ are subgroups of $\oonestar$, $\Z_{6}$ is a subgroup of $\tonestar$ and of $\dic{12}$, $\Z_{3}$ is a subgroup of $\Z_{6}$, $\Z_{8}$ is a subgroup  of $\quat[16]$, $\dic{4m}$ has exactly two subgroups isomorphic to $\dic{2m}$ for all $m\geq 4$ even, the subgroups of $\tonestar\cong \quat \rtimes \Z_{3}$ isomorphic to $\Z_{4}$ are pairwise conjugate. If we now consider the copies of these subgroups embedded in $B_{3}(\rp)$, we see that $\quat^{(2)}$ and $\Z_{4}^{(4)}$ are respectively the copies of $\quat$ and $\Z_{4}$ lying in $\tonestar$, $\Z_{4}^{(4)}$ is the copy of $\Z_{4}$ lying in $\quat^{(2)}$, and $\quat^{(1)}$ and $\quat^{(2)}$ are the two copies of $\quat$ lying in $\quat[16]$. Since $q(\oonestar)= q(\dic{12}^{(2)})= \ang{(\overline{1}, \overline{1})}$ and $q(\dic{12}^{(1)})= \ang{(\overline{0}, \overline{1})}$, $\dic{12}^{(2)}$ is the copy of $\dic{12}$ lying in $\oonestar$. In the same manner, since $q(\quat^{(1)})= q(\Z_{4}^{(1)})= \ang{(\overline{1}, \overline{1})}$, we see that $\quat^{(1)}$ and $\oonestar$ contain $\Z_{4}^{(1)}$ and $\Z_{4}^{(4)}$ as their copies of $\Z_{4}$, and $\dic{12}^{(2)}$ contains $\Z_{4}^{(1)}$ as its copy of $\Z_{4}$. It follows from this and the fact that most of the finite subgroups admit a single conjugacy class that none of the abstractly non-maximal subgroups of $\oonestar$ are maximal.

\item Abstractly, $\dic{24}$ possesses two subgroups isomorphic to $\dic{12}$, and one subgroup isomorphic to each of $\Z_{3}, \Z_{6},\Z_{12}$ and $\quat$. From case~(\ref{it:OstarB3}), the two copies of $\dic{12}$ in $B_{3}(\rp)$ are non conjugate and by Table~\ref{tab:ccB3}, the copy of $\quat$ is $\quat^{(3)}$, and $Z_{12}$ has $\Z_{4}^{(2)}$ as its unique subgroup of order $4$. As in case~(\ref{it:OstarB3}), $\dic{12}^{(2)}$ contains the copy of $\Z_{6}$, and its unique subgroup of order $4$ is $\Z_{4}^{(1)}$. The subgroup $\dic{12}^{(2)}$ contains $\Z_{6}$, and by case~(\ref{it:OstarB3}) and Table~\ref{tab:ccB3}, we see that it contains $\Z_{4}^{(3)}$ as its unique subgroup of order $4$. Finally, $q(\quat^{(3)})\cong \Z_{2} \oplus \Z_{2}$, and we conclude that $\quat^{(3)}$ contains $\Z_{4}^{(i)}$ for $i=1,2,3$. Once more, we see that none of the abstractly non-maximal subgroups of $\oonestar$ are maximal.
\end{enumerate}

From this analysis, we obtain Figure~\ref{fig:omegai}, where for the sake of clarity, we have separated the subgroup lattices of $\oonestar$ and of $\dic{24}$ in $B_{3}(\rp)$. Some groups (in coloured text) appear in both, and two other groups ($\Z_4^{(3)}$ and $\Z_4^{(4)}$) are also coloured for later use. The lines joining the different groups represent inclusions between the corresponding conjugacy classes of the representatives appearing in Table~\ref{tab:ccB3}, and are not necessarily honest inclusions. 
As a consequence, we conclude that in $B_{3}(\rp)$, every maximal finite subgroup is abstractly maximal.

\section{The subgroup structure of $B_{4}(\rp)$}\label{sec:subgroupsb4}

In this section we prove several results about the elements and subgroups of $B_4(\rp)$ of finite order. In \resec{genconsid}, we recall the conjugacy classes of the former and isomorphism classes of the latter, and in \resec{centralB4}, we compute the centralisers of the most of the elements of finite order.
 In \resec{q8b4}, we construct five (resp.\ three, two) copies of $\quat$ (resp.\ $\quat[16]$, \dic{12}), which using~\cite{BCG} yields representatives of all of their conjugacy classes  in $B_{4}(\rp)$. In \resec{polyB4}, we provide explicit constructions of $\tonestar$ and $\oonestar$ as subgroups of $B_4(\rp)$, and in \resec{latticeB4}, we bring together all of this information to determine the lattice of the conjugacy classes of the finite subgroups of $B_{4}(\rp)$.
 
\subsection{The elements and subgroups of $B_{4}(\rp)$ of finite order}\label{sec:genconsid}

From~\cite{M} and~\cite[Proposition~8 and Corollary~12]{GG3}, every finite order element of $B_{4}(\rp)$ of order at least three is conjugate to a power of one of the following elements: 
\begin{enumerate}[(a)]
\item $\sigma_1\sigma_2\sigma_3$ (order $8$, cycle type $(\cdot\;\cdot\;\cdot\;\cdot)$ and Abelianisation $(\overline{1},\overline{0})$).
\item $(\rho_2\sigma_1)^2\sigma_3$ (order $4$, cycle type $(\cdot\;\cdot)$ and Abelianisation $(\overline{1}, \overline{0})$).
\item $a=\rho_4\sigma_3\sigma_2\sigma_1$ (order $16$, cycle type $(\cdot\;\cdot\;\cdotp \;\cdot)$ and Abelianisation $(\overline{1}, \overline{1})$).
\item $\sigma_1\sigma_2\sigma_3\sigma_1$ (order $6$, cycle type $(\cdot\;\cdot \; \cdot)$ and Abelianisation $(\overline{0},\overline{0})$).
\item $\rho_1\sigma_2\sigma_3\sigma_2$ (order $4$, cycle type $(\cdot\;\cdot)$ and Abelianisation $(\overline{1}, \overline{1})$).
\item $b=\rho_3\sigma_2\sigma_1$ (order $12$, cycle type $(\cdot\;\cdot\; \cdot)$ and Abelianisation $(\overline{0}, \overline{1})$).
\end{enumerate}
Up to conjugacy, the representatives of the elements of order $4$ of $B_{4}(\rp)$ are:
\begin{enumerate}[(a)]
\item $(\rho_2\sigma_1)^2\sigma_3$ (cycle type $(\cdotp \;\cdot)$ and Abelianisation $(\overline{1}, \overline{0})$).

\item $a^4$ (cycle type $\id$ and Abelianisation $(\overline{0},\overline{0})$).
\item $\rho_1\sigma_2\sigma_3\sigma_2$ (cycle type $(\cdotp \;\cdot)$ and Abelianisation $(\overline{1},\overline{1})$).
\item $(\sigma_1\sigma_2\sigma_3)^2$ (cycle type $(\cdotp \;\cdot)(\cdotp \;\cdotp)$ and Abelianisation $(\overline{0},\overline{0})$).
\item $b^3$ (cycle type $\id$ and Abelianisation $(\overline{0}, \overline{1})$).
\end{enumerate}
Since the cycle type and Abelianisation are conjugacy invariants, they may be used to distinguish the conjugacy classes of the six Murasugi elements of $B_{4}(\rp)$ given above, as well as those of the elements of order $4$. In particular, the above five elements of order $4$ and the cyclic subgroups that they generate are pairwise non-conjugate. This yields representatives of the five conjugacy classes of $\Z_{4}$ given by \reco{ccorder4}(\ref{it:ccorder4a}).

By \reth{finitebnrp2}, the isomorphism classes of the non-trivial finite subgroups of $B_{4}(\rp)$ are given by $\oonestar, \quat[32], \dic{24}, \tonestar, \quat[16], \dic{12}, \quat$ and $\Z_{m}$ for $m\in\brak{2,3,4,6,8,12,16}$. By \reco{ccorder4}, for the isomorphism classes of the cyclic subgroups, there is a single conjugacy class of each, with the exception of $\Z_{8}$ (resp.\ $\Z_{4}$), for which there are two (resp.\ five) conjugacy classes, and from \repr{conjdicmax}, there is a single conjugacy class of subgroups isomorphic to $\quat[32]$ or to $\dic{24}$. For the remaining cases, the number of conjugacy classes is given in \reco{clasq8} and 
\repr{conjclass345} that make use of~\cite{BCG}, and is equal to one for $\oonestar$ and $\tonestar$, to three for $\quat[16]$, to two for $\dic{12}$, and to five for $\quat$. We will construct explicit representatives of these conjugacy classes later in this section.

\subsection{Centralisers of the finite order elements of $B_{4}(\rp)$}\label{sec:centralB4}

In this section, we will compute the centraliser of most of the elements of finite order of $B_{4}(\rp)$. Up to conjugacy, it suffices to consider powers of the five representatives given in \resec{genconsid}. First note that if $1 \leq k \leq 7$ (resp.\ $k\in \brak{1,3,5}$), \repr{centraliserx}(\ref{it:centraliser}) implies that the centraliser in $B_4(\rp)$ of $a^{k}$ (resp.\ of $b^{k}$) is equal to $\ang{a}$ (resp.\ to $\ang{b}$). By \reth{ccrp2}(\ref{it:ccrp2b}), $b^{2}$ is conjugate to $\sigma_1\sigma_2\sigma_3\sigma_1$, and from \reco{ccorder4}(\ref{it:ccorder4b}), there are two conjugacy classes of subgroups isomorphic to $\Z_{8}$, one of which is represented by $\ang{a^{2}}$, and the other by $\ang{\sigma_1\sigma_2\sigma_3}$ (the two generators are non conjugate). Also, since $b^{4}=\ft[4] b^{-2}$, the centraliser of $b^{4}$ is equal to that of $b^{2}$. So it remains to determine the centralisers of the following five elements: $b^{2}$ (order $6$), $\sigma_1\sigma_2\sigma_3$ (order $8$), and three elements of order $4$, $(\sigma_1\sigma_2\sigma_3)^2$, $(\rho_2\sigma_1)^2\sigma_3$ and $\rho_1 \sigma_2\sigma_3 \sigma_2$. In \repr{conj4} we will show that the centraliser of $b^{2}$ is equal to $\ang{b}$, and in \repr{centeltsB4}, we will see that the centralisers of $(\rho_2\sigma_1)^2\sigma_3$ and $\rho_1\sigma_2\sigma_3 \sigma_2$ are both isomorphic to $\Z \oplus \Z_4$.

\begin{prop}\label{prop:conj4}
In $B_{4}(\rp)$, $Z_{B_{4}(\rp)}(b^{2})= \ang{b}$.
\end{prop}

\begin{proof}
Clearly $\ang{b} \subset Z_{B_{4}(\rp)}(b^{2})$. Conversely, let $u\in Z_{B_{4}(\rp)}(b^{2})$. Note that $b=\rho_{3} \sigma_{2}\sigma_{1}$ belongs to the mixed braid subgroup:
\begin{equation*}
B_{3,1}(\rp)=\setr{\beta \in B_{4}(\rp)}{\pi(\beta)(4)=4}.
\end{equation*}
Since $\pi(b)=(1,2,3)$, we have $\pi(u)\in \ang{(1,3,2)}$, and thus $Z_{B_{4}(\rp)}(b^{2})\subset B_{3,1}(\rp)$.  We have the following short exact sequence:
\begin{equation}\label{eq:sesb31}
\begin{tikzcd}
1\arrow{r} & \pi_{1}(\rp\setminus\brak{x_{1},x_{2},x_{3}}) \arrow{r} & B_{3,1}(\rp) \arrow{r}{f} & B_{3}(\rp) \arrow{r} & 1,
\end{tikzcd}
\end{equation} 
where $f$ is the homomorphism given geometrically by forgetting the last string. Then $f(b)=a$ in $B_{3}(\rp)$ and $f(u)\in Z_{B_{3}(\rp)}(a^{2})$. Now $Z_{B_{3}(\rp)}(a^{2})=\ang{a}$ by \repr{centregarside}(\ref{it:centrea2}), so it follows from~\reqref{sesb31} that $u=vb^{j}$, where $v\in \ker{f}\cap Z_{B_{4}(\rp)}(b^{2})$ and $j\in \brak{0,1,\ldots,11}$. It thus suffices to prove that $v=1$. The group $\ker{f}$ is free of rank $3$ with basis $\brak{A_{1,4},A_{2,4},\rho_{4}}$. By~\reqref{rn2} and~\reqref{defaij}, we have $\rho_{4}^{-2}= A_{1,4}A_{2,4} A_{3,4}$. Further, by \relem{revrhocomm}(\ref{it:revrhocommb}), we have $\rho_{3} A_{3,4}\rho_{3}^{-1}= \rho_{4}^{-1} A_{3,4}^{-1}\rho_{4}$, and using \repr{present}, we obtain: 
\begin{equation*}
b\sigma_{3}b^{-1} = \rho_{3} \sigma_{2}\sigma_{1} \sigma_{3} \sigma_{1}^{-1} \sigma_{2}^{-1} \rho_{3}^{-1}= \rho_{3} \sigma_{3}^{-1} \sigma_{2} \sigma_{3}  \rho_{3}^{-1}= \sigma_{3} \rho_{4} \sigma_{2} \rho_{4}^{-1} \sigma_{3}^{-1}= \sigma_{3} \sigma_{2}  \sigma_{3}^{-1}.
\end{equation*}
Hence: 
\begin{align*}
b\rho_{4}b^{-1} &= \rho_{3} \sigma_{2}\sigma_{1} \rho_{4} \sigma_{1}^{-1} \sigma_{2}^{-1} \rho_{3}^{-1}=\rho_{3} \rho_{4} \rho_{3}^{-1}=A_{3,4}\rho_{4}= A_{2,4}^{-1} A_{1,4}^{-1} \rho_{4}^{-1}\\
b A_{1,4} b^{-1} &= \rho_{3} \sigma_{2}\sigma_{1} \sigma_{1}^{-1} \sigma_{2}^{-1} \sigma_{3}^{2} \sigma_{2}\sigma_{1}\sigma_{1}^{-1} \sigma_{2}^{-1} \rho_{3}^{-1}= \rho_{3} A_{3,4}\rho_{3}^{-1}=\rho_{4}^{-1} A_{3,4}^{-1}\rho_{4}= \rho_{4} A_{1,4} A_{2,4} \rho_{4}\\
b A_{2,4} b^{-1} &= b \sigma_{3} \sigma_{2}^{2} \sigma_{3}^{-1}  b^{-1}= \sigma_{3} \sigma_{2}  \sigma_{3}^{-1} \sigma_{1}^{2} \sigma_{3} \sigma_{2}^{-1}  \sigma_{3}^{-1}=A_{1,4},
\end{align*}
where we have also used Propositions~\ref{prop:present} and~\ref{prop:permute}(\ref{it:conjbsigma}). It follows that $b^{2}\rho_{4}b^{-2}= A_{1,4}^{-1}\rho_{4}^{-1}$, $b^{2}A_{1,4} b^{-2}= \rho_{4} A_{1,4} A_{2,4}^{-1} A_{1,4}^{-1} \rho_{4}^{-1}$, hence $b^{4}\rho_{4}b^{-4}= \rho_{4} A_{1,4} A_{2,4}$. Then $\setl{b^{2i} \rho_{4} b^{-2i}}{i=0,1,2}$ is a basis of $\ker{f}$ whose elements are permuted cyclically by the action of conjugation by $b^{2}$, from which we conclude that $\ker{f}\cap Z_{B_{4}(\rp)}(b^{2})$ is trivial. So $v=1$ and $u\in \ang{b}$, and it follows that $Z_{B_{4}(\rp)}(b^{2})= \ang{b}$.
\end{proof}

\begin{rem}
In $B_{4}(\rp)$, one may check using the surface relation $\rho_{3} \sigma_{2}\sigma_{1}^{2}\sigma_{2} \rho_{3}=\sigma_{3}^{2}$ that $\rho_{2}^{-1} b^{2} \rho_{2}=\sigma_{3}^{2} \sigma_{1} \sigma_{2}$, and then that $\lambda b^{2} \lambda^{-1}=\sigma_{1} \sigma_{2} \sigma_{3}\sigma_{1}$, where $\lambda=\garside[4] \sigma_{1}^{-1} \sigma_{3}^{-1}\rho_{2}^{-1}$, from which it follows using \repr{conj4} that $Z_{B_{4}(\rp)}(\sigma_1\sigma_2\sigma_3\sigma_1)= \lambda\ang{b} \lambda^{-1}$. 
\end{rem}


   %

We now compute the centralisers of $(\rho_2\sigma_1)^2 \sigma_3$ and $\rho_1\sigma_2\sigma_3\sigma_2$. Note that they are infinite since the two elements are products of two commuting blocks each of which is of infinite order, and as we now see, they are generated by the `obvious' elements with respect to these blocks. 

\begin{prop}\label{prop:centeltsB4}\mbox{}
\begin{enumerate}
\item\label{it:centr1s2} The centraliser $Z_{B_{4}(\rp)}(\rho_{1} \sigma_{2}\sigma_{3}\sigma_{2})$ of $\rho_{1} \sigma_{2}\sigma_{3}\sigma_{2}$ is equal to $\ang{\rho_{1},\sigma_{2}\sigma_{3}\sigma_{2}}$, and is isomorphic to $\Z\oplus \Z_{4}$.  
\item\label{it:rsr} The centraliser $Z_{B_{4}(\rp)}((\rho_{2} \sigma_{1})^{2} \sigma_{3})$ of $(\rho_{2} \sigma_{1})^{2} \sigma_{3}$ is equal to $\ang{\rho_{2} \sigma_{1}, \sigma_{3}}$, and is isomorphic to $\Z\oplus \Z_{4}$.
 \end{enumerate}
\end{prop}

\begin{proof}\mbox{}
Let $u'=u=\rho_{1}$ and $v=\sigma_{2}\sigma_{3} \sigma_{2}$ (resp.\ $u'=(\rho_{2} \sigma_{1})^{2}$, $u=\rho_{2} \sigma_{1}$ and $v=\sigma_{3}$). Since $u$ and $v$ commute, it is clear that $Z_{B_{4}(\rp)}(u'v)\supset \ang{u,v}$. Conversely, let $z\in Z_{B_{4}(\rp)}(u'v)$. Then $\pi(z)\in Z_{\sn[4]}((2,4))$ (resp.\ $\pi(z)\in Z_{\sn[4]}((3,4))$), so $\pi(z)\in \ang{(1,3), (2,4)}$ (resp.\ $\pi(z)\in \ang{(1,2),(3,4)}$). We analyse these two cases as follows:
\begin{enumerate}
\item assume that we are in the first case, that of $\rho_{1} \sigma_{2}\sigma_{3}\sigma_{2}$, so $\pi(z)\in \ang{(1,3), (2,4)}$. Then $z,uv\in B_{2,2}(\rp)$, where the first (resp.\ second) index of the mixed braid group $B_{2,2}(\rp)$ corresponds to the first and third (resp.\ second and fourth) strings. Consider the projection $\map{g}{B_{2,2}(\rp)}[B_{2}(\rp)]$ given geometrically by forgetting the second and fourth strings. Then $g(u'v)= g(\rho_{1} \sigma_{2} \sigma_{3} \sigma_{2})=\rho_{1}$, and thus $g(z)\in Z_{B_{2}(\rp)}(\rho_{1})$. Now $\rho_{1}^{2}= \sigma_{1}^{2} = \ft[2]$ by~\reqref{surfacerp2}, so $\rho_{1}$ is of order $4$ in $B_{2}(\rp)$, and since $q(\ang{a})=\ang{(\overline{1}, \overline{1})}$ and $q(\rho_1)=(\overline{0}, \overline{1})$, we see that $\rho_{1} \in B_{2}(\rp) \setminus \ang{a}$. Now $B_2(\rp)=\ang{a,\garside[2]}\cong \quat[16]$, and it follows from standard properties of dicyclic groups that $Z_{B_{2}(\rp)}(\rho_{1})= \ang{\rho_{1}}\subset P_{2}(\rp)$, and so $g(z)\in P_{2}(\rp)$, whence $\pi(z) \in \ang{(2,4)}$. 

\item suppose that we are in the second case, that of $(\rho_{2} \sigma_{1})^{2} \sigma_{3}$, so $\pi(z)\in \ang{(1,2), (3,4)}$. Now $\rho_{2} \sigma_{1}\in Z_{B_{4}(\rp)}(u'v)$, and multiplying $z$ by this element if necessary, we may suppose that $\pi(z)\in \ang{(3,4)}$.
\end{enumerate}
We now return to the analysis of the two cases simultaneously. We have just shown that $\pi(z) \in \ang{(2,4)}$ (resp.\ $\pi(z) \in \ang{(3,4)}$). Then $z\in B_{2,1,1}(\rp)$, where the first index corresponds to the second and fourth strings (resp.\ to the third and fourth strings). Consider the projection $\map{h}{B_{2,1,1}(\rp)}[B_{2,1}(\rp)]$ given geometrically by forgetting the first string. Then $h(u'v)= w$, where $w=\garside[3]$ (resp.\ $w= \rho_{1}\sigma_{2}$) (for the second case, we have used the fact that $(\rho_{2} \sigma_{1})^{2}= \rho_{2}\rho_{1}$). Since $h(z)\in Z_{B_3(\rp)}(w)$, it follows from \repr{centregarside}(\ref{it:centregar3}) (resp.\ \repr{centralr1s2b3}) that $h(z)\in K$, where $K=\ang{\garside[3]}$ (resp.\ $K=\ang{\rho_{1}, \sigma_{2}}$). Let $L=\ang{\rho_{1} \sigma_{2}\sigma_{3}\sigma_{2}}$ (resp.\ $L=\ang{(\rho_{2}\sigma_{1})^{2}, \sigma_{3}}$). Since $h(\rho_{1} \sigma_{2}\sigma_{3}\sigma_{2})=\garside[3]$ (resp.\ $h((\rho_{2}\sigma_{1})^{2})= \rho_{1}$ and $h(\sigma_{3})= \sigma_{2}$), there exist $k\in \ker{h}$ and $l\in L$ such that $z=kl$. Now $\ker{h}$ may be identified with 
$P_{1}(\rp\setminus\brak{x_{2},x_{3},x_{4}},x_{1})$, which is a free group of rank $3$ with basis $\brak{\rho_{1}, A_{1,2}, A_{1,3}}$. The action by conjugation of certain elements of $B_{4}(\rp)$ on these basis elements is given as follows: 
\begin{center}\renewcommand{\arraystretch}{1.4}
\begin{tabular}{|>{$}c<{$}||*{4}{>{$}c<{$}|}}\hline
& \sigma_{2} & \sigma_{3} & \rho_{1} & \rho_{2}\\ \hline\hline
\rho_{1}  & \rho_{1}  & \rho_{1}  & \rho_{1}  & A_{1,2}^{-1} \rho_{1}\\ \hline
A_{1,2} & A_{1,3} & A_{1,2} &  \rho_{1} A_{1,2} \rho_{1}^{-1} & \rho_{1}^{-1}  A_{1,2}^{-1} \rho_{1} \\ \hline
A_{1,3} & A_{1,3}^{-1}A_{1,2}A_{1,3} & A_{1,3}^{-1}A_{1,2}^{-1} \rho_{1}^{2} & \rho_{1} A_{1,3} \rho_{1}^{-1} & \rho_{1}^{-1}  A_{1,2} \rho_{1}  A_{1,3} \rho_{1}^{-1}  A_{1,2}^{-1} \rho_{1}\\ \hline
\end{tabular}
\end{center}
The action of $\sigma_{3}$ on $A_{1,3}$ is obtained using~\reqref{defaij} and the surface relation $\rho_{1}^{2}= A_{1,2} A_{1,3} A_{1,4}$, 
and the action of $\rho_2$ (resp.\ of $\sigma_2$) on $\rho_1$, $A_{1,2}$ and $A_{1,3}$ follows from parts~(\ref{it:revrhocomma}), (\ref{it:revrhocommb}) and~(\ref{it:revrhocommc}) respectively of \relem{revrhocomm} (resp.\ from~\reqref{defaij}). 
Using this table and the fact that $(\rho_2 \sigma_1)^2=\rho_2 \rho_1$, we obtain the action of $\rho_{1}\sigma_{2}\sigma_{3}\sigma_{2}$ and $(\rho_{2} \sigma_{1})^{2} \sigma_{3}$ on the basis elements as follows: 
\begin{center}\renewcommand{\arraystretch}{1.4}
\begin{tabular}{|>{$}c<{$}||*{2}{>{$}c<{$}|}}\hline
& \rho_{1}\sigma_{2}\sigma_{3}\sigma_{2} & (\rho_{2} \sigma_{1})^{2} \sigma_{3}\\ \hline\hline
\rho_{1}  & \rho_{1} & A_{1,2}^{-1} \rho_{1}\\ \hline
A_{1,2} & \rho_{1} A_{1,3}^{-1} A_{1,2} ^{-1}\rho_{1} 
& A_{1,2}^{-1}
\\ \hline
A_{1,3} & \rho_{1}^{-1} A_{1,2} A_{1,3}A_{1,2}^{-1} \rho_{1} 
& \rho_{1} A_{1,3}^{-1} A_{1,2}^{-1}\rho_{1}
\\ \hline
\end{tabular}
\end{center}
Now  $\brak{\rho_{1},A_{1,2},\rho_{1} A_{1,3}^{-1} A_{1,2} ^{-1}\rho_{1}^{-1}}$ (resp.\ $\brak{\rho_{1}, A_{1,2}^{-1} \rho_{1}, A_{1,3}^{-1} \rho_{1}}$) is a basis of $\ker{p}$, and the action by conjugation by $\rho_{1}\sigma_{2} \sigma_{3} \sigma_{2}$ (resp.\ $(\rho_{2} \sigma_{1})^{2} \sigma_{3}$)  induces an automorphism of $\ker{h}$ that fixes $\rho_{1}$ and exchanges $A_{1,2}$ and $\rho_{1} A_{1,3}^{-1} A_{1,2} ^{-1}\rho_{1}$ (resp.\ that exchanges $\rho_{1}$ and $A_{1,2}^{-1} \rho_{1}$ and sends $\rho_{1} A_{1,3}^{-1} $ to $(\rho_{1} A_{1,3}^{-1})^{-1}$). In the first case, the elements of $\ker{h}$ that commute with $\rho_{1}\sigma_{2} \sigma_{3} \sigma_{2}$ are the elements of $\ang{\rho_{1}}$, and it follows from above that $z\in \ang{\rho_{1}, \rho_{1}\sigma_{2} \sigma_{3} \sigma_{2}}$, but then $z\in \ang{\rho_{1},\sigma_{2} \sigma_{3} \sigma_{2}}$, and so $Z_{B_{2,1,1}(\rp)}(\rho_{1} \sigma_{2}\sigma_{3}\sigma_{2}) \subset \ang{\rho_{1},\sigma_{2} \sigma_{3} \sigma_{2}}$. In the second case, the automorphism has no non-trivial fixed points, , and we conclude from above that $z\in \ang{(\rho_{2} \sigma_{1})^{2}, \sigma_{3}}$. This shows that $Z_{B_{4}(\rp)}(u'v)\subset \ang{u,v}$, which proves the first part of the statements of~(\ref{it:centr1s2}) and~(\ref{it:rsr}).

It remains to prove that both of the given centralisers are isomorphic to $\Z\oplus \Z_{4}$. It follows from the first part of the proof that $Z_{B_{4}(\rp)}(u'v)= \ang{u,v}=\ang{u'v,u}$. Since $u'v$ is of order $4$ and commutes with $u$, and $u$ is of infinite order, we see that the map $\map{\psi}{\Z\oplus \Z_{4}}[Z_{B_{4}(\rp)}(u'v)]$ given by $\psi(1,\overline{0})=u$ and $\psi(0,\overline{1})=u'v$ is a surjective homomorphism, and that the quotient $Q$ of $Z_{B_{4}(\rp)}(u'v)=\ang{u'v,u}$ by the subgroup $\ang{u'v}$ is isomorphic to $\Z$. Let $\map{\alpha}{Z_{B_{4}(\rp)}(u'v)}[Q]$ denote the corresponding canonical projection, let $\map{\overline{\psi}}{\Z_{4}}[\ang{u'v}]$ be the isomorphism that is the restriction of $\psi$ to the $\Z_{4}$-factor, and let $\map{\beta}{\Z}[Q]$ denote the isomorphism given by $\beta(1)=\alpha(u)$. We thus obtain the following commutative diagram of short exact sequences:
\begin{equation*}\label{eq:diag}
\begin{tikzcd}
1 \arrow[r] & \Z_4 \arrow[r] \arrow[d, "\overline{\psi}"] & \Z\oplus \Z_{4} \arrow[r] \arrow[d, "\psi"] & \Z \arrow[d, "\beta"] \arrow[r] & 1\\
1 \arrow[r] & \ang{u'v}   \arrow[r] & Z_{B_{4}(\rp)}(u'v) \arrow[r, "\alpha"] & Q \arrow[r] & 1.
\end{tikzcd}
\end{equation*}
Since $\overline{\psi}$ and $\beta$ are isomorphisms, it follows that $\psi$ is too, and this completes the proof of the proposition.
\end{proof}

With respect to the discussion at the beginning of this section, it remains to determine the centralisers of $\sigma_1\sigma_2\sigma_3$ and of $(\sigma_1\sigma_2\sigma_3)^2$. That of the former element is more difficult since there exists no projection onto a (mixed) subgroup on a smaller number of strings, and that of the latter element also seems to be complicated. The presentation of $B_n(\rp)$ given in \reth{newpresrp2} may be useful in this respect.

\subsection{Conjugacy classes of the dicyclic subgroups of $B_{4}(\rp)$}\label{sec:q8b4}


From the discussion in \resec{genconsid}, the isomorphism classes of the dicyclic subgroups of $B_{4}(\rp)$ are given by $\quat[32], \dic{24}, \quat[16], \dic{12}$ and $\quat$, and there is a single conjugacy class of the first two. We also know from~\cite{BCG} that there are three conjugacy classes for $\quat[16]$, to two for $\dic{12}$, and to five for $\quat$. We now construct explicit representatives of these conjugacy classes.

Recall that the standard copies of $\quat[32]$ and $\dic{24}$ are given by $\ang{a,\garside[4]}$ and $\ang{b,\garside[4] a^{-1}}$ respectively, from which we obtain three copies  $\ang{a^{4},\garside[4]}$, $\ang{a^{4}, a\garside[4]}$ and $\ang{b^{3},\garside[4] a^{-1}}$ of $\quat$. Let $\alpha_{0}= \sigma_1\sigma_2\sigma_3$. By~\reqref{muragen}, $\alpha_{0}$ is the Murasugi element $\delta_{0,0}$ of order $8$, and using the constructions of \relem{conjinv} and \repr{nomaxcyclic}, we see that $\ang{\alpha_{0}, a^{4}\garside[4]}$ is isomorphic to $\quat[16]$. It possesses two subgroups $\ang{\alpha_{0}^{2}, a^{4}\garside[4]}$ and $\ang{\alpha_{0}^{2}, \alpha_{0} a^{4}\garside[4]}$ that are isomorphic to $\quat$. By studying the images of these five copies of $\quat$ under Abelianisation and the permutation homomorphism, we may see that they are pairwise non conjugate, and so are representatives of the five conjugacy classes of $\quat$ in $B_{4}(\rp)$. We sum up this information in the following table.
\begin{center}\renewcommand{\arraystretch}{1.4}
\begin{tabular}{|>{$}c<{$}||*{2}{>{$}c<{$}|}}\hline
\text{copy of $\quat$} & \text{Abelianisation} & \text{permutation subgroup}\\ \hline\hline
\quat^{(1)}= \ang{\alpha_{0}^2, \alpha_{0} a^4\garside[4]} & \ang{(\overline{1}, \overline{0})} & \ang{(13)(24), (24)}\\ \hline
\quat^{(2)}= \ang{\alpha_{0}^2, a^4\garside[4]} & \ang{(\overline{0}, \overline{0})} & \ang{(13)(24), (14)(23)}\\ \hline
\quat^{(3)}= \ang{b^3,\garside[4]a^{-1}} & \Z_2\oplus \Z_2 & \ang{(24)}\\ \hline
\quat^{(4)}= \ang{a^4,a\garside[4]} & \ang{(\overline{1}, \overline{1})} & \ang{(13)}\\ \hline
\quat^{(5)}=\ang{a^4,\garside[4]} & \ang{(\overline{0}, \overline{0})} & \ang{(14)(23)}\\ \hline
\end{tabular}
\end{center}
For future reference, we label these subgroups as in the table.

\begin{rems}\mbox{}
\begin{enumerate}
\item During our investigations, we found two other copies $\quat^{(6)}= \ang{\rho_2\rho_1\sigma_3,a^2\garside[4]}$ and $\quat^{(7)}=\ang{\rho_2\rho_1\sigma_3, \sigma_{3}a^2\garside[4]}$ of $\quat$, but they are conjugate to $\quat^{(4)}$ because they have the same invariants. One may check that $\lambda \quat^{(4)} \lambda^{-1}= \quat^{(6)}$ and $\lambda' \quat^{(6)} \lambda'^{-1}= \quat^{(7)}$, where $\lambda=\sigma_1^{-1} a \sigma_1^{-1}$ and $\lambda'=\sigma_1^{-1} \rho_{2}^{-1}$.
\item It follows from the constructions of these copies of $\quat$ that none of them is maximal.
\end{enumerate}
\end{rems}

We now analyse the conjugacy classes of $\quat[16]$ in $B_{4}(\rp)$. By 
\repr{conjclass345} of 
\ref{app:quatbnrp2}, there are three such classes, and they are represented by the two subgroups $\ang{a^2,\garside[4]}$ and $\ang{a^2,a\garside[4]}$ of the standard copy $\ang{a,\garside[4]}$ of $\quat[32]$, and by the third copy $\ang{\alpha_{0}, a^{4}\garside[4]}$ mentioned above. Once more, they may be seen to be pairwise non conjugate by studying their images under Abelianisation and the permutation homomorphism. We sum up this information in the following table, and again we label the copies for future use.
\begin{center}\renewcommand{\arraystretch}{1.4}
\begin{tabular}{|>{$}c<{$}||*{2}{>{$}c<{$}|}}\hline
\text{copy of $\quat[16]$} & \text{Abelianisation} & \text{permutation subgroup}\\ \hline\hline
\quat[16]^{(1)}= \ang{\alpha_{0}, a^{4}\garside[4]} & \ang{(\overline{1}, \overline{0})} & \sn[4]\\ \hline
\quat[16]^{(2)}= \ang{a^2,\garside[4]} & \ang{(\overline{0}, \overline{0})} & \ang{(13)(24), (14)(23)}\\ \hline
\quat[16]^{(3)}= \ang{a^2,a\garside[4]} & \ang{(\overline{1}, \overline{1})} & \ang{(13)(24), (13)}\\ \hline
\end{tabular}
\end{center}
By construction, the first two copies are non maximal. We will see in \resec{latticeB4} that the third copy is also non maximal.

Finally, by 
\repr{conjclass345} of~\ref{app:quatbnrp2}, there are two conjugacy classes of subgroups of $B_{4}(\rp)$ isomorphic to $\dic{12}$. They are represented by the (non conjugate) subgroups $\dic{12}^{(1)}= \ang{b^2,\garside[4]a^{-1}}$ and $\dic{12}^{(2)}= \ang{b^2, b\garside[4]a^{-1}}$ of the standard copy $\ang{b,\garside[4]a^{-1}}$ of $\dic{24}$.

\subsection{Realisation of the binary polyhedral groups in $B_{4}(\rp)$}\label{sec:polyB4}

By \reth{finitebnrp2}, $B_{4}(\rp)$ possesses subgroups isomorphic to $\tonestar$ and $\oonestar$, and from~\cite{BCG}, there is one conjugacy class of each. In this section, we give explicit constructions of these two subgroups.

\begin{prop}\label{prop:realTOB4}
In $B_{4}(\rp)$, we have:
\begin{enumerate}[(a)]
\item\label{it:consostarb4} $\ang{a^{2} \sigma_{3}\alpha_{0}^{-3}, \alpha_{0}^{2}, a^{4}\garside[4]}\cong \tonestar$.
\item $\ang{\alpha_{0}^{2}a^{4} \garside[4]\alpha_{0}^{-1}, a^{2} \sigma_{3}\alpha_{0}^{-3}, \alpha_{0}^{2}, a^{4}\garside[4]}\cong \oonestar$.
\end{enumerate}
\end{prop}

\begin{proof}\mbox{}
\begin{enumerate}[(a)]
\item Let $H=\ang{w, \alpha_{0}^{2}, a^{4}\garside[4]}$, where $w=a^{2} \sigma_{3}\alpha_{0}^{-3}$, and let $w'=a^{2} \sigma_{3}\alpha_{0}$. Note that the subgroup $\ang{\alpha_{0}^{2}, a^{4}\garside[4]}$ of $H$ is the subgroup $\quat^{(2)}$ of \resec{q8b4}. Conjugating the surface relation $\rho_{1}^{2}=\sigma_{1} \sigma_{2} \sigma_{3}^{2} \sigma_{2} \sigma_{1}$ by $\garside[4]$ and using~\reqref{conjsigma} and~\reqref{conjrho}, we have $\rho_{4}^{-2}= \sigma_{3}\sigma_{2}\sigma_{1}^{2}\sigma_{2}\sigma_{3}$, and thus:
\begin{align*}
w' &= \rho_{4} \sigma_{3}\sigma_{2}\sigma_{1}\rho_{4} \sigma_{3}\sigma_{2} \sigma_{1}^{2}\sigma_{2}\sigma_{3} \sigma_{2}= \rho_{4} \sigma_{3}\sigma_{2}\sigma_{1}\rho_{4}^{-1} \sigma_{2}= \rho_{4} \sigma_{3}\sigma_{1}\sigma_{2}\sigma_{1} \rho_{4}^{-1}\\
& = \rho_{4} \sigma_{3}(\sigma_{1}\sigma_{2}\sigma_{1} \sigma_{3})\sigma_{3}^{-1}\rho_{4}^{-1} = \rho_{4} \sigma_{3} \delta_{0,1} \sigma_{3}^{-1}\rho_{4}^{-1},
\end{align*}
which is of order $6$ by~\reqref{muragen}. Hence $w=a^{2} \sigma_{3}\alpha_{0}^{-3}= w'\ft[4]$ is of order $3$ since $\ft[4]=\alpha_{0}^{4}$, and $\alpha_{0}$ is of order $8$. Now:
\begin{equation*}
\sigma_{3} \alpha_{0}^{2} \sigma_{3}^{-1}=\sigma_{3} \sigma_{1} \sigma_{2} \sigma_{3} \sigma_{1} \sigma_{2} =\sigma_{1} \sigma_{2} \sigma_{3} \sigma_{2} \sigma_{1} \sigma_{2} = \sigma_{1} \sigma_{2} \sigma_{3} \sigma_{1} \sigma_{2} \sigma_{1}=\garside[4],
\end{equation*}
so 
\begin{align*}
w\alpha_{0}^{2}w^{-1}=a^{2} \sigma_{3} \alpha_{0}^{2} \sigma_{3}^{-1} a^{-2} =a^{2} \garside[4] a^{-2} =a^{4} \garside[4].
\end{align*}
by \req{conja}, and making use of~\reqref{conjsigma}--\reqref{conjrho} and \repr{permute}(\ref{it:permsi}), we obtain:
\begin{align*}
w a^{4} \garside[4] w^{-1} &= w'  a^{4} \garside[4] w'^{-1}= a^{2} \sigma_{3}\alpha_{0} a^{4} \garside[4] \alpha_{0}^{-1} \sigma_{3}^{-1} a^{-2} \garside[4]^{-1} \garside[4]\\
&= a^{2} \sigma_{3}\alpha_{0} a^{4} \sigma_{1}^{-1} \sigma_{2}^{-1}\sigma_{3}^{-1}\sigma_{1}^{-1} a^{-4}a^{6} \garside[4]= a^{2} \sigma_{3}\alpha_{0} \sigma_{1} \sigma_{2} \sigma_{3} \sigma_{1} a^{6} \garside[4] = a^{2} \sigma_{3}\alpha_{0}^{2} \sigma_{1} a^{6} \garside[4]\\
&= (a^{2} \sigma_{3}\alpha_{0}^{2}  \sigma_{3}^{-1} a^{-2}) a^{2} \sigma_{3}\sigma_{1} a^{6} \garside[4]= a^{4} \garside[4] a^{2} \sigma_{3}\sigma_{1} a^{6} \garside[4]=\garside[4] (a^{-2} \sigma_{3}\sigma_{1} a^{2}) a^{4} \garside[4]\\
&= \garside[4] \sigma_{1}^{-1}\sigma_{3} a^{4} \garside[4] = \sigma_{3}^{-1}\garside[4] \sigma_{3} a^{4} \garside[4]= \alpha_{0}^{2}\ldotp a^{4} \garside[4].
\end{align*}
Setting $X=w$, $P=\alpha_{0}^{2} a^{4}\garside[4]$ and $Q=\alpha_{0}^{2}$, we see that the relations of~\reqref{presostar} involving these three elements are satisfied, so $H$ is a quotient of $\tonestar$. Now $H$ contains $\quat^{(2)}$ and $w$, so it is of order at least $24$, and we conclude that $H\cong \tonestar$.

\item Let $K=\ang{\alpha_{0}^{2}a^{4} \garside[4]\alpha_{0}^{-1}, w, \alpha_{0}^{2}, a^{4}\garside[4]}$. With respect to the presentation~\reqref{presostar} of $\oonestar$, set $R=\alpha_{0}^{2}a^{4} \garside[4]\alpha_{0}^{-1}$. By part~(\ref{it:consostarb4}), it suffices to check the validity of the relations of that presentation involving $R$. Making use of~\reqref{conjsigma}--\reqref{conjrho} and \repr{permute}(\ref{it:permsi}), the fact that $\alpha_{0}^{-1}= \alpha_{0}^{5}= \alpha_{0} \ft[4]$ and the relation $a^{4} \garside[4]\alpha_{0} \garside[4]^{-1}a^{-4}= \alpha_{0}^{-1}$, we have:
\begin{align*}
R^{2} &= \alpha_{0}^{2}a^{4} \garside[4]\alpha_{0}^{-1}\alpha_{0}^{2}a^{4} \garside[4]\alpha_{0}^{-1}= \alpha_{0}^{2} a^{4}\sigma_{3} \sigma_{2} \sigma_{1} a^{-4} \alpha_{0}^{-1} \ft[4]=\ft[4]=P^{2}=Q^{2}\\
RXR^{-1} &= \alpha_{0}^{2}a^{4} \garside[4]\alpha_{0}^{-1} a^{2} \sigma_{3}\alpha_{0} \ft[4]\alpha_{0} \garside[4]^{-1} a^{-4} \alpha_{0}^{-2}= \alpha_{0}^{3} a^{-2}  \sigma_{1}^{-1} \alpha_{0}^{-2} \alpha_{0}^{-2} \ft[4]=\alpha_{0}^{3} \sigma_{3}^{-1} a^{-2}= X^{-1}\\
RPR^{-1} &=\alpha_{0}^{2}a^{4} \garside[4]\alpha_{0}^{-1} \alpha_{0}^{2} a^{4}\garside[4] \alpha_{0} \garside[4]^{-1} a^{-4} \alpha_{0}^{-2}=\alpha_{0} a^{4}\garside[4] \alpha_{0}^{-3} =\alpha_{0}^{2} \ldotp \alpha_{0}^{2} a^{4}\garside[4]= QP\\
RQR^{-1} &=\alpha_{0}^{2}a^{4} \garside[4]\alpha_{0}^{-1} \alpha_{0}^{2} \alpha_{0} \garside[4]^{-1} a^{-4} \alpha_{0}^{-2}=\alpha_{0}^{2}a^{4} \garside[4] \alpha_{0}^{2} \garside[4]^{-1} a^{-4} \alpha_{0}^{-2}=\alpha_{0}^{-2}= Q^{-1}.
\end{align*}
Thus $K$ is a quotient of $\oonestar$. But $K$ contains $H$, and $R^{-1}P=\alpha_{0}\in K$, which is of order $8$, and so $R^{-1}P \in K \setminus H$. It follows that $K$ is of order at least $48$, from which we deduce that $K\cong \oonestar$.\qedhere
\end{enumerate}
\end{proof}

\begin{rem}
Using the constructions of~\cite{GG5}, we may obtain explicit realisations of $\tonestar$ and $\oonestar$ in $B_{8}(\St)$.
\end{rem}

\subsection{The subgroup lattice of $B_{4}(\rp)$}\label{sec:latticeB4}

In a manner similar to that of \resec{latticeB3}, we are now able to determine the subgroup lattice of the finite subgroups of $B_{4}(\rp)$ up to conjugacy. By \reth{finitebnrp2}, the abstractly maximal subgroups of $B_{4}(\rp)$ are $\oonestar$, $\dic{24}$ and $\quat[32]$, and applying \repr{conjdicmax} and~\cite{BCG}, there is a single conjugacy class of each, represented by the subgroups $\ang{\alpha_{0}^{2}a^{4} \garside[4]\alpha_{0}^{-1}, a^{2} \sigma_{3}\alpha_{0}^{-3}, \alpha_{0}^{2}, a^{4}\garside[4]}$, $\ang{b, \garside[4] a^{-1}}$ and $\ang{a,\garside[4]}$ respectively. The subgroup lattices of $\oonestar$ and $\dic{24}$ in $B_{4}(\rp)$ have the same form as those of Figure~\ref{fig:omegai}, and we just need to check that the copies of the subgroups that appear in these lattices are correctly labelled with respect to the notation of \resec{subgroupsb4}. 

Recall that for $m\in \brak{3,6,12,16}$ (resp.\ $m=8$), there is a single conjugacy class (resp.\ two conjugacy classes) of subgroups isomorphic to $\Z_{m}$ in $B_{4}(\rp)$. With respect to \resec{genconsid}, let $\Z_{4}^{(1)}=\ang{(\rho_2\sigma_1)^2\sigma_3}$, $\Z_{4}^{(2)}= \ang{b^3}$, $\Z_{4}^{(3)}= \ang{\rho_1\sigma_2\sigma_3\sigma_2}$, $\Z_{4}^{(4)}=\ang{\alpha_{0}^2}$ and $\Z_{4}^{(5)}=\ang{a^4}$, let $\Z_{16}=\ang{a}$, $\Z_{12}=\ang{b}$, $\Z_{8}^{(1)}=\ang{\alpha_{0}}$, $\Z_{8}^{(2)}=\ang{a^{2}}$, $\Z_{6}=\ang{\sigma_1\sigma_2\sigma_3\sigma_1}$ and let $\Z_{3}= \ang{\sigma_1\sigma_2\sigma_3\sigma_1}$. 

By \repr{realTOB4} and \resec{q8b4}, we have $q(\oonestar)= q(\dic{12}^{(2)})= \ang{(\overline{1}, \overline{0})}$ and $q(\dic{12}^{(1)})= \ang{(\overline{1}, \overline{1})}$, so $\dic{12}^{(2)}$ is the copy of $\dic{12}$ lying in $\oonestar$. From the proof of \repr{realTOB4}, $\alpha_{0} \in \oonestar$, and since $\alpha_{0}^{2}a^{4} \garside[4]\alpha_{0}^{-1} \in \oonestar$, it follows that $a^{4} \garside[4]\in \oonestar$, and thus $\quat[16]^{(1)}$ is the copy of $\quat[16]$ lying in $\oonestar$. Since $\quat[16]^{(1)}= \ang{\alpha_{0}, a^{4}\garside[4]}$, we see that $\quat^{(1)}$ and $\quat^{(2)}$ are its subgroups isomorphic to $\quat$, and as we saw in the proof of \repr{realTOB4}, $\tonestar$ contains the copy $\quat^{(2)}$ of $\quat$. Clearly $\quat^{(1)}$, $\quat^{(2)}$ and $\Z_{8}^{(1)}$ all contain $\Z_{4}^{(4)}$. Since $q(b\garside[4] a^{-1})=q(\alpha_{0} a^{4} \garside[4])= q(\alpha_{0}^{3} a^{4} \garside[4])= (\overline{1}, \overline{0})$, it follows from \resec{genconsid} that $\dic{12}^{(2)}$ and $\quat^{(1)}$ contain a copy of $\Z_{4}^{(1)}$. So with the given notation for the finite subgroups of $B_{4}(\rp)$, the subgroup lattice of $\oonestar$ is that of Figure~\ref{fig:omegai}.

We now consider the subgroup lattice of $\dic{24}$. Since $\dic{24}= \ang{b, \garside[4] a^{-1}}$, its copy of $\quat$ is $\ang{b^3, \garside[4] a^{-1}}$, which is $\quat^{(3)}$. Both $\quat^{(3)}$ and $\Z_{12}$ contain $\Z_{4}^{(2)}$. Since $q(b^3 \garside[4] a^{-1})= (\overline{1}, \overline{0})$ (resp.\ $q(\garside[4] a^{-1})= (\overline{1}, \overline{1})$), the subgroup $\ang{b^3 \garside[4] a^{-1}}$ (resp.\ $\ang{\garside[4] a^{-1}}$) of $\quat^{(3)}$ is represented by $\Z_{4}^{(1)}$ (resp.\ by $\Z_{4}^{(3)}$). Once more, with the given notation for the finite subgroups of $B_{4}(\rp)$, the subgroup lattice of $\dic{24}$ is that of Figure~\ref{fig:omegai}.

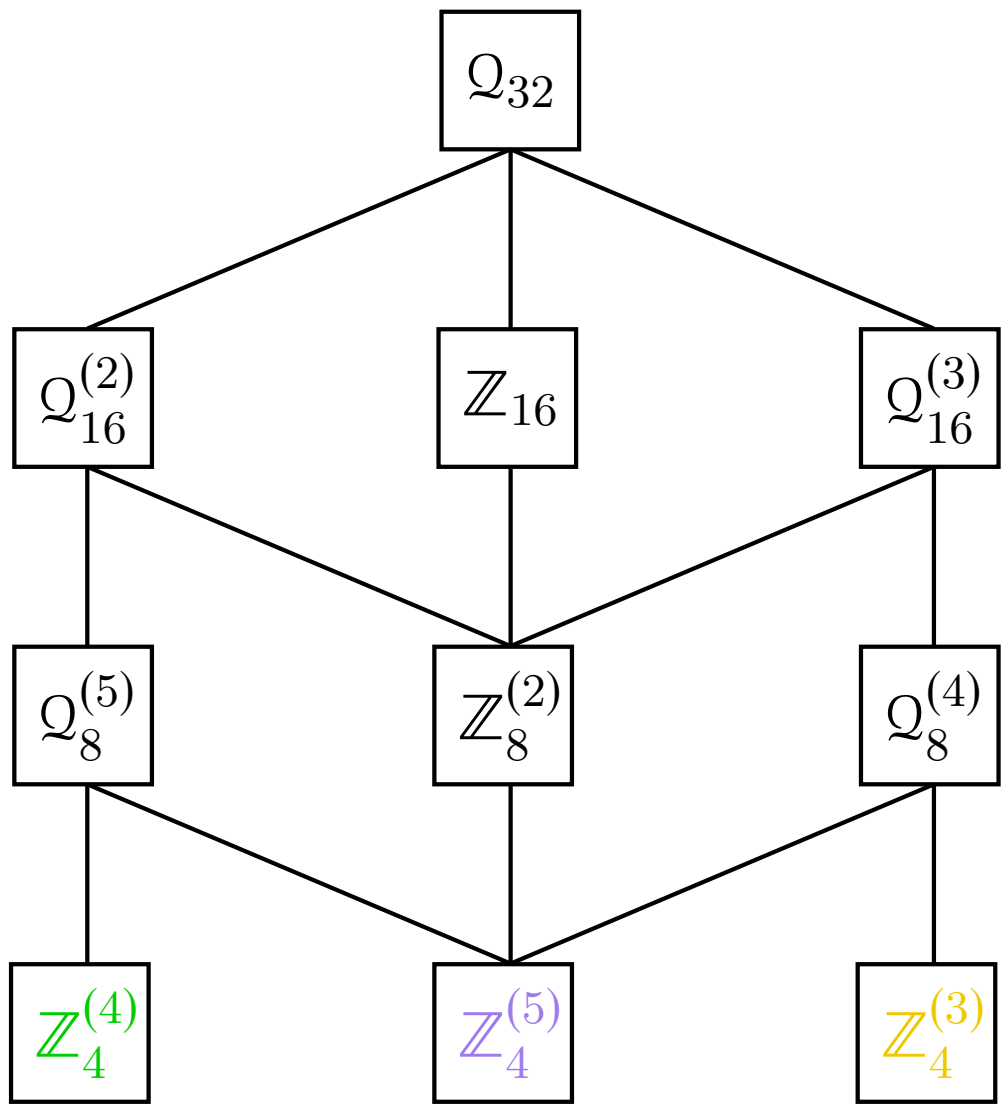
\begin{figure}[!h]
\hfill
\begin{tikzpicture}[scale=0.7]
\tikzstyle{every node} =[node contents={}, minimum height=0.9cm,minimum width=0.9cm]
\node (a) at (5,11) {$\quat[32]$};

\node (b) at (1,8){$\quat[16]^{(2)}$};
\node (c) at (5,8){$\Z_{16}$};
\node (d) at (9,8){$\quat[16]^{(3)}$};

\node (e) at (1,5){$\quat^{(5)}$};
\node (f) at (5,5){$\Z_{8}^{(2)}$};
\node (g) at (9,5){$\quat^{(4)}$};

\node (h) at (1,2){\textcolor{Green3}{$\Z_{4}^{(4)}$}};
\node (i) at (5,2){\textcolor{MediumPurple2}{$\Z_{4}^{(5)}$}};
\node (j) at (9,2){\textcolor{Gold2}{$\Z_{4}^{(3)}$}};

\foreach \toa in {b,c,d} \draw[thick] (a.south)--(\toa.north);
\foreach \tob in {e,f} \draw[thick] (b.south)--(\tob.north);
\foreach \tof in {c,d} \draw[thick] (\tof.south)--(f.north);
\foreach \tog in {d} \draw[thick] (\tog.south)--(g.north);
\foreach \toi in {e,f,g} \draw[thick] (\toi.south)--(i.north);

\foreach \toh in {e} \draw[thick] (\toh.south)--(h.north);
\foreach \toj in {g} \draw[thick] (\toj.south)--(j.north);

\foreach \boxx in {a,b,c,d,e,f,g,h,i,j} \draw[thick] (\boxx.south west) rectangle ++(1.3,1.3); 

\end{tikzpicture}\hspace*{\fill}
\caption{The lattice of the conjugacy classes of the subgroups of $\quat[32]$ in $B_{4}(\rp)$. 
}\label{fig:omegaii}
\end{figure}

Finally, the subgroup lattice of $\quat[32]=\ang{a,\garside[4]}$ is as shown in Figure~\ref{fig:omegaii}. The group $\ang{a,\garside[4]}$ has three subgroups $\quat[16]^{(2)}= \ang{a^{2},\garside[4]}$, $\quat[16]^{(3)}= \ang{a^{2},a\garside[4]}$ and $\Z_{16}=\ang{a}$ of order $16$, and the two subgroups $\quat^{(5)}= \ang{a^{4},\garside[4]}$ and $\ang{a^{4}, a^{2}\garside[4]}$ (resp.\ $\quat^{(4)}= \ang{a^{4}, a\garside[4]}$ and $\ang{a^{4}, a^{3}\garside[4]}$) of $\quat[16]^{(2)}$ (resp.\ of $\quat[16]^{(3)}$) isomorphic to $\quat$ are conjugate in $B_{4}(\rp)$. The group $\Z_{16}=\ang{a}$ contains $\Z_{8}^{(2)}=\ang{a^{2}}$, and the groups $\quat^{(4)}$, $\quat^{(5)}$ and $\Z_{8}^{(2)}$ contain $\Z_{4}^{(5)}=\ang{a^{4}}$. The other subgroups of $\quat^{(4)}$ (resp.\ $\quat^{(5)}$) of order $4$ are $\ang{a\garside[4]}$ and $\ang{a^{5}\garside[4]}$ (resp.\ $\ang{\garside[4]}$ and $\ang{a^{4}\garside[4]}$), so they are conjugate in $B_{4}(\rp)$, and using the results of \resec{genconsid}, are conjugate to $\Z_{4}^{(3)}$ (resp.\ to to $\Z_{4}^{(4)}$). Note that the groups $\Z_{4}^{(3)}$ and $\Z_{4}^{(4)}$ also appear in the lattice of $\dic{24}$ and $\oonestar$ respectively.
Observe that all of the conjugacy classes of the finite subgroups of $B_{4}(\rp)$ appear in at least one of Figures~\ref{fig:omegai} and~\ref{fig:omegaii}. This implies that no non abstractly maximal subgroup of $B_{4}(\rp)$ is maximal. As in Figure~\ref{fig:omegai}, the two groups $\Z_4^{(3)}$ and $\Z_4^{(4)}$ are in coloured text to indicate that they also appear in one of the lattices of that figure, and the group $\Z_4^{(5)}$ is also coloured for later use.


\section{The subgroup structure of $B_5(\rp)$}\label{sec:subgroupsB5}
 
By \reth{finitebnrp2}, the abstractly maximal finite subgroups of $B_{5}(\rp)$ are $\dic{40}$ and $\quat[32]$, and so the isomorphism classes of the non-trivial finite subgroups of $B_{5}(\rp)$ are $\dic{4l}$ for $l\in \brak{2,4,5,8,10}$ and $\Z_{m}$ for $m\in \brak{2,4,5,8,10,16,20}$. In particular, $B_{5}(\rp)$ has no binary polyhedral group as a finite subgroup. From \repr{conjdicmax}, there is a single conjugacy class of subgroups isomorphic to $\dic{40}$ or to $\quat[32]$ and these subgroups are maximal. By \reqref{muragen}, up to conjugacy and powers, the finite-order elements of $B_{5}(\rp)$ of order at least $4$ are as follows:
\begin{enumerate}[(a)]
\item $\sigma_1\sigma_2\sigma_3\sigma_4$ (order $10$)
\item $\rho_1(\sigma_2\sigma_3\sigma_4)^2$ (order $4$)
\item $(\rho_3\sigma_2\sigma_1)^3\sigma_4$ (order $4$)
\item $a=\rho_5\sigma_4\sigma_3\sigma_2\sigma_1$ (order $20$)
\item $\sigma_1\sigma_2\sigma_3\sigma_4\sigma_1$ (order $8$)
\item $(\rho_2\sigma_1)^2\sigma_3\sigma_4\sigma_3$ (order $4$)
\item $b=\rho_4\sigma_3\sigma_2\sigma_1$ (order $16$)
\end{enumerate}
As in the previous sections, we label the copies of the representatives for future use. Taking appropriate powers of the above elements, we obtain the six (pairwise non conjugate) subgroups of order $4$ given in the following table.
\begin{center}\renewcommand{\arraystretch}{1.4}
\begin{tabular}{|>{$}c<{$}||*{2}{>{$}c<{$}|}}\hline
\text{copy of $\Z_{4}$} & \text{Abelianisation} & \text{permutation subgroup}\\ \hline\hline
\Z_{4}^{(1)}= \ang{a^5} & \ang{(\overline{0}, \overline{1})} & \id\\ \hline
\Z_{4}^{(2)}= \ang{(\sigma_1\sigma_2\sigma_3\sigma_4\sigma_1)^2} & \ang{(\overline{0}, \overline{0})} & \ang{(14)(35)}\\ \hline
\Z_{4}^{(3)}= \ang{(\rho_2\sigma_1)^2\sigma_3\sigma_4\sigma_3} & \ang{(\overline{1}, \overline{0})} & \ang{(35)}\\ \hline
\Z_{4}^{(4)}= \ang{\rho_1(\sigma_2\sigma_3\sigma_4)^2} & \ang{(\overline{0}, \overline{1})} & \ang{(24)(35)}\\ \hline
\Z_{4}^{(5)}= \ang{b^4} & \ang{(\overline{0}, \overline{0})} & \id\\ \hline
\Z_{4}^{(6)}= \ang{(\rho_3\sigma_2\sigma_1)^3\sigma_4} & \ang{(\overline{1}, \overline{1})} & \ang{(45)}\\ \hline
\end{tabular}
\end{center}
By \reco{ccorder4}(\ref{it:ccorder4b}), there are two conjugacy classes of subgroups isomorphic to $\Z_{8}$, $\Z_{8}^{(1)}=\ang{\sigma_1\sigma_2\sigma_3\sigma_4\sigma_1}$ and $\Z_{8}^{(2)}= \ang{b^{2}}$, and one conjugacy class of subgroups isomorphic to $\Z_{m}$ for $m\in \brak{16, 20}$, so $\Z_{16}=\ang{b}$ and $\Z_{20}=\ang{a}$,  and by \reth{ccrp2}(\ref{it:ccrp2b}), there is one conjugacy class of subgroups isomorphic to $\Z_{m}$ for $m\in \brak{5,10}$. By \reco{clasq8} and
\repr{conjclass345},
$B_{5}(\rp)$ has three conjugacy classes for $\quat[16]$, two conjugacy classes for $\dic{20}$, and six conjugacy classes for $\quat$. For $\dic{20}$, we may take the representatives to be the two subgroups $\dic{20}^{(1)}= \ang{a^2,\garside[5]}$ and $\dic{20}^{(2)}= \ang{a^2,a\garside[5]}$ of the standard copy $\ang{a,\garside[5]}$ of $\dic{40}$. In the same way, for $\quat[16]$, two of the representatives may be taken to be the subgroups $\quat[16]^{(2)}= \ang{b^2,\garside[5]a^{-1}}$, $\quat[16]^{(3)}= \ang{b^2,b\garside[5]a^{-1}}$ of the standard copy $\ang{b,\garside[5]a^{-1}}$ of $\quat[32]$. A third copy may be obtained by making use of the above element $\sigma_1\sigma_2\sigma_3\sigma_4\sigma_1$ of order $8$. Noting that it is the Murasugi element $\delta_{0,1}$ of $B_{5}(\rp)$, taking $y=\sigma_{2}\sigma_{3} a^{5} \garside[5]$ in the proof of \relem{conjinv} and following the proof of \repr{nomaxcyclic}, we see that $\quat[16]^{(1)}= \ang{\delta_{0,1}, y}\cong \quat[16]$. These three copies are non conjugate, since $q(\quat[16]^{(2)})=\ang{(\overline{0}, \overline{1})}$, $q(\quat[16]^{(3)})=\ang{(\overline{1}, \overline{0})}$ and $q(\quat[16]^{(1)})\cong \Z_{2}\oplus \Z_{2}$. 

To obtain representatives of five of the six conjugacy classes of $\quat$, we may make use of the subgroups of $\ang{\delta_{0,1}, y}$ and of the standard copies of $\quat[32]$ and $\dic{40}$ that are isomorphic to $\quat$.  The sixth copy may be constructed as follows. Consider the Murasugi element $\delta_{2,1}= (\rho_2\sigma_1)^2\sigma_3\sigma_4\sigma_3=\rho_2\rho_1\sigma_3\sigma_4\sigma_3$, which is of order $4$ from, and the subgroup $\ang{\delta_{2,1}, a^3\garside[5]^{-1}}$. By~\reqref{conja}, $(a^3\garside[5]^{-1})^{2}=\ft[5]$, so $a^3\garside[5]^{-1}$ is of order $4$. Further, $q(\delta_{2,1})=(\overline{1}, \overline{0})$ and $q(a^3\garside[5]^{-1})= (\overline{0}, \overline{1})$, so $a^3\garside[5]^{-1} \notin \ang{\delta_{2,1}}$, and $a^3\garside[5]^{-1} \delta_{2,1} (a^3\garside[5]^{-1})^{-1}= \delta_{2,1}^{-1}$ using~\reqref{conjsigma},~\reqref{conjrho}, \repr{permute}(\ref{it:permsi}) and~(\ref{it:permrhi}), and the Artin relations. We conclude that $\ang{\delta_{2,1}, a^3\garside[5]^{-1}}\cong \quat$. 
We thus have the following copies of $\quat$ in $B_{5}(\rp)$:
\begin{enumerate}[(a)]
\item $\quat^{{(1)}}= \ang{\delta_{0,1}^{2}, \delta_{0,1} y}$  (Abelianisation $\ang{(\overline{1}, \overline{1})}$, permutation  subgroup $\ang{(14)(35), (35)}$).
\item $\quat^{{(2)}}= \ang{a^5,\garside[5]}$ (Abelianisation $\ang{(\overline{0}, \overline{1})}$, permutation subgroup $\ang{(15)(24)}$).
\item $\quat^{{(3)}}= \ang{\delta_{0,1}^{2}, y}$ (Abelianisation $\ang{(\overline{0}, \overline{1})}$, permutation subgroup $\ang{(14)(35), (15)(34)}$).
\item $\quat^{{(4)}}= \ang{b^4,b\garside[5]a^{-1}}$ (Abelianisation $\ang{(\overline{1}, \overline{0})}$, permutation subgroup $\ang{(13)}$).
\item $\quat^{{(5)}}= \ang{b^4,\garside[5]a^{-1}}$ (Abelianisation $\ang{(\overline{0}, \overline{1})}$, permutation subgroup $\ang{(14)(23)}$).
\item $\quat^{{(6)}}= \ang{\delta_{2,1} ,a^3\garside[5]^{-1}}$ (Abelianisation $\Z_2\oplus \Z_2$, permutation subgroup $\ang{(35), (12)(35)}$).
\end{enumerate}
The given invariants allow us to conclude that these six copies of $\quat$ are pairwise non conjugate with the possible exception of 
$\quat^{{(2)}}$ and $\quat^{{(5)}}$ that have the same invariants. 
Suppose on the contrary that these two subgroups are conjugate in $B_{5}(\rp)$. Now $\quat^{{(2)}}\cap P_{5}(\rp)= \ang{a^{5}}$ and $\quat^{{(5)}}\cap P_{5}(\rp)= \ang{b^{4}}$, and since $P_{5}(\rp)$ is normal in $B_{5}(\rp)$, these two intersections would be conjugate to, but this is not possible because $q(a^{5})= (\overline{0}, \overline{1})$ and $q(b^{4})= (\overline{0}, \overline{0})$. So the above six copies are representatives of the six conjugacy classes of $\quat$ in $B_{5}(\rp)$. 

\begin{figure}[!h]
\hfill
\begin{tikzpicture}[scale=0.7]
\tikzstyle{every node} =[node contents={}, minimum height=0.9cm,minimum width=0.9cm]
\node (aa) at (18.5,11) {$\dic{40}$};
\node (ab) at (15.5,8) {$\dic{20}^{(1)}$};
\node (ac) at (12,6) {$\quat^{(2)}$};
\node (ad) at (18.5,8) {$\Z_{20}$};
\node (ae) at (21.5,8) {$\dic{20}^{(2)}$};
\node (af) at (18.5,5) {$\Z_{10}$};
\node (ag) at (18.5,2) {$\Z_{5}$};
\node (ah) at (21.5,2) {\textcolor{Green3}{$\Z_{4}^{(4)}$}}; 
\node (ai) at (15.5,2) {$\Z_{4}^{(1)}$}; 
\node (aj) at (12,2) {\textcolor{Orange1}{$\Z_{4}^{(2)}$}}; 
\foreach \tob in {ab,ad,ae} \draw[thick] (aa.south)--(\tob.north);
\foreach \tob in {aj,af} \draw[thick] (ab.south)--(\tob.north);
\foreach \tob in {ag} \draw[thick] (af.south)--(\tob.north);

\path[thick] (aa.south) edge [out=-165, in=90] (ac.north);

\draw[line width=5pt,white] ($0.9*(ac.south) + 0.1*(ah.north)$)--($0.1*(ac.south) + 0.9*(ah.north)$);

\foreach \tob in {ah} \draw[thick] (ac.south)--(\tob.north);
\foreach \tob in {ai} \path[line width=5pt,white] (ad.south) edge [out=-165, in=90] (\tob.north);
\foreach \tob in {ah,af} \draw[thick] (ae.south)--(\tob.north);

\draw[line width=5pt,white] ($0.9*(ac.south) + 0.1*(ai.north)$)--($0.1*(ac.south) + 0.9*(ai.north)$);
\foreach \tob in {ai} \draw[thick] (ac.south)--(\tob.north);
\foreach \tob in {ai} \draw[thick] (ad.south) edge [out=-165, in=90] (\tob.north);
\foreach \tob in {aj} \draw[thick] (ac.south)--(\tob.north);
\foreach \tob in {af} \draw[thick] (ad.south)--(\tob.north);
\foreach \boxx in {aa,ab,ae} \draw[thick] (\boxx.south west) rectangle ++(1.8,1.3); 
\foreach \boxx in {ac,ad,af,ag,ah,ai,aj} \draw[thick] (\boxx.south west) rectangle ++(1.3,1.3); 

\node (a) at (5,11) {$\quat[16]^{(1)}$};

\node (b) at (1,8){$\quat^{(3)}$};
\node (c) at (5,8){$\Z_{8}^{(1)}$};
\node (d) at (9,8){$\quat^{(1)}$};

\node (e) at (1,5){\textcolor{Green3}{$\Z_{4}^{(4)}$}};
\node (f) at (5,5){\textcolor{Orange1}{$\Z_{4}^{(2)}$}};
\node (g) at (9,5){\textcolor{VioletRed2}{$\Z_{4}^{(6)}$}};

\foreach \toa in {b,c,d} \draw[thick] (a.south)--(\toa.north);
\foreach \tob in {e,f} \draw[thick] (b.south)--(\tob.north);
\foreach \tof in {c,d} \draw[thick] (\tof.south)--(f.north);
\foreach \tog in {d} \draw[thick] (\tog.south)--(g.north);
\foreach \boxx in {a,b,c,d,e,f,g} \draw[thick] (\boxx.south west) rectangle ++(1.3,1.3); 

\node (ba) at (5,3) {$\quat^{(6)}$};

\node (bb) at (1,0.5){\textcolor{Gold2}{$\Z_{4}^{(3)}$}};
\node (bc) at (5,0.5){\textcolor{Green3}{$\Z_{4}^{(4)}$}};
\node (bd) at (9,0.5){\textcolor{VioletRed2}{$\Z_{4}^{(6)}$}};

\foreach \toba in {bb,bc,bd} \draw[thick] (ba.south)--(\toba.north);

\foreach \boxx in {ba,bb,bc,bd} \draw[thick] (\boxx.south west) rectangle ++(1.3,1.3); 

\end{tikzpicture}\hspace*{\fill}
\caption{The lattices of the conjugacy classes of $\quat[16]^{(1)}$, $\quat^{(6)}$ and $\dic{40}$ in $B_{5}(\rp)$. 
}\label{fig:omegaiii}
\end{figure}
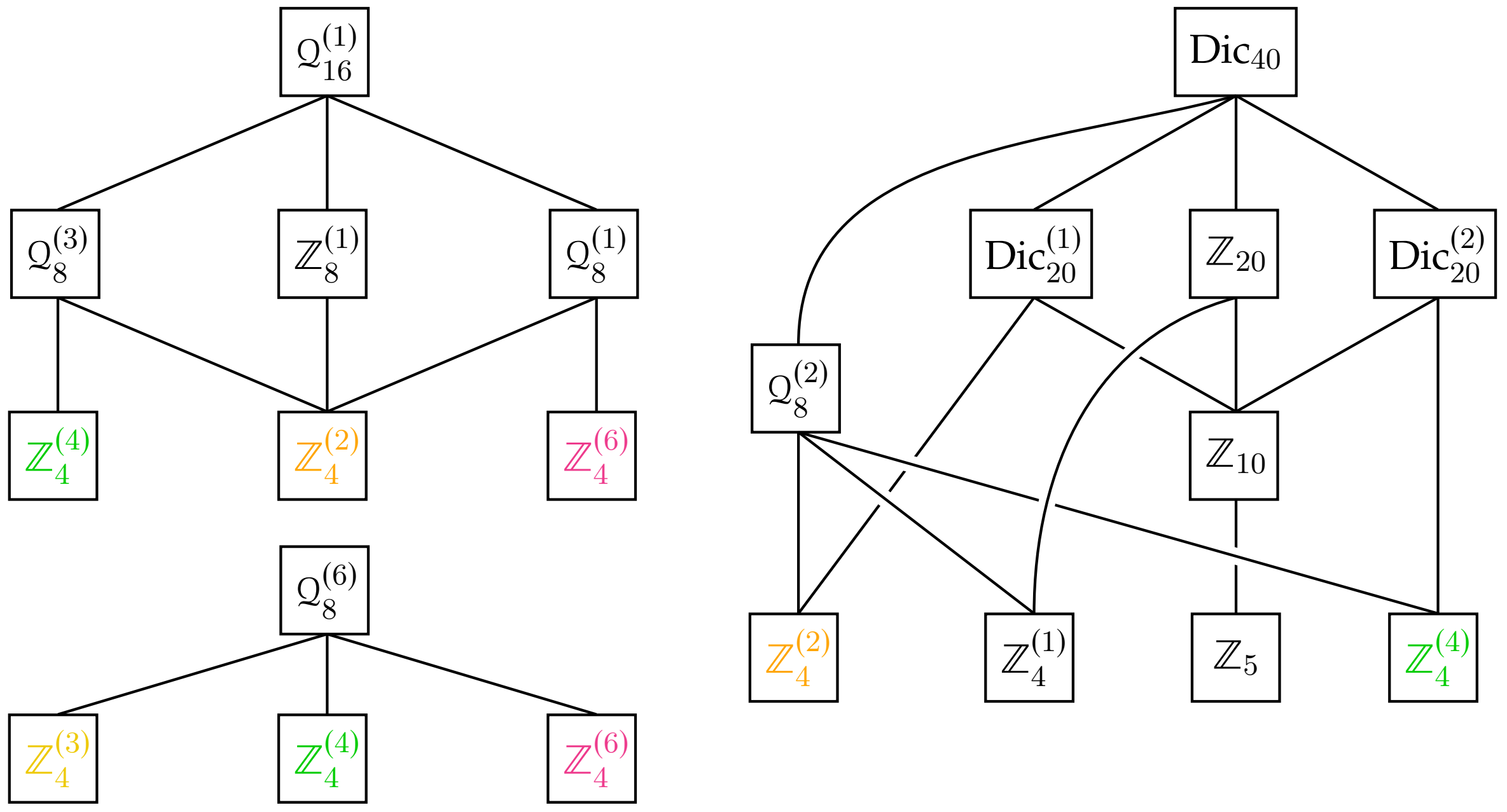

To determine the subgroup lattice for $B_{5}(\rp)$, we first analyse that of the two abstractly maximal subgroups $\quat[32]$ and $\dic{40}$. The lattice of $\quat[32]=\ang{b, \garside[5]a^{-1}}$ is given by Figure~\ref{fig:omegaii} where the labels of the subgroups correspond to those in this section. By the above table of the cyclic subgroups of order $4$, one may check that $\ang{\garside[5]a^{-1}}$ (resp.\ $\ang{b\garside[5]a^{-1}}$) is indeed conjugate to $\Z_{4}^{(4)}$ (resp.\ to $\Z_{4}^{(3)}$). The subgroup lattice of $\dic{40}$ is given in Figure~\ref{fig:omegaiii}. In the same manner, one may check that $\ang{\garside[5]}$ (resp.\ $\ang{a^{5}\garside[5]}$) is conjugate to $\Z_{4}^{(2)}$ (resp.\ to $\Z_{4}^{(4)}$).
One copy of $\quat[16]$, $\quat[16]^{(1)}=\ang{\delta_{0,1}, y}$, does not appear in either of the subgroup lattices of $\quat[32]$ and $\dic{40}$. Clearly its subgroups isomorphic to $\quat$ are $\quat^{(1)}$ and $\quat^{(3)}$, and the subgroup isomorphic to $\Z_{8}$ is $\Z_{8}^{(1)}$. These three subgroups all contain $\ang{\delta_{0,1}}$, which is $\Z_{4}^{(2)}$, and since $q(y)= (\overline{0}, \overline{1})$ and $\pi(y)=(15)(24)$ (resp.\ $q(\delta_{0,1}y)= (\overline{1}, \overline{1})$), it follows that $\ang{y}$ (resp.\ $\ang{\delta_{0,1}y}$) is conjugate to $\Z_{4}^{(2)}$ (resp.\ to $\Z_{4}^{(6)}$). The subgroup lattice of $\quat[16]^{(1)}$ is illustrated in Figure~\ref{fig:omegaiii}. One copy of $\quat$, $\quat^{{(6)}}= \ang{\delta_{2,1} ,a^3\garside[5]^{-1}}$, does not appear in any of the subgroup lattices of $\quat[32]$, $\dic{40}$ or $\quat[16]^{(1)}$. Since $q(\delta_{2,1})= (\overline{1}, \overline{0})$, $q(a^3\garside[5]^{-1})= (\overline{0}, \overline{1})$, $\pi(a^3\garside[5]^{-1})=(12)(35)$ and $q(\delta_{2,1} a^3\garside[5]^{-1})=(\overline{1}, \overline{1})$, it follows that $\ang{a^3\garside[5]^{-1}}$ is conjugate to $\Z_{4}^{(3)}$, $\ang{\delta_{2,1}}$ is conjugate to $\Z_{4}^{(4)}$, and $\ang{\delta_{2,1}a^3\garside[5]^{-1}}$ is conjugate to $\Z_{4}^{(6)}$. The subgroup lattice of $\quat^{{(6)}}$ is also shown in Figure~\ref{fig:omegaiii}. It follows that $\quat[16]^{(1)}$ and $\quat^{(6)}$ are both maximal but not abstractly maximal. Once more, the groups in coloured text appear in at least two of the four lattices given in Figures~\ref{fig:omegaii} and~\ref{fig:omegaiii}.

\begin{rem}
Note that $\quat[16]^{(1)}$ is the copy of $\quat[16]$ that appears in the proof of \repr{nomax}. The results of Sections~\ref{sec:subgroupsb3}--\ref{sec:subgroupsB5} show that $n=5$ is the smallest value of $n$ for which there exist maximal finite subgroups of $B_{n}(\rp)$ that are not abstractly maximal.
\end{rem}

\section{A new presentation of $B_{n}(\rp)$}\label{sec:newpres}


To finish the main part of this paper, in this section, we give a presentation of $B_{n}(\rp)$ for which the generators $\rho_{i}$, $i=1,\ldots,n$ of \repr{present} are replaced by the generators $\tau_{i}$ that appear in the presentation of $P_{n}(\rp)$ given in~\cite[Theorem~4]{GGdedicata} (see also~\cite[Remark~9(c)]{GGfibre}). Recall that for $i=1,\ldots,n$, $\tau_{i}$ is defined by:
\begin{equation}\label{eq:deftau}
\tau_{i}=\rho_{i}^{-1} A_{i,i+1}\cdots A_{i,n},
\end{equation}
where for $1\leq i< j \leq n$, the elements $A_{i,j}$ are defined in~\reqref{defaij}.


\begin{thm}\label{th:newpresrp2}
Let $n\in \N$.The following constitutes a presentation of the group $B_n(\rp)$:
\begin{enumerate}
\item[\underline{\textbf{generators:}}] $\si{1},\ldots,\si{n-1},\tau_{1},\ldots,\tau_{n}$.

\item[\underline{\textbf{relations:}}]
\begin{align}
\si{i}\si{j} &=\si{j}\si{i}\quad\text{if $\lvert i-j\rvert\geq 2$}\label{eq:newpres1}\\
\si{i}\si{i+1}\si{i}&=\si{i+1}\si{i}\si{i+1} \quad\text{for $1\leq i\leq n-2$}\label{eq:newpres2}\\
\si{i}\tau_{j}&=\tau_{j}\si{i}\quad\text{for $j\neq i,i+1$}\label{eq:newpres3}\\
\tau_{i+1}&=\si{i}\tau{_i}\sigma_{i}^{-1} \quad\text{for $1\leq i\leq n-1$}\label{eq:newpres4}\\
[\tau_{i+1}, \tau_{i}^{-1}]&= \si[2]{i} \quad\text{for $1\leq i\leq n-1$}\label{eq:newpres5}\\
\tau_{1}^{2} &=\si{1}\si{2}\cdots\si{n-2}\si[2]{n-1} \si{n-2}\cdots\si{2}\si{1}.\label{eq:newpres6}
\end{align}
\end{enumerate}
\end{thm}

\begin{rem}
The presentation of \reth{newpresrp2} is similar to that of Van Buskirk, with the difference that the $\tau_{i}$ are pairwise conjugate in $B_{n}(\rp)$, which is not the case for the $\rho_{i}$.
\end{rem}

\begin{proof}[Proof of \reth{newpresrp2}]
Let $1\leq i<j \leq n$. By~\reqref{defaijalt}, note that: 
\begin{equation}\label{eq:prodAij}
A_{i,i+1}\cdots A_{i,j}=\sigma_{i}\cdots \sigma_{j-2} \sigma_{j-1}^{2} \sigma_{j-2} \cdots \sigma_{i}.
\end{equation}
Taking $j=n$, we see from~\reqref{deftau} that $\rho_{i}=\sigma_{i}\cdots \sigma_{n-2} \sigma_{n-1}^{2} \sigma_{n-2} \cdots \sigma_{i}\tau_{i}^{-1}$, and it follows from \repr{present} that the elements $\si{1},\ldots,\si{n-1},\tau_{1},\ldots,\tau_{n}$ generate $B_{n}(\rp)$. To obtain the relations of \reth{newpresrp2}, we substitute the $\tau_{i}$ into those of~\repr{present}. The two Artin relations~\reqref{newpres1}--\reqref{newpres2} are clear, and relation~\reqref{newpres3} is a consequence of~\reqref{sirj} and~\reqref{newpres1}--\reqref{newpres2}. 
For the surface relation~\reqref{newpres6}, by~\req{surfacerp2},~\reqref{deftau} and~\reqref{prodAij}, we have:
\begin{equation*}
\tau_{1}= \rho_{1}^{-1} A_{1,2}\cdots A_{1,n}=  \rho_{1}^{-1} \sigma_{1}\cdots \sigma_{n-2} \sigma_{n-1}^{2} \sigma_{n-2} \cdots \sigma_{1}=\rho_{1},
\end{equation*}
so $\tau_{1}^{2}=\rho_{1}^{2}= \sigma_{1}\cdots \sigma_{n-2} \sigma_{n-1}^{2} \sigma_{n-2} \cdots \sigma_{1}$.
To obtain~\reqref{newpres4}, for $1\leq i\leq n-1$, by~\req{sirisi},~\reqref{deftau} and~\reqref{prodAij}, we have:
\begin{align*}
\tau_{i+1}&= \rho_{i+1}^{-1} A_{i+1,i+2}\cdots A_{i+1,n}= \rho_{i+1}^{-1}  \sigma_{i+1}\cdots \sigma_{j-2} \sigma_{j-1}^{2} \sigma_{j-2} \cdots \sigma_{i+1}\\
&= \rho_{i+1}^{-1} \sigma_{i}^{-1} \ldotp \sigma_{i} \sigma_{i+1}\cdots \sigma_{j-2} \sigma_{j-1}^{2} \sigma_{j-2} \cdots \sigma_{i+1}\sigma_{i} \ldotp \sigma_{i}^{-1}= \rho_{i+1}^{-1} \sigma_{i}^{-1} \rho_{i} \tau_{i} \sigma_{i}^{-1}= \sigma_{i} \tau_{i} \sigma_{i}^{-1}.
\end{align*}
Finally, let $1\leq i\leq n-1$. From \req{rhocomm}, ~\reqref{deftau} and~\reqref{prodAij}, we have:
\begin{equation*}
\sigma_{i}^{2}= [\rho_{i+1}^{-1}, \rho_{i}^{-1}]= [\tau_{i+1} \sigma_{i+1}^{-1}\cdots \sigma_{n-1}^{-1}\sigma_{n-1}^{-1} \cdots \sigma_{i+1}^{-1}, \tau_{i} \sigma_{i}^{-1}\cdots \sigma_{n-1}^{-1}\sigma_{n-1}^{-1} \cdots \sigma_{i}^{-1}].
\end{equation*}
Now $\sigma_{i+1}^{-1}\cdots \sigma_{n-1}^{-1}\sigma_{n-1}^{-1} \cdots \sigma_{i+1}^{-1}$ commutes with both $\tau_{i}$ and $ \sigma_{i}^{-1}\cdots \sigma_{n-1}^{-1}\sigma_{n-1}^{-1} \cdots \sigma_{i}^{-1}$, so:
\begin{align*}
\sigma_{i}^{2}&= \tau_{i+1} \tau_{i} \sigma_{i}^{-1}\cdots \sigma_{n-1}^{-1}\sigma_{n-1}^{-1} \cdots \sigma_{i}^{-1} \tau_{i+1}^{-1} \sigma_{i} \cdots \sigma_{n-1} \sigma_{n-1} \cdots \sigma_{i} \tau_{i}^{-1}\\
&= \tau_{i+1} \tau_{i} \sigma_{i}^{-1}\cdots \sigma_{n-1}^{-1}\sigma_{n-1}^{-1} \cdots \sigma_{i+1}^{-1} \tau_{i}^{-1} \sigma_{i+1} \cdots \sigma_{n-1} \sigma_{n-1} \cdots \sigma_{i} \tau_{i}^{-1}\\
&= \tau_{i+1} \tau_{i} \sigma_{i}^{-1} \tau_{i}^{-1} \sigma_{i} \tau_{i}^{-1}.
\end{align*}
It follows from this that:
\begin{equation*}
1=\sigma_{i}^{-1} \tau_{i+1} \sigma_{i} \sigma_{i}^{-2} \sigma_{i} \tau_{i} \sigma_{i}^{-1} \tau_{i}^{-1} \sigma_{i} \tau_{i}^{-1} \sigma_{i}^{-1} = \tau_{i} \sigma_{i}^{-2} \tau_{i+1} \tau_{i}^{-1}\tau_{i+1}^{-1},
\end{equation*}
and we conclude that $\sigma_{i}^{2}= [\tau_{i+1}, \tau_{i}^{-1}]$ as required. 
\end{proof}

This presentation has some other nice properties when compared with that of \repr{present}.

\begin{prop}\label{prop:goodpropsnew} 
Let $n\geq 2$. Then:
\begin{enumerate}[(a)]
\item for all $1\leq i<j\leq n$, $[\tau_{j}, \tau_{i}^{-1}]=A_{i,j}$.
\item for all $i=1,\ldots, n$, $\tau_{i}^{2}= A_{1,i} \cdots A_{i-1,i} A_{i,i+1}\cdots A_{i,n}$.
\item\label{it:goodpropsnewc} Let $\alpha_{0}=\sigma_{1}\cdots \sigma_{n-1}$. Then for all $k=1,\ldots, n$, $\alpha_{0} \tau_{k} \alpha_{0}^{-1}=\tau_{k+1}$, where the index is taken modulo $n$.
\item The subgroup $\ang{\tau_{1},\ldots, \tau_{n}}$ is normal in $B_{n}(\rp)$ and it coincides with $P_{n}(\rp)$.
\end{enumerate}
\end{prop}


\begin{rem}
The result of \repr{goodpropsnew}(\ref{it:goodpropsnewc}) may prove to be useful in the study of the centraliser of $\alpha$ mentioned at the end of \resec{centralB4}.
\end{rem}

\begin{proof}[Proof of \repr{goodpropsnew}]\mbox{}
\begin{enumerate}[(a)]
\item Let $1\leq i<j\leq n$. Then by~\reqref{newpres5}, $[\tau_{j}, \tau_{j-1}^{-1}]=\sigma_{j-1}^{2}=A_{j-1,j}$ and from~\reqref{defaij},~\reqref{newpres3} and~\reqref{newpres5}, we have:
\begin{align*}
A_{i,j} &=\sigma_{i}^{-1} \cdots \sigma_{j-2}^{-1} \sigma_{j-1}^{2} \sigma_{j-2}\cdots \sigma_{i}=
\sigma_{i}^{-1} \cdots \sigma_{j-2}^{-1} [\tau_{j}, \tau_{j-1}^{-1}] \sigma_{j-2}\cdots \sigma_{i}\\
&=  [\tau_{j}, \sigma_{i}^{-1} \cdots \sigma_{j-2}^{-1}\tau_{j-1}^{-1}\sigma_{j-2}\cdots \sigma_{i}]=  [\tau_{j}, \tau_{i}^{-1}].
\end{align*}

\item First note that if $i=1$ then the result follows from~\reqref{newpres6} and~\reqref{prodAij}. Suppose by induction that it holds for some $i=1,\ldots, n-1$.
Then by~\reqref{newpres4},~\reqref{prodAij} and the fact that $\sigma_{i-1}\cdots \sigma_{2}\sigma_{1}^{2} \sigma_{2} \cdots \sigma_{i-1}=A_{1,i} \cdots A_{i-1,i}$, we have:
\begin{align*}
\tau_{i+1}^{2} &= \sigma_{i} \tau_{i}^{2} \sigma_{i}^{-1}= \sigma_{i} A_{1,i} \cdots A_{i-1,i} A_{i,i+1}\cdots A_{i,n} \sigma_{i}^{-1}\\
&= \sigma_{i}\sigma_{i-1}\cdots \sigma_{2}\sigma_{1}^{2} \sigma_{2} \cdots \sigma_{i-1} \sigma_{i} \ldotp\sigma_{i+1}\cdots \sigma_{n-2}\sigma_{n-1}^{2} \sigma_{n-2} \cdots \sigma_{i+1}\sigma_{i} \sigma_{i}^{-1}\\
&= A_{1,i+1} \cdots A_{i,i+1} \ldotp A_{i+1,i+2} \cdots A_{i+1,n},
\end{align*}
which yields the result.

\item First let $i=1,\ldots, n-1$. Then by~\reqref{newpres3} and~\reqref{newpres4}, we have:
\begin{align*}
\alpha_{0} \tau_{i} \alpha_{0}^{-1} &= \sigma_{1} \cdots \sigma_{i-1} \sigma_{i} \tau_{i} \sigma_{i}^{-1} \sigma_{i-1}^{-1} \cdots \sigma_{1}^{-1} =\tau_{i+1}. 
\end{align*}
Finally, since $\alpha_{0}^{n}= \ft$, we see that $\tau_{1} = \alpha_{0}^{n} \tau_{1} \alpha_{0}^{-n}= \alpha_{0} \tau_{n} \alpha_{0}^{-1}$.

\item 
Let $H=\ang{\tau_{1},\ldots, \tau_{n}}$. It suffices to show that $\sigma_{i} \tau_{j} \sigma_{i}^{-1} \in H$ for all $i=1,\ldots, n-1$ and $j=1,\ldots,n$. If $j\neq i+1$ then this follows from~\reqref{newpres3} and~\reqref{newpres4}. If $j=i+1$ then by~\reqref{newpres4} and~\reqref{newpres5}, we have $\sigma_{i} \tau_{i+1} \sigma_{i}^{-1} =\sigma_{i}^{2} \tau_{i} \sigma_{i}^{-2}= [\tau_{i+1}, \tau_{i}^{-1}] \tau_{i} [\tau_{i+1}, \tau_{i}^{-1}]^{-1} \in H$, so $H$ is normal in $B_{n}(\rp)$. Since $\tau_i\in P_{n}(\rp)$ for $i=1,2,...,n$ it follows that $H\subset P_{n}(\rp)$. One may see from the presentation of $B_n(\rp)$ given in~\reth{newpresrp2} that the quotient $B_{n}(\rp)/H$ is isomorphic to the symmetric group $\sn$, and then that the induced homomorphism $B_{n}(\rp)/H \to B_{n}(\rp)/P_{n}(\rp)$ between copies of $\sn$ is surjective, and so injective. It follows that $H=P_n(\rp)$.\qedhere
\end{enumerate}
\end{proof}

\appendix

\renewcommand{\thesection}{Appendix~\Roman{section}}

\section{Conjugacy classes of finite-order elements of $P_{3}(\rp)$}\label{app:conjclass}

Let $n\geq 2$. With the notation of \req{defab}, let $A=a^{n}$ and $B=b^{n-1}$ in $B_{n}(\rp)$. From~\cite[Proposition~26 and Remark~27]{GG2}, $A=\rho_{n}\cdots \rho_{1}$ and $B=\rho_{n-1}\cdots \rho_{1}$ and these two elements are both pure braids of order $4$. Recall from~\cite[Theorem~6(d)]{GGzeit} that $P_{n}(\rp)$ is (torsion) generated by the set:
\begin{equation*}
\brak{A, B, a^{-1} B a,\ldots, a^{-(n-2)} B a^{n-2}}
\end{equation*}
of elements that are of order $4$, and that this set is minimal with respect to the property of normally generating $P_{n}(\rp)$. In particular, the elements of this set are pairwise non conjugate, and so are representatives of $n$ distinct conjugacy classes in $P_{n}(\rp)$ of elements of order $4$, which provides a lower bound for
the number of conjugacy classes of elements of order $4$ in $P_{n}(\rp)$. In~\cite[Proposition 21(3)]{GG2}, we stated that the set of elements of $P_{3}(\rp)$ of order $4$ form three distinct conjugacy classes. However this is not correct. Unfortunately this result was used in the proof 
of~\cite[Proposition~17]{GG3}. In this appendix, we take the opportunity to correct the proofs of these two results.
 
We first introduce some notation. 
Recall from the proof of \repr{centregarside} that there is an isomorphism $P_{3}(\rp)\cong \F[2] \rtimes \quat$~\cite{vB}, where $\F[2]=\F[2](\rho_{3}, A_{2,3})$, and $\quat=\ang{x,y}$, where $x=b^{2}= \rho_{2}\rho_{1}$ and $y=\rho_{2}A_{1,2}\rho_{3}^{-1}$. Let $z=xy$, and consider the basis 
$(\alpha,\beta)$ of $\F[2]$, where $\alpha=A_{2,3} \rho_{3}$ and $\beta=\rho_{3}$. We are now able to give the statement of the correct version of~\cite[Proposition 21(3)]{GG2}, as well as a result regarding the number of conjugacy classes in $P_{3}(\rp)$ of cyclic subgroups of order $4$.

\begin{prop}\label{prop:5conjclasses}\mbox{}
\begin{enumerate}
\item\label{it:5conja} The set of elements of $P_{3}(\rp)$ of order $4$ form five distinct conjugacy classes. A set of representatives of these classes is given by $\brak{(1,x), (1,y),(1,z),( \alpha,x), ( \alpha,x^{-1})}$. 
\item\label{it:5conjb} The group $P_{3}(\rp)$ possesses four conjugacy classes of subgroups isomorphic to $\Z_{4}$. A set of representatives of these classes is given by $\brak{\ang{(1,x)}, \ang{(1,y)},\ang{(1,z)},\ang{( \alpha,x)}}$
\end{enumerate}
\end{prop}

\begin{proof}\mbox{}
\begin{enumerate}
\item In what follows, we identify $P_{3}(\rp)$ with $\F[2] \rtimes \quat$. Let $T=\brak{x,y,z}$. If $t\in T$, we start by recalling the action $\phi_{t}$ of conjugation by $t$ on $\F[2]$ that may be obtained using~\reqref{conjxy}:
\begin{equation}\label{eq:phiaction}
\phi_{x}:\left\{
\begin{aligned}
\alpha &\mapsto \alpha^{-1}\\
\beta &\mapsto \beta^{-1}
\end{aligned}
\right. 
\qquad
\phi_{y}:\left\{
\begin{aligned}
\alpha &\mapsto \beta\\
\beta &\mapsto \alpha
\end{aligned}
\right. \qquad
\phi_{z}:\left\{
\begin{aligned}
\alpha &\mapsto \beta^{-1}\\
\beta &\mapsto \alpha^{-1}.
\end{aligned}
\right.
\end{equation}
Since $t^{2}=\ft[3]$, we have $\phi_{t^{2}}=\id$ and $\phi_{t}= \phi_{t^{-1}}$. 

Let $(\gamma, s)$ be an element of $P_{3}(\rp)$ of order $4$. If $\gamma=1$ then $s\in \brak{x^{\pm 1}, y^{\pm 1}, z^{\pm 1}}$. Since the second factor is $\quat$, it is clear that $(1,t^{-1})$ is conjugate in $P_{3}(\rp)$ to $(1,t)$ for all $t\in T$. Conversely, if $t,t'\in T$, comparison of the second factor shows that $(1,t)$ and $(1,t')$ are conjugate in $P_{3}(\rp)$ if and only if $t=t'$, so the elements $(1,t)$, where $t\in T$, are of order $4$, they are pairwise non-conjugate in $P_{3}(\rp)$, and they are representatives of the conjugacy classes of the elements of the form $(1,s)$ of $P_{3}(\rp)$ of order $4$, where $s\in \quat$.

So suppose that $\gamma \neq 1$. Then $s\notin\ang{\ft[3]}$ for otherwise $(\gamma, s)^{2}=(\gamma^{2},1)$, and $(\gamma, s)$ could not be of finite order. Thus $s$ is of order $4$ and belongs to $\brak{x^{\pm 1}, y^{\pm 1}, z^{\pm 1}}$. Since:
\begin{equation*}
1=(\gamma, s)^4=\left( (\gamma, s)\ldotp (\gamma, s)\right)^2= (\gamma \ldotp\phi_{s}(\gamma),s^2)^2=((\gamma\ldotp\phi_{s}(\gamma))^2, 1).
\end{equation*}
and $\gamma\ldotp\phi_{s}(\gamma)\in \F[2]$, it follows that $\gamma\ldotp\phi_{s}(\gamma)=1$, so:
\begin{equation}\label{eq:formw}
\phi_{s}(\gamma)=\gamma^{-1}.
\end{equation}
Set $\gamma =a_{1} \cdots a_{m}$, where $m\in \N$, each $a_{i}$ is a non-zero power of $\alpha$ or $\beta$, and $\gamma$ is written in reduced form with no cancellation or combination between successive terms. If $m$ is even then set $u=a_{1}\cdots a_{m/2}$, $v=a_{m/2+1}\cdots a_{m}$ and $\epsilon=0$. If $m$ is odd, then set $a_{(m+1)/2}=z^{2r+\epsilon}$, where $z\in \brak{\alpha,\alpha^{-1},\beta, \beta^{-1}}$, $r\geq 0$ and $\epsilon\in \brak{0,1}$, and $u=a_{1}\cdots a_{(m-1)/2}z^{r}$ and $v=z^{r}a_{(m+3)/2}\cdots a_{m}$. In both cases, $\gamma= u z^{\epsilon}v$ is in reduced form (there may be combination of terms if $m$ is odd, but there is no cancellation). Equation~\reqref{formw} then implies that:
\begin{equation}\label{eq:phitw}
v^{-1}z^{-\epsilon}u^{-1}=\gamma^{-1}=\phi_{s}(\gamma)=\phi_{s}(u) \phi_{s}(z^{\epsilon})\phi_{s}(v).
\end{equation}
From the form of $\phi_{s}$ given by~\reqref{phiaction}, there is no cancellation in the right-hand term of~\reqref{phitw}, and since the set $\brak{\alpha,\alpha^{-1}, \beta, \beta^{-1}}$ is invariant under the action $\phi_{s}$, it follows that $\phi_{s}(u)=v^{-1}$ (which implies also that $\phi_{s}(v)=u^{-1}$) and $\phi_{s}(z^{\epsilon})= z^{-\epsilon}$. Setting $\rho=(u,1)$, we thus have:
\begin{align*}
\rho (z^{\epsilon},s) \rho^{-1} &=(u,1)\ldotp (z^{\epsilon},s)\ldotp (u^{-1},1)= (u z^{\epsilon},s)\ldotp (u^{-1},1)=(u z^{\epsilon}\phi_{s}(u^{-1}),s)\\
&= (uz^{\epsilon}v,s)=(\gamma, s),
\end{align*}
and so $(\gamma, s)$ is conjugate in $P_{3}(\rp)$ to $(z^{\epsilon},s)$. If $\epsilon=0$ then we already saw that there are three conjugacy classes of the form $(1,s)$, where $s\in \brak{x^{\pm 1}, y^{\pm 1}, z^{\pm 1}}$. So assume that $\epsilon=1$. Since $\phi_{s}(z^{\epsilon})= z^{-\epsilon}$,~\reqref{phiaction} implies that $s\in \brak{x,x^{-1}}$. 
So $(\gamma, s)$ is conjugate in $P_{3}(\rp)$ to one of the eight pairs of the form $(z,s)$, where $z\in \brak{\alpha,\alpha^{-1},\beta, \beta^{-1}}$ and $s\in \brak{x,x^{-1}}$. The fact that $(1,x)\ldotp (z,s) \ldotp (1,x^{-1})=(z^{-1},s)$ implies that $(\gamma, s)$ is conjugate in $P_{3}(\rp)$ to one of the four pairs of the form $(z,s)$, where $z\in \brak{\alpha,\beta}$ and $s\in \brak{x,x^{-1}}$. Further, $(1,y)\ldotp (\beta,s) \ldotp (1,y^{-1})=(\alpha,s^{-1})$, from which we conclude that $(\gamma, s)$ is conjugate to one of $(\alpha,x)$ or $(\alpha,x^{-1})$.

Let $\mathcal{Z}=\brak{(\alpha,x^{-1}), (1,x), (1,y), (1,z)}$. We claim that $(\alpha,x)$ is not conjugate to any element of $\mathcal{Z}$. Indeed, if $(\gamma',t)\in P_{3}(\rp)$ then let $\tau=(\gamma',t) \ldotp (\alpha,x)\ldotp (\gamma',t)^{-1}$. Then:
\begin{align*}
\tau &= (\gamma',t) \ldotp (\alpha,x) \ldotp (\phi_{t^{-1}}(\gamma'^{-1}),t^{-1})= (\gamma'\phi_{t}(\alpha) \phi_{txt^{-1}}(\gamma'^{-1}), txt^{-1}).
\end{align*}
If $t\in \ang{x^{2}}$ then $\tau= (\gamma' \alpha \gamma'^{-1},x)$, which is not equal to an element of $\mathcal{Z}$ for any choice of $\gamma'\in \F[2]$. So suppose that $t\notin \ang{x^{2}}$. We now project onto the Abelianisation $\Z^{2}$ of $\F[2]$. Letting $(m,n)\in \Z^{2}$ denote the Abelianisation of $\gamma'$, and using the fact that $txt^{-1} \in \brak{x,x^{-1}}$, the first factor of $\tau$ is sent to $2(m,n)+\nu$, where $\nu= (-1,0)$ if $t\in \brak{x,x^{-1}}$, $\nu= (0,1)$ if $t\in \brak{y,y^{-1}}$, and $\nu= (0,-1)$ if $t\in \brak{z,z^{-1}}$,
So $2(m,n)+\nu \neq (0,0)$, and thus $(\alpha,x)$ cannot be conjugate to $(1,t')$ for any $t'\in T$. If $\tau$ is equal to the remaining element $(\alpha,x^{-1})$ of $\mathcal{Z}$ then $2(m,n)+\nu = (1,0)$, so $t \in \brak{x,x^{-1}}$, but then the second factor of $\tau$ cannot be equal to $x^{-1}$, which yields a contradiction. So $(\alpha,x)$ is not conjugate to $(\alpha,x^{-1})$ in $P_{3}(\rp)$. This proves the claim. Repeating the first part of the above argument, we see that $(\alpha,x^{-1})$ is not conjugate to any element of $\brak{(1,x),(1,y), (1,z)}$. We thus conclude that $P_{3}(\rp)$ contains exactly five conjugacy classes of elements of order~$4$, these classes being represented by the elements of $\brak{(\alpha ,x)}\cup \mathcal{Z}$, which proves the first part of the proposition. 

\item It suffices to study the subgroups generated by the five elements of order $4$ given by part~(\ref{it:5conja}). First note that if two subgroups of order $4$ are conjugate then the generator of the first subgroup is conjugate either to a given generator of the second subgroup or to its inverse. Since $(1,t)^{-1}=(1,t^{-1})$, and $(1,t)$ and $(1,t^{-1})$ are conjugate in $P_{3}(\rp)$ for all $t\in T$, it follows from part~(\ref{it:5conja}) that the three subgroups of the form $\ang{(1,t)}$, where $t\in T$, are pairwise non conjugate. Further, since $(\alpha,x)^{-1}= (\alpha,x^{-1})$, we have $\ang{(\alpha,x)}= \ang{(\alpha,x^{-1})}$, and we see once more using part~(\ref{it:5conja}) that this subgroup is not conjugate in $P_{3}(\rp)$ to any of the subgroups of the form $\ang{(1,t)}$, where $t\in T$. 
We conclude that $P_{3}(\rp)$ possesses exactly four conjugacy classes of subgroups isomorphic to $\Z_{4}$.\qedhere
\end{enumerate}
\end{proof}

We now prove~\cite[Proposition~17]{GG3} in light of the corrected version of 
In fact, we can give two proofs: one that follows the methods of that of the original proposition, and the other that makes use of \repr{centraliserx} and \resec{ccb3rp}.
 
\begin{prop}\label{prop:p3noz4z}
The pure braid group $P_{3}(\rp)$ has no subgroup isomorphic to $\Z_{4}\times \Z$.
\end{prop}


\begin{proof}[First proof of \repr{p3noz4z}]
The first part of the proof follows that given in~\cite{GG3}. Suppose on the contrary that $P_{3}(\rp)$ possesses a subgroup $G$ isomorphic to $\Z_{4}\times \Z$, and let $H$ denote its subgroup that is isomorphic to $\Z_{4}$. Conjugating $G$ by an element of $P_{3}(\rp)$ if necessary, identifying $P_{3}(\rp)$ with $\F[2]\rtimes \quat$ and applying \repr{5conjclasses}(\ref{it:5conjb}), we may suppose that $H=\ang{(\rho,t)}$, where $t\in T$, $\rho=1$ if $t\in \brak{y,z}$ and $\rho\in \brak{1, \alpha}$ if $t= x$. Let $(u,\xi)$ generate the $\Z$-factor of $G$ for some $u\in \F[2]$, where $u\neq 1$, and $\xi\in \quat$. Since $G=H \times \ang{(u,\xi)}$ we have $(\rho,t)\ldotp (u,\xi) = (u,\xi)\ldotp (\rho,t)$. In particular, $t\xi=\xi t$, hence $\xi\in \ang{t}$ because $\xi,t \in \quat$, and $\rho \phi_{t}(u)=u\phi_{\xi}(\rho)$. If $\rho=1$ then $\phi_{t}(u)=u$, but this is impossible since by \req{phiaction}, the free group automorphism $\phi_{t}$ has no non-trivial fixed points. Hence $\rho= \alpha$ and $t=x$, and since $\xi\in \ang{x}$, we have $\phi_{\xi}(\alpha)= \alpha^{\epsilon}$ where $\epsilon=1$ if $\xi\in \ang{\ft[3]}$ and $\epsilon=-1$ if $\xi\in \ang{x}\setminus \ang{\ft[3]}$.
Setting $u= \alpha^{\epsilon_{1}} \beta^{\gamma_{1}}\cdots  \alpha^{\epsilon_{m}} \beta^{\gamma_{m}}$, where for $i=1,\ldots,m$, $\epsilon_{i}$ and $\gamma_{i}$ are non-zero integers, with the possible exception of $\epsilon_{1}$ and $\gamma_{m}$, the relation $ \alpha\phi_{x}(u)=u\phi_{\xi}( \alpha)$ gives rise to:
\begin{equation*}
 \alpha^{1-\epsilon_{1}} \beta^{-\gamma_{1}}\cdots  \alpha^{-\epsilon_{m}} \beta^{-\gamma_{m}}=  \alpha^{\epsilon_{1}} \beta^{\gamma_{1}}\cdots  \alpha^{\epsilon_{m}} \beta^{\gamma_{m}} \alpha^{\epsilon}.
\end{equation*}
If $\gamma_{1}\neq 0$ then $1-\epsilon_{1}=\epsilon_{1}$, which is impossible. So $\gamma_{1}=0$, and thus $m=1$, $u=\alpha^{\epsilon_{1}}$ and $1-\epsilon=2\epsilon_{1}$. If $\epsilon=1$ then $\epsilon_{1}=0$ and $u=1$, which contradicts the fact that $u\neq 1$. Hence $\epsilon=-1$, $\epsilon_{1}=1$ and $(u,\xi)\in \brak{( \alpha,x),( \alpha,x^{-1})}$, but then $(u,\xi)$ is of order $4$ by \repr{5conjclasses}(\ref{it:5conja}), which contradicts the fact that it generates a copy of $\Z$. We conclude that $P_{3}(\rp)$ has no subgroup $G$ isomorphic to $\Z_{4}\times \Z$.
\end{proof}

\begin{proof}[Second proof of \repr{p3noz4z}]
Let $g\in P_3(\rp)$ be an element of order $4$. Using the information regarding the conjugacy classes of the elements of $B_3(\rp)$ of order $4$ given in \resec{ccb3rp}, we see that $g$ is conjugate in $B_3(\rp)$ to $a^3$ or to $b^2$. But by \repr{centraliserx}(\ref{it:centraliser}), the centraliser of $g$ in $B_3(\rp)$ is finite, in particular, there is no element of $P_3(\rp)$ of infinite order that commutes with $g$.
\end{proof}

As we saw in \repr{5conjclasses}, in $P_{3}(\rp)$, there are five conjugacy classes of elements of order $4$, and four conjugacy classes of cyclic subgroups of order $4$. In~\cite[Proposition~11]{GGzeit}, we proved that if $n\geq 2$, there are $(n-2)!(2n-1)$ conjugacy classes of elements of order four in $P_n(\rp)$ (note that in the statement of that proposition, one should read `in $P_{n}(\rp)$' and not `in $B_{n}(\rp)$'). It remains to determine the number of cyclic subgroups of order $4$ in $P_n(\rp)$ if $n\geq 4$.

 
\begin{prop}
Let $n\geq 4$. The number of conjugacy classes of cyclic subgroups of order $4$ in $P_n(\rp)$ is equal to $(2n-1)(n-2)!/2$.
\end{prop}
  
\begin{proof} 
Let $z\in P_{n}(\rp)$ be of order $4$. From the comments preceding the statement of the proposition, it suffices to show that $z$ is not conjugate in $P_{n}(\rp)$ to $z^{-1}$. Let $x=a^{n}$ (resp.\ $x=b^{n-1}$), and if $u\in P_{n}(\rp)$, let $[u]$ denote its conjugacy class in $P_{n}(\rp)$. Further, let $\tau= \garside$ (resp.\ $\tau=\garside a^{-1}$), and let $c=\pi(a)= (1,2,\ldots,n)$ (resp.\ $c=\pi(b)= (1,2,\ldots,n-1)$). One may check that $\pi(\tau) \notin \ang{c}$. From the proof of~\cite[Proposition~11]{GGzeit}, there exists $\gamma\in B_{n}(\rp)$ such that $z=\gamma x \gamma^{-1}$, and so $\gamma \tau x (\gamma \tau)^{-1}= (\gamma \tau \gamma^{-1}) (\gamma x \gamma^{-1}) (\gamma \tau \gamma^{-1})^{-1}= \gamma x^{-1} \gamma^{-1}= z^{-1}$ by~\reqref{conja} and~\reqref{conjb}. Using the computation at the bottom of~\cite[p.~676]{GGzeit}, it follows that $[z]= [z^{-1}]$ if and only if $\pi(\gamma) \equiv \pi(\gamma \tau) \bmod{\ang{c}}$ in $\sn$, which is equivalent to saying that $\pi(\tau) \in \ang{c}$, which yields a contradiction. So in $P_{n}(\rp)$, $z$ and $z^{-1}$ belong to different conjugacy classes.
\end{proof}

\section{Conjugacy classes of $\quat$ in $B_n(\rp)$ and other finite subgroups in $B_n(\rp)$, $n\in \brak{3,4,5}$}\label{app:quatbnrp2} 

By~\repr{BraMap}, the number of conjugacy classes of a finite subgroup $G$ of $B_n(\rp)$ is equal to the number of conjugacy classes of $\phi(G)$ in $\mcg{\rp}$. Following~\cite{BCG}, to each finite subgroup $H$ of $\mcg{\rp}$ one may associate a collection of objects called \emph{marked signatures}. By~\cite[Theorem~3.1]{BCG}, the number of conjugacy classes in $\mcg{\rp}$ of the isomorphism class of $H$ is equal to the number of equivalence classes of marked signatures of this isomorphism class, where the notion of equivalence is defined in~\cite[p.~268]{BCG} (two marked signatures are equivalent if one may be obtained from the other by means of a cyclic permutation and/or reversal). In this appendix, we compute the number of conjugacy classes of certain finite groups of $B_{n}(\rp)$ in the following cases:
\begin{enumerate}[(a)]
\item the subgroups of $B_{n}(\rp)$ that are isomorphic to $\quat$.
\item the remaining non abstractly maximal non Abelian finite subgroups of $B_{n}(\rp)$ for $n=3,4,5$.
\end{enumerate}
  
\subsection{The conjugacy classes of $\quat$ in $B_n(\rp)$}\label{app:quat}
 
The case where the isomorphism class of a finite subgroup of $B_n(\rp)$ is $\quat$ corresponds that of the subgroup $\dih{4}$ of $\mcg{\rp}$, for which the marked signatures are given on~\cite[p.~269]{BCG}, and are as follows:
\begin{enumerate}
\item\label{it:mark1} $\sigma_1 = (0; +; [-]; \{(-)^r , (2^{2p_1} , \underline{2}, 2^{2p_2} , \underline{2}, 2^{2p_3} , \underline{2})\}) $, where $n = 4r + 2p$
\item\label{it:mark2} $\sigma_2 = (0; +; [-]; \{(-)^r , (2^{2p_1}, \underline{2}, \underline{2}, 2^{2p_2}, \underline{2}, 2^{2p_3} , \underline{2})\})$, where $n = 4r + 2p + 1$
\item\label{it:mark3} $\sigma_3 = (0; +; [-]; \{(-)^r , (2^{2p_1}, \underline{2}, \underline{2}, 2^{2p_2}, \underline{2}, \underline{2}, 2^{2p_3}, \underline{2})\})$, where $n= 4r + 2p + 2$
\item\label{it:mark6} $\sigma_6 = (0; +; [-]; \{(-)^r , (2^{2p_1}, \underline{2}, \underline{2}, 2^{2p_2}, \underline{2}, \underline{2}, 2^{2p_3}, \underline{2}, \underline{2})\})$, where $n= 4r + 2p + 3$
\end{enumerate}
where $p=p_1+p_2+p_3$ and  $r,p,p_1,p_2, p_3\geq 0$. With respect to~\cite{BCG}, we use $n$ for the number of marked points rather than $g$. To determine the number of conjugacy classes of $\quat$ in $B_n(\rp)$, we first need to solve the following equations in $r$ and $p$:
\begin{enumerate}[(1)]
\item\label{it:eq1} $n=4r+2p$
\item\label{it:eq2} $n=4r+2p+1$
\item\label{it:eq3} $n=4r+2p+2$
\item\label{it:eq4} $n=4r+2p+3$.
\end{enumerate}
It is clear from the parity that if $n$ is even (resp.\ odd), the only equations that admit solutions are~(\ref{it:eq1}) and~(\ref{it:eq3}) (resp.\ (\ref{it:eq2}) and~(\ref{it:eq4})). If $j=1,2,3,4$, then let: 
\begin{equation}
\Sigma_j(n)=\setl{p\geq 0}{\text{there exists $r\geq 0$ such that $(r,p)$ is a solution of equation~($j$)}}.
\end{equation}
If $n=2s$ is even (resp.\ $n=2s+1$ is odd), where $s\in \N$ (resp.\ $s\geq 0$), then equations~(\ref{it:eq1}) and~(\ref{it:eq3}) (resp.\ equations~(\ref{it:eq2}) and~(\ref{it:eq4})) may be written as $s=2r+p$ and $s-1=2r+p$ respectively, and so $\Sigma_1(n)$ (resp.\ $\Sigma_2(n)$) is equal to $\brak{s,s-2,\ldots, (1+(-1)^{s-1})/2}$, and $\Sigma_3(n)$ (resp.\ $\Sigma_4(n)$) is equal to $\brak{s-1,s-3,\ldots, (1+(-1)^{s-2})/2}$.
In particular, every $0\leq p\leq s$ 
is a solution of exactly one of the two equations~(\ref{it:eq1}) and~(\ref{it:eq3}) (resp.\ (\ref{it:eq2}) and~(\ref{it:eq4})).
For $p\geq 0$, consider the equation $p_1+p_2+p_3=p$, which we denote by $E_{p}$. We first describe the solutions of $E_{p}$, and for such solutions, we determine the number of equivalence classes of each the marked signatures $\sigma_1, \sigma_2, \sigma_3, \sigma_6$. 
Each solution of $E_{p}$ is represented by a triple $(p_1,p_2,p_3)$ of non-negative integers, which we refer to as a \emph{solution} of $E_{p}$. For each such solution, the set $\brak{p_1,p_2,p_3}$ will be referred to as a \emph{set solution} of $E_{p}$, and the number of set solutions of $E_{p}$ will be denoted by $I(p)$. The cardinality of a set solution is equal to $1,2$ or $3$, and the number of solutions $(p_1, p_2, p_3)$ of $E_{p}$ that yield the same set solution is equal to $6$ (resp.\ $3$) if the set solution is of cardinality $3$ (resp.\ $2$), and is equal to $1$ (resp.\ $0$) if the set solution is of cardinality $1$ and if $p$ is divisible by $3$ (resp.\ $p$ is not divisible by $3$). 
For $i=1,2,3$, let $S_i(p)$ denote the number of set solutions of $E_p$ of cardinality $i$.
Note that $I(p)= \sum_{i=1}^3 S_i(p)$. 
Let $I'(p)= \sum_{i=1}^3 i S_i(p)$. Note that: 
\begin{equation}\label{eq:defIprime}
I'(p)= 3I(p)-S_2(p)- 2S_1(p).
\end{equation}
The functions $I$ and $I'$ will be used later to describe the number of equivalence classes of the marked signatures $\sigma_i$, where $i\in \brak{1,2,3,6}$.


The aim of the following lemma is to determine $I(p)$ as a function of $p$.
To state the result, given $m\geq 0$, consider the arithmetic sequence whose first term is $m$, whose common difference is $-3$, and whose last term is $\epsilon$, where $\epsilon\in \{0,1,2\}$. Letting $\Psi(m)$ denote the sum of the terms of this arithmetic sequence, and writing $m=3w+\epsilon$, where $w\geq 0$, we have $\Psi(m)= (3w+2\epsilon)(w+1)/2= (m+\epsilon)(m+3-\epsilon)/6$. By convention, we set $\Psi(-1)=0$. Let $\floor{x}$ denote the greatest integer that is less or equal to $x$. 
 
\begin{lem}\label{lem:solset} 
Let $p\geq 0$, and let $p=2s+\delta$, where $s\geq 0$ and $\delta\in \{0,1\}$. Then:
\begin{equation*}
I(p)=\frac{(s+2)(s+1)}{2}-\Psi (s-\delta).
\end{equation*}
\begin{enumerate}[(a)]
\item\label{it:solseta} 
 $S_{1}(p)=1$ if $3$ divides $p$, and $S_{1}(p)=0$ otherwise.

\item\label{it:solsetb} 
$S_{2}(p)=\floor{p/3}+\floor{p/6}$ if $3$ divides $p$, $S_2(p)= \floor{p/3}+\floor{p/6}+1$ if $p\equiv 1 \bmod{6}$, and $S_2(p)= \floor{p/3}+\floor{p/6}+2$ otherwise.


\item\label{it:solsetc}   
$S_{3}(p)=I(p)-S_2(p)-S_1(p)$.
\end{enumerate}
\end{lem}

\begin{proof} 
Let $p\geq 0$. 
First note that $I(p)$ 
is equal to the number of solutions $(p_1, p_2, p_3)$ of $E_p$ for which $p_1\leq p_2\leq p_3$, which we now determine. Note that the maximal such value of $p_1$ is equal to $\floor{p/3}$. For $0\leq p_1 \leq \floor{p/3}$, 
let $I_{p_1}(p)$ denote the number of pairs $(p_2, p_3)$ for which $(p_1, p_2, p_3)$ is a solution of $E_p$, where $p_1\leq p_2\leq p_3$. Let $p=2s+\delta$ and $p_1=2k+\gamma$, where $s,k\geq 0$ and $\delta,\gamma\in \brak{0,1}$, and let $l=l(p_1)$ be the largest value of $p_2$ for a given value of $p_1$ for which there is a solution $(p_1, p_2, p_3)$ satisfying $p_1\leq p_2\leq p_3$. If $p_2=l$ then by maximality of $l$, we have $p_3= l+\tau$, where $\tau\in \brak{0,1}$. So for fixed $p_1$ (which is equivalent to fixing $k$), we have $p_2\in \brak{p_1,p_1+1, \ldots, l}$, and thus $I_{p_1}(p)=l-p_1+1=l-(2k+\gamma)+1$. On the other hand, since $p=p_1+p_2+p_3$, we have $2s+\delta= 2k+\gamma+ l+l+\tau$,  so $l=s-k+\Delta$, where $\Delta=(\delta-\gamma- \tau)/2$, and hence $I_{p_1}(p)= s-3k+\Delta-\gamma+1$. Now $l$ is an integer, hence:
\begin{equation}\label{eq:triples}
(\delta, \gamma, \tau)\in \brak{(0,0,0),(1,1,0), (1,0,1),(0,1,1)},
\end{equation}
for which the respective values of $\Delta$ are $0,0,0$ and $-1$. 
To determine the allowed values of $k$, note that the maximal value of $p_1$ is $\floor{p/3}$, and since $k=(p_1-\gamma)/2$, the maximal value of $k$, which we denote by $m_{\gamma}$, is equal to $\floor{\bigl(\floor{p/3}- \gamma\bigr)/2}$. 
If $\gamma=0$ then $m_0=\floor{p/6}$, and so $0\leq k\leq m_0$. So suppose that $\gamma=1$. Let $s=3u+v$, where $u\geq 0$ and $v\in \brak{0,1,2}$. 
Then for $\delta\in \brak{0,1}$, we have $p/6=(6u+2v+\delta)/6$, and thus $\floor{p/6}=u$. Hence:
\begin{equation*}
m_1=\floor{\frac{\floor{p/3}-\gamma}{2}}= \floor{\frac{\floor{(6u+2v+\delta)/3}-1}{2}}=u+\floor{\frac{\floor{(2v+ \delta)/3}-1}{2}},
\end{equation*}
from which it follows that:
\begin{equation*}
\text{$m_1=\begin{cases}
\floor{p/6}-1 & \text{if $v=0,1$}\\
\floor{p/6} & \text{if $v=2$}
\end{cases}$ if $\delta=0$, and $m_1=\begin{cases}
\floor{p/6}-1 & \text{if $v=0$}\\
\floor{p/6} & \text{if $v=1,2$}
\end{cases}$ if $\delta=1$,}
\end{equation*}
and $0\leq k\leq m_1$ in this case. 
For the maximal value of $k$, the number of solutions is either $1,2$ or $3$. We now obtain the given formula for $I(p)$ by separating the cases where $p$ is even or odd:
\begin{enumerate}[(i)]
\item\label{it:peven} if $p$ is even then $\delta=0$. If $\gamma=0$ (resp.\ $\gamma=1$) then $\Delta=0$ (resp.\ $\Delta=-1$) from~\reqref{triples}, and so for fixed $p_1$, which is equivalent to fixing $k$, $I_{p_1}(p)=l-2k+1=s-3k+1$ (resp.\ $I_{p_1}(p)=l-2k=s-3k-1$). 

\item if $p$ is odd then $\delta=1$, and it follows from~\reqref{triples} that $\Delta=0$. If $\gamma=0$ (resp.\ $\gamma=1$), for fixed $p_1$, we have $I_{p_1}(p)=l-2k+1=s-3k+1$ (resp.\ $I_{p_1}(p)= l-2k=s-3k$).
\end{enumerate}
Combining the two cases, we obtain $I_{p_1}(p)=s-3k+1$ if $\gamma=0$, and $I_{p_1}(p)=s-3k+\delta-1$ if $\gamma=1$, and thus:
\begin{equation*}
I(p)=\sum_{k_0=0}^{m_0} (s-3k_0+1)+\sum_{k_1=0}^{m_1} (s-3k_1+\delta-1)= \Psi(s+1)+\Psi(s+\delta-1)= \frac{(s+2)(s+1)}{2}-\Psi(s-\delta),
\end{equation*}
which yields the required formula.

We now prove parts~(\ref{it:solseta})--(\ref{it:solsetc}) of the statement. For part~(\ref{it:solseta}), if $\brak{p_1,p_2,p_3}$ is a set solution of $E_p$ of cardinality $1$ then $p_1=p_2=p_3$, and so $p=3p_1$. Thus $S_1(p)=1$ if $p$ is divisible by $3$, and $S_1(p)=0$ otherwise. 

To prove part~(\ref{it:solsetb}), let $\brak{p_1,p_2,p_3}$ be a set solution of $E_p$ of cardinality $2$ such that $p_1 \leq p_2 \leq p_3$. Then either $p_1 = p_2 < p_3$, in which case $p_3$ and $p$ have the same parity, or $p_1 < p_2 = p_3$, in which case $p_1$ and $p$ have the same parity. First suppose that $3$ divides $p$, so there exists $v\geq 0$ such that $p=3v$. Thus $p$ and $v$ have the same parity, and hence $p_1$ and $v$ have the same parity. 
\begin{enumerate}[(i)]
\item For the case $p_1 = p_2 < p_3$, we have $0\leq p_1 \leq v-1$, since if $p_1=p_2=v$, then $p_3=v$ also, and thus the number of solutions is equal to $v=\floor{p/3}$.

\item For the case $p_1 < p_2 = p_3$, let $v=2u+\epsilon$, where $n\geq 0$ and $\epsilon \in \brak{0,1}$. Then $p=3v= 6u+ 3\epsilon$. 
 Let $\epsilon=0$ (resp.\ $\epsilon=1$). Since $p_1$ and $v$ have the same parity, they are both even (resp.\ odd), and so $0\leq p_1 \leq v-2$ (resp.\ $1\leq p_1 \leq 2u-1$) because $v-1$ is odd (resp.\ $2u$ is even), and if $p_1=v$ then $p_1=p_2=p_3=v$ (resp.if $p_1=2u+1$ then $p_1=p_2=p_3=2u+1$), which is not permitted. So the number of solutions is equal to $v/2$ (resp.\ to $u$), which in both cases, may be seen to be equal to $\floor{p/6}$.
\end{enumerate}
So if $3$ divides $p$, we have $S_2(p)= \floor{p/3}+\floor{p/6}$ as required.

 
Now suppose that $3$ does not divide $p$, so $p=3v+\epsilon$, where $v\geq 0$ and $\epsilon\in \{1,2\}$. 
\begin{enumerate}[(i)]
\item For the case $p_1 = p_2 < p_3$, we have $0\leq p_1 \leq v$, and the case is realised, since if $p_1=p_2=v$ then $p_3=v+\epsilon$. So the number of solutions in this case is equal to $v+1$, which may be seen to be equal to $\floor{p/3}+1$. 

\item 
For the case $p_1 < p_2 = p_3$, let $v=2u+\delta$, where $u\geq 0$ and $\delta\in \brak{0,1}$. So $p=6u+ 3\delta+ \epsilon$. Since $p$ and $p_1$ have the same parity, $p_1$ has the same parity as $3\delta+ \epsilon$. Let $\delta=0$. If $\epsilon=1$ (resp.\ $\epsilon=2$), then $1 \leq p_1 \leq 2u-1$ (resp.\ $0 \leq p_1 \leq 2u$), and if $p_1=2u-1$ (resp.\ $p_1=2u$) then $p_2=p_3=2u+1$. It follows that there are $u$ solutions (resp.\ $u+1$ solutions), which may be seen to be equal to $\floor{p/6}$ (resp.\ to $\floor{p/6}+1$).

Now let $\delta=1$. If $\epsilon=1$ (resp.\ $\epsilon=2$), then $0 \leq p_1 \leq 2u$ (resp.\ $1 \leq p_1 \leq 2u+1$), and if $p_1=2u$ (resp.\ $p_1=2u+1$) then $p_2=p_3=2u+2$. It follows that there are $u+1$ solutions, which may be seen to be equal to $\floor{p/6}+1$ in both cases.
\end{enumerate}
So if $3$ does not divide $p$, we have $S_2(p)= \floor{p/3}+\floor{p/6}+1$ if $p\equiv 1 \bmod{6}$, and $S_2(p)= \floor{p/3}+\floor{p/6}+2$ otherwise. This completes the proof of part~(\ref{it:solsetb}). Finally, part~(\ref{it:solsetc}) follows from parts~(\ref{it:solseta}) and~(\ref{it:solsetb}), and the fact that $I(p)= \sum_{i=1}^3 S_i(p).$
 \end{proof}

The following results will allow us to give formul\ae\ for $I(n)$ as a function of $n$.

\begin{lem}\label{lem:countIn}
For all $n\geq 0$, $I(n+6)=I(n)+n+6$, 
and if $n=6m+\epsilon$, where $m\geq 0$ and $\epsilon=0,1,\ldots,5$, then $I(n)=I(\epsilon) + 3m(m+1)+m\epsilon$.
\end{lem}

\begin{proof}
Let $n\geq 0$. To prove the first part of the statement, let $n=2s+\delta$ and $s=3k+\epsilon$, where $\delta\in \{0,1\}$ and $\epsilon\in \{0,1,2\}$. Then $n+6=2s+\delta+6= 6k+2\epsilon +\delta +6$. From the definition of $\Psi$, we have $\Psi(s+3-\delta)-\Psi (s-\delta)=s+3-\delta=3k+\epsilon+3-\delta$, and it follows from \relem{solset} that: 
\begin{equation*}
I(n+6)-I(n)= \begin{pmatrix}
3k+5+\epsilon\\
2
\end{pmatrix}+ \Psi(s+3-\delta)- \begin{pmatrix}
3k+2+\epsilon\\
2
\end{pmatrix}- \Psi(s-\delta)=6k+6+2\epsilon +\delta=n+6,
\end{equation*}
which proves the first part. To prove the second part, let $n=6m+\epsilon$, where $m\geq 0$ and $\epsilon=0,1,\ldots,5$. We argue by induction on $m$. If $m=0$, then the result is clearly true. If the equation holds for some $m\geq 0$ then by the first part, we have:
\begin{align*}
I(6(m+1)+\epsilon)&= I(6m+\epsilon) + 6m+\epsilon+6= I(\epsilon) + 3m(m+1)+m\epsilon + 6m+\epsilon+6\\
&= I(\epsilon) + 3(m+1)(m+2)+(m+1)\epsilon,
\end{align*}
which yields the result.
%
\end{proof}

For $n=0,1,\ldots,5$, the values of $I(n)$ may be computed directly, and are given by $I(0)=1$, and $I(n)=n$ for all $n=1,\ldots,5$. The formula given in \relem{countIn} then allows us to determine $I(n)$ for all $n\geq 0$.

We are now able to give an explicit formula for $I'(n)$, and obtain relations between $I$ and $I'$.

\begin{lem}\mbox{}\label{lem:contIprime}
\begin{enumerate}
\item\label{it:contIprimea} Let $n\geq 0$. If $n$ is even (resp.\ odd) then $I'(n)=(n+2)^2/4$ (resp.\ $I'(n)= (n+1)(n+3)/4$).

\item\label{it:contIprimeb} For all $n\geq 1$, $I'(n-1)+I(n)=I(2n)$ and $I(n-1)+I'(n)=I(2n+1)$.
\end{enumerate}
\end{lem}

\begin{proof}\mbox{}
\begin{enumerate}
\item Let $n=6m+\epsilon$, where $m\geq 0$ and $\epsilon \in \brak{0,1,\ldots,5}$. Recall from~\reqref{defIprime} that $I'(n)=3I(n)-S_2(n)-2S_1(n)$. One may check that $\floor{n/6}+\floor{n/3}=3m+\floor{\epsilon/3}$, and using \relem{solset}(\ref{it:solseta}) and~(\ref{it:solsetb}), we see that $2S_1(n)+S_2(n)=3m+\floor{\epsilon/3}+\rho$, where $\rho=1$ if $\epsilon=1$ and $\rho=2$ otherwise. Using \relem{countIn}, it follows that:
\begin{equation*}
I'(n)=3(I(\epsilon) + 3m(m+1)+m\epsilon)-3m -\floor{\epsilon/3}- \rho= 9m^2+(6+3\epsilon)m+ 3I(\epsilon)- \floor{\epsilon/3}-\rho.
\end{equation*}
We now analyse the cases $\epsilon=0,1,\ldots,5$ in turn:
\begin{enumerate}
\item $\epsilon=0$: then $I'(n)= 9m^2+6m+1=(3m+1)^2$.
\item $\epsilon=1$: then $I'(n)= 9m^2+9m+2=(3m+1)(3m+2)$.
\item $\epsilon=2$: then $I'(n)= 9m^2+12m+4=(3m+2)^2$.
\item $\epsilon=3$: then $I'(n)= 9m^2+15m+6=(3m+2)(3m+3)$.
\item $\epsilon=4$: then $I'(n)= 9m^2+18m+9=(3m+3)^2$.
\item $\epsilon=5$: then $I'(n)= 9m^2+21m+12=(3m+3)(3m+4)$.
\end{enumerate}
If $n$ is even (resp.\ odd) then we observe that $I'(n)=(n+2)^2/4$ (resp.\ $I'(n)= (n+1)(n+3)/4$) as required.
\item Let $n\geq 1$. Note that the relations may be expressed in the form $I(2n+\epsilon)-I(n-\epsilon)=I'(n-1+\epsilon)$, where $\epsilon\in \brak{0,1}$. We argue by induction on $n$.
One may check (see Table~\ref{tab:IGamman}) that the relations hold for $n=1,\ldots,6$. 
Suppose that they hold for all integers up to and including some $n\geq 6$, and let us show that $I(2(n+1)+\epsilon)-I(n+1-\epsilon)=I'(n+\epsilon)$. By part~(\ref{it:contIprimea}), if $n+\epsilon$ is even (resp.\ odd) then $I'(n+\epsilon)-I'(n-6+\epsilon) = ((n+2+\epsilon)/2)^2-((n-4+\epsilon)/2)^2$ (resp.\ $I'(n+\epsilon)-I'(n-6+\epsilon) =(n+1+\epsilon)(n+3+\epsilon)/4-(n-5+\epsilon)(n-3+\epsilon)/4$), and
in both cases we may check that $I'(n+\epsilon)-I'(n-6+\epsilon)=3(n+\epsilon-1)$. By \relem{countIn}, we have $I(n+1-\epsilon)=I(n-5-\epsilon)+n+1-\epsilon$ and $I(2(n+1)+\epsilon) =I(2n-10+\epsilon) +4n-2+2\epsilon$, and hence:
\begin{align*}
I(2(n+1)+\epsilon)- I(n+1-\epsilon) &= I(2n-10+\epsilon) +4n-2+2\epsilon-(I(n-5-\epsilon)+n+1-\epsilon)\\
&=  I(2n-10+\epsilon)-I(n-5-\epsilon) +3n-3+3\epsilon\\
&= I(2(n-5)+\epsilon)-I(n-5-\epsilon) + I'(n+\epsilon)-I'(n-6+\epsilon)\\
&=I'(n+\epsilon),
\end{align*}
using the fact that $I(2(n-5)+\epsilon)-I(n-5-\epsilon)=I'(n-6+\epsilon)$ by the induction hypothesis, which yields the result.\qedhere
%
%
\end{enumerate}
\end{proof}


For a given value of $p$, we now compute the number of equivalence classes of the marked signatures $\sigma_i$ mentioned at the beginning of the section, where $i\in \brak{1,2,3,6}$. 

\begin{lem}\label{lem:countingcc}
Let $n\geq 0$. If $p\in \Sigma_1(n)$ or $p\in \Sigma_4(n)$ (resp.\ $p\in \Sigma_2(n)$ or $p\in \Sigma_3(n)$) then the number of equivalence classes of the marked signatures $\sigma_1$ and $\sigma_6$ (resp.\ $\sigma_2$ and $\sigma_3$) is equal to $I(p)$
(resp.\ to $I'(p)$).
\end{lem}

\begin{proof} 
Let $n,p\geq 0$, where $p\in \Sigma_i(n)$ for some $j=1,\ldots,4$. Recall from the beginning of this Appendix that two marked signatures are equivalent if one may be obtained from the other by means of cyclic permutations and/or reversal. Up to these operations and a permutation of the $p_i$, for each set solution $\{p_1, p_2,p_3\}$ of $E_p$, the marked signatures $\sigma_i$, where $i\in \brak{1,2,3,6}$, may be written as $(2^{2p_1}, a, 2^{2p_2}, b, 2^{2p_3} , b)$, where $a$ and $b$ are one of the two sequences $\underline{2}, \underline{2}$ or $\underline{2}$, and for each of the four possibilities for  $(a,b)$, there 
are $6$ marked signatures of the form $\sigma_{i,\tau}= (2^{2p_{\tau(1)}}, a, 2^{2p_{\tau(2)}}, b, 2^{2p_{\tau(3)}}, b)$, where $\tau\in \sn[3]$. If $i\in \brak{1,6}$ then $a=b$, and the six marked signatures $\sigma_{i,\tau}$, where $\tau\in \sn[3]$, are pairwise equivalent, hence each set solution $\{p_1, p_2,p_3\}$ of $E_p$ gives rise to a single equivalence class of marked signatures, and so the number of equivalence classes of marked signatures is equal to $\sum_{i=1}^3 S_i(p)$, which is equal to $I(p)$.

Now suppose that $i\in \brak{2,3}$, so $a\neq b$. If $\tau\in \sn[3]$, the elements in the same equivalence class as $\sigma_{i,\tau}$ may be obtained by permuting cyclically the six elements of $\sigma_{i,\tau}$ by powers of the $6$-cycle $\alpha=(1,2, 3, 4, 5, 6)$, and/or by reversing the order of these elements, which corresponds to the permutation $\beta=(1,6)(2,5)(3,4)$. Noting that the conjugate of $\beta \alpha \beta^{-1}= \alpha^{-1}$, we see that $\ang{\alpha, \beta}\cong \dih{12}$, so the elements of $\ang{\alpha, \beta}$ are of the form $\alpha^{i} \beta^{j}$ for $0 \leq i \leq 5$ and $0\leq j \leq 1$. From this, one may check that the marked signatures $(2^{2p_{\tau(1)}}, a, 2^{2p_{\tau(2)}}, b, 2^{2p_{\tau(3)}}, b)$ and $(2^{2p_{\tau(2)}}, a, 2^{2p_{\tau(2)}}, b, 2^{2p_{\tau(3)}}, b)$ are equivalent. Let $\{p_1, p_2,p_3\}$ be a set solution of $E_p$ of cardinality $l$, where $l\in \brak{1,2,3}$. If $l=1$ then clearly there is a single equivalence class of marked signatures. So suppose that $l=2$. Permuting the $p_i$ if necessary, we may suppose that $p_1=p_2$, and in this case, there are two equivalence classes of marked signatures, represented by $(2^{2p_1}, a, 2^{2p_2}, b, 2^{2p_3}, b)$ and $(2^{2p_1}, a, 2^{2p_3}, b, 2^{2p_2}, b)$. Finally if $l=3$ then there are three equivalence classes of marked signatures, represented by $(2^{2p_1}, a, 2^{2p_2}, b, 2^{2p_3}, b)$ and $(2^{2p_3}, a, 2^{2p_1}, b, 2^{2p_2}, b)$ and $(2^{2p_2}, a, 2^{2p_3}, b, 2^{2p_1}, b)$. It follows that the number of equivalence classes of marked signatures is equal to $S_1(p)+ 2S_2(p) + 3S_3(p)$, which is equal to $I'(p)$.
\end{proof}

For $n\geq 2$, let $\Gamma_n$ denote the number of conjugacy classes of subgroups of $B_n(\rp)$ isomorphic to $\quat$. 
In what follows, we will give an explicit formula for $\Gamma_n$, and to this end, we set $\Gamma_0=\Gamma_1=0$. Making use of the comments before Lemma~\ref{lem:solset} as well as Lemmas~\ref{lem:solset} and~\ref{lem:countingcc}, the following result allows us to compute $\Gamma_n$. 


\begin{thm}\label{th:formGamman}
Let $n\geq 2$.
\begin{enumerate}
\item\label{it:formGammana} if $n$ is even (resp.\ odd) then $\Gamma_n= \sum_{p\in \Sigma_1(n)} I(p)+\sum_{p\in \Sigma_3(n)} I'(p)$ (resp.\ $\Gamma_n=\sum_{p\in \Sigma_4(n)} I(p)+\sum_{p\in \Sigma_2(n)} I'(p)$).
\item\label{it:formGammanb} $\Gamma_{n+4}= \Gamma_{n}+ I(n+4)$.
\end{enumerate}
\end{thm}

\begin{proof} 
Let $n\geq 2$ be even (resp.\ odd). 
\begin{enumerate}
\item Using~\cite[Theorem~3.1]{BCG} 
and \repr{BraMap}, $\Gamma_n$ is equal to the number of equivalence classes of marked signatures $\sigma_1$ and $\sigma_3$ (resp.\ $\sigma_2$ and $\sigma_6$), for the values of $p$ belonging to $\Sigma_1(n)$ and $\Sigma_3(n)$ (resp.\ to $\Sigma_2(n)$ and $\Sigma_4(n)$) respectively. The result then follows from \relem{countingcc}.
\item 
From the beginning of the section, we see that $\Sigma_1(n+4)=\Sigma_1(n) \cup \brak{2+n/2}$ and $\Sigma_3(n+4)=\Sigma_3(n) \cup \brak{1+n/2}$ (resp.\ $\Sigma_2(n+4)=\Sigma_2(n) \cup \brak{(n+3)/2}$ and $\Sigma_4(n+4)=\Sigma_4(n) \cup \brak{(n+1)/2}$). It follows from part~(\ref{it:formGammana}) and \relem{contIprime}(\ref{it:contIprimeb}) that $\Gamma_{n+4}=\Gamma_{n}+ I(2+n/2)+ I'(1+n/2)= \Gamma_{n}+I(n+4)$ (resp.\ $\Gamma_{n+4}=\Gamma_{n}+ I((n+1)/2)+ I'((n+3)/2)= \Gamma_{n}+I(n+4)$) as required.\qedhere
\end{enumerate}
\end{proof}


Using the above results, in Table~\ref{tab:IGamman}, we determine the values of $S_i(n)$, $i=1,2,3$, $I(n)$, $I'(n)$, $\Sigma_j(n)$, $j\in \brak{1,2,3,4}$, and $\Gamma_n$ for $0\leq n\leq 17$. From these computations, we obtain the following corollary that gives the number of conjugacy classes of subgroups isomorphic to $\quat$ for the values of $n=3,4,5$ studied in Sections~\ref{sec:subgroupsb3}--\ref{sec:subgroupsB5}.


\begin{cor}\label{cor:clasq8} The number of conjugacy classes of subgroups isomorphic to $\quat$ is equal to two in $B_2(\rp)$, to three in $B_3(\rp)$, to five in $B_4(\rp)$, and to six in $B_5(\rp)$.
\end{cor}

\begin{table}[h!]
\begin{center}
\resizebox{\textwidth}{!}{\renewcommand{\arraystretch}{1.4}
\setlength{\tabcolsep}{2pt}\begin{tabular}{|>{$}c<{$}|*{18}{>{$}c<{$}|}}\hline
n & 0 & 1 & 2 & 3 & 4 & 5 & 6 & 7 & 8 & 9 & 10 & 11 & 12 & 13 & 14 & 15 & 16 & 17\\ \hline
S_1(n) & 1 & 0 & 0 & 1 & 0 & 0 & 1 & 0 & 0 & 1 & 0 & 0 & 1 & 0 & 0  & 1 & 0 & 0 \\ 
\hline
S_2(n) & 0 & 1 & 2 & 1 & 3 & 3 & 3  & 4 &  5 & 4  & 6 & 6  & 6 &   7 & 8 & 7 & 9 & 9 \\ \hline
S_3(n) & 0 & 0 & 0 & 1 & 1 & 2 & 3  & 4 &  5 & 7  & 8 & 10  & 12 & 14 & 16 & 19 & 21 & 24 \\ \hline
I(n) & 1 & 1 & 2 & 3 & 4  & 5 & 7  & 8  & 10 & 12 & 14 & 16 & 19  & 21 & 24 & 27 &30  & 33 \\ \hline
I'(n) & 1 & 2 & 4 & 6 & 9 & 12 & 16 & 20 & 25  & 30  & 36 & 42   & 49  & 56 & 64 &72 & 81& 90 \\
 \hline
\Gamma_n & 0 & 0 & 2 & 3 & 5 & 6 & 9 & 11 & 15 & 18 & 23 & 27 & 34 & 39 & 47 & 55 & 64 & 72\\ \hline
\Sigma_1(n) & \brak{0} & & \{1\}  &  & \{ 0, 2 \} & & \{1, 3\} &  & \{0,2,4\} & 
 & \{1,3,5\}  & & \{0,2,4,6\}  &  & \{1,3,5,7\}   &   & \{0,2,4,6,8 \} & \\
 \hline
\Sigma_3(n) & \varnothing & & \{0\} & & \{1\}  &  & \{ 0, 2 \} & & \{1, 3\} &  & \{0,2,4\} & 
 & \{1,3,5\}  & & \{0,2,4,6\}  &  & \{1,3,5,7\}   &  \\
 \hline 
\Sigma_2(n) & & \brak{0} & & \{1\} & & \{0,2\}  &  & \{ 1, 3 \} & & \{0,2,4\} &  & \{1, 3, 5\} & 
 & \{0,2,4,6\}  & & \{1,3,5,7\}  &  & \{0,2,4,6,8\}     \\
 \hline 
\Sigma_4(n) & & \varnothing & & \{0\}  &  & \{ 1 \} & & \{0, 2\} &  & \{1,3\} & 
 & \{0,2,4\}  & & \{1, 3, 5 \}  &  & \{0, 2, 4, 6 \}   &   & \{1,3,5,7 \}   \\
 \hline 
\end{tabular}} \end{center} \caption{Values of $S_i(n)$, $i\in \brak{1,2,3}$, $I(n),I'(n)$, $\Sigma_j(n)$, $j\in \brak{1,2,3,4}$, and $\Gamma_n$ for $n\leq 17$.}\label{tab:IGamman}\end{table}

To end this section, we give an explicit formula to compute $\Gamma_n$. 
\begin{thm}\label{th:GammanQ8}
Let $n\geq 0$, let $n=6m+\epsilon$, where $m\geq 0$ and $\epsilon\in \brak{0,1,\ldots,5}$, and let $\alpha=\epsilon$ (resp.\ $\alpha=6+\epsilon$) if $m$ is even (resp.\ odd). Let 
$\epsilon'=3$ if $\epsilon=0,3$, and let $\epsilon'= \epsilon+1$ if $\epsilon=1,2,4,5$.
Then:
\begin{equation*}
\Gamma_n= \begin{cases}
\displaystyle \frac{(6m^2+3(5+\epsilon)m +2(3+3\epsilon+ \epsilon'))m}{4}+ \Gamma_{\alpha}& \text{if $m$ is even}\\
\displaystyle \frac{(6m^2+3(5+\epsilon)m +2(3+3\epsilon+ \epsilon'))m-(27+9\epsilon+2\epsilon')}{4} + \Gamma_{\alpha}& \text{if $m$ is odd.}
\end{cases}
\end{equation*}
 \end{thm}

\begin{proof}
Let $n,m,\epsilon,\epsilon'$ and $\alpha$ be as in the statement. If $n\geq 12$, then using \reth{formGamman}(\ref{it:formGammanb}), we obtain:
\begin{equation}\label{eq:Gamman12}
\Gamma_n= \Gamma_{n-4}+I(n)=\Gamma_{n-8}+I(n-4)+I(n)= \Gamma_{n-12} + I(n-8)+ I(n-4)+ I(n).
\end{equation}
If $n\geq 12$ and $\epsilon=0$ then $n-4=6(m-1)+2$ and $n-8=6(m-2)+4$, and by~\req{Gamman12} and \relem{countIn}, we see that:
\begin{align*}
\Gamma_n-\Gamma_{n-12}=& I(n-8)+ I(n-4)+ I(n)\\
=& (I(4)+3(m-2)(m-1)+4(m-2))+\\
& (I(2)+3(m-1)m+ 2(m-1))+(I(0)+3m(m+1))\\
=&9m^2-3m+3.
\end{align*}
Similarly, if $\epsilon=1,2,3,4,5$, one may check that $\Gamma_n-\Gamma_{n-12}$ is equal to $9m^{2}+2$, 
$9m^2+3m+3$, $9m^2+6m+3$, $9m^2+9m+5$ and $9m^2+12m+6$ respectively. These results may be summarised as $\Gamma_n-\Gamma_{n-12}=P(m)$, where $P(m)= 9m^2+3(\epsilon-1)m+\epsilon'$. Using the fact that $n=12(m/2)+\epsilon$ (resp.\ $n=12((m-1)/2)+ 6+\epsilon$) if $m$ is even (resp.\ odd), we obtain: 
\begin{equation}\label{eq:gammanalpha}
\Gamma_n- \Gamma_{\alpha}= 
\begin{cases}
\displaystyle 
\sum_{i=0}^{m/2-1} (\Gamma_{n-12i}- \Gamma_{n-12(i+1)}) =
\sum_{\substack{i=2 \\ \text{$i$ even}}}^{m} P(i) & \text{if $m$ is even}\\
\displaystyle 
\sum_{i=0}^{(m-3)/2} (\Gamma_{n-12i}- \Gamma_{n-12(i+1)}) = \sum_{\substack{i=3 \\ \text{$i$ odd}}}^{m} P(i) & \text{if $m$ is odd.}
\end{cases}
\end{equation}
Making use of the standard formul\ae\ for the sum of (squares of) even and odd integers, if $m$ is even, it follows from~\reqref{gammanalpha} that: 
\begin{align*}
\Gamma_n- \Gamma_{\alpha}&= \sum_{\substack{i=2 \\ \text{$i$ even}}}^{m} (9i^2+3(\epsilon-1)i+ \epsilon')\\
&= 9\left[ \frac{2\ldotp \frac{m}{2}\left( \frac{m}{2}+1 \right) 
\left( 2\ldotp \frac{m}{2}+1\right)}{3}\right] +3(\epsilon-1) \frac{m}{2}\left( \frac{m}{2}+1\right)+ \frac{m}{2}\epsilon'\\
&= \frac{(6m^2+3(5+\epsilon)m +2(3+3\epsilon+ \epsilon'))m}{4},
\end{align*}
while if $n$ is odd, we have:
\begin{align*}
\Gamma_n- \Gamma_{\alpha}&= \sum_{\substack{i=3 \\ \text{$i$ odd}}}^{m} (9i^2+3(\epsilon-1)i+ \epsilon')\\
&= 9\left[ \frac{\frac{m+1}{2}\left( 2\ldotp \frac{m+1}{2}-1 \right) 
\left( 2\ldotp \frac{m+1}{2}+1\right)}{3}-1 \right] +3(\epsilon-1)\left( \left(\frac{m+1}{2} \right)^2-1 \right) +\left( \frac{m+1}{2}-1 \right) \! \epsilon'\\
&=\frac{(6m^2+3(5+\epsilon)m +2(3+3\epsilon+ \epsilon'))m-(27+9\epsilon+2\epsilon')}{4}
\end{align*}
as required.
\end{proof}

\subsection {The number of conjugacy classes of the finite subgroups of $B_n(\rp)$ for $n\in \brak{3,4,5}$}

In this section, we apply the results of~\cite{BCG} to \repr{BraMap} to determine the number of conjugacy classes of the finite subgroups of $B_n(\rp)$ for $n\in \brak{3,4,5}$. As we saw in \resec{subgroupsb3}, this may be achieved without reference to~\cite{BCG}, but the computation here enables us to check our results. 

Let $n \in \brak{3,4,5}$. The number of conjugacy classes in $B_{n}(\rp)$ of the finite cyclic subgroups may be obtained using \reth{ccrp2}(\ref{it:ccrp2b}) and \reco{ccorder4}. The number of conjugacy classes in $B_{n}(\rp)$ subgroups isomorphic to $\quat$ has already been considered in \ref{app:quat}. 
The case of the dicyclic  subgroups $\dic{8n}$ and $\dic{8(n-1)}$ also has been considered subgroups since they are maximal.
It remains to compute the number of conjugacy classes of the dicyclic groups $\dic{4m}$ of order greater than eight and different from distinct $8(n-1)$ and $8n$, 
and the binary polyhedral groups. The dicyclic group $\dic{4m}$ in $B_n(\rp)$ corresponds to the subgroup  $\dih{2m}$ of $\mcg{\rp}$, which is denoted $D_m$ in~\cite{BCG}.

\begin{prop}\label{prop:conjclass345}\mbox{}
\begin{enumerate}
\item In $B_{3}(\rp)$, there are two conjugacy classes of subgroups isomorphic to $\dic{12}$, and a single conjugacy classes of subgroups isomorphic to $\tonestar$ or to $\oonestar$.

\item In $B_{4}(\rp)$, there are three conjugacy classes of subgroups isomorphic to $\quat[16]$, two conjugacy classes of subgroups isomorphic to $\dic{12}$, and a single conjugacy classes of subgroups isomorphic to $\tonestar$ or to $\oonestar$.

\item In $B_{5}(\rp)$, there are three conjugacy classes of subgroups isomorphic to $\quat[16]$, and two conjugacy classes of subgroups isomorphic to $\dic{20}$.
\end{enumerate}
\end{prop}

\begin{proof}
Using \rerems{finitemcgbnrp}, under the homomorphism $\phi$ of \req{bnrpmcg}, the finite subgroups of $B_{n}(\rp)$ isomorphic to $\dic{4m}$, where $m\geq 2$, $\tonestar, \oonestar$ and $\istar$ are sent to finite subgroups of $\mcg{\rp}$ isomorphic to $\dih{2m}, \an[4], \sn[4]$ and $\an[5]$ respectively. For $n\in \brak{3,4,5}$, $\istar$ is not realised as a finite subgroup of $B_{n}(\rp)$, and $\tonestar$ and $\oonestar$ are not realised as finite subgroups of $B_{5}(\rp)$. We now analyse the cases $n=3,4,5$ separately.

\begin{enumerate}
\item Let $n=3$. We consider the finite subgroups of $B_{3}(\rp)$ that are isomorphic to $\dic{12}, \tonestar$ or $\oonestar$.
\begin{enumerate}[(i)]
\item The group $\dic{12}$ corresponds to $\dih{6}$ in $\mcg[3]{\rp}$. Taking $m=3$ in~\cite[p.~269]{BCG}, the number of conjugacy classes of $\dih{6}$ in $\mcg[3]{\rp}$ is given by the number of marked signatures (up to equivalence) of the form $\sigma_1 =(0;+;[2];\{(-)^r,(2^{2p},\underline{3})\})$, where $3=3(2r+p)$, and $\sigma_3 =(0;+;[-];\{(-)^{r+1},(2^{2p},\underline{3})\})$, where $n=m(2r+p+1)$. In the first (resp.\ second) case, the only solution is $(r,p)=(0,1)$ (resp.\ $(r,p)=(0,0)$), which give rise to a single equivalence class of marked signatures in each case. So in $\mcg[3]{\rp}$, there are two conjugacy classes of subgroups isomorphic to $\dih{6}$.

\item  The group $\tonestar$ (resp.\ $\oonestar$) corresponds to $\an[4]$ (resp.\ to $\sn[4]$) in $\mcg[3]{\rp}$. By~\cite[p.~269]{BCG} (resp.~\cite[p.~270]{BCG}), the number of conjugacy classes of each of these groups is given by the number of marked signatures (up to equivalence) of the form $\sigma_2 =(0;+;[3];\{(-)^r,(2^{2p},\underline{2}, \underline{2})\})$, where $3=12r+6p+3$ (resp.\ the number of marked signatures (up to equivalence) of the form $\sigma_4 = (0;+;[-];\{(-)^r,(2^{2p_1},\underline{2},2^{2p_2},\underline{3} ,2^{2p_3},\underline{2},\underline{2})\})$, where $3= 24r + 12p + 3$). In both cases, the only solution is $(r,p)=(0,0)$, and so in $\mcg[3]{\rp}$, there is a single conjugacy class of subgroups isomorphic to $\an[4]$ (resp.\ to $\sn[4]$).
\end{enumerate}

\item Let $n=4$. We consider the finite subgroups of $B_{4}(\rp)$ that are isomorphic to $\dic{12}, \quat[16], \tonestar$ or $\oonestar$. They correspond to the subgroups of $\mcg[4]{\rp}$ that are isomorphic to $\dih{6}, \dih{8}, \an[4]$ and $\sn[4]$ respectively.

\begin{enumerate}[(i)]
\item For $\dic{12}$, the calculation is similar to that for $B_{3}(\rp)$, where we use the marked signatures of the form $\sigma_2 =(0;+;[2];\{(-)^r,(2^{2p},\underline{2}, \underline{2})\})$ where $4=3(2r+p)+1$, and $\sigma_4 = (0;+;[-];\{(-)^{r+1},(2^{2p},\underline{2}, \underline{2})\})$, where $4=3(2r+p+1)+1$. Again, in $\mcg[4]{\rp}$, there are two conjugacy classes of subgroups isomorphic to $\dih{6}$. 

\item For $\quat[16]$, we take $m=4$ in~\cite[p.~267]{BCG}. The number of conjugacy classes of $\dih{8}$ is given by the number of marked signatures (up to equivalence) of the form $\sigma_1= (0;+; [-]; \{(-)^r, (2^{2p_1} , \underline{2}, 2^{2p_2} ,\underline{2}, 2^{2p_3},\underline{4})\})$, where $4= 4(2r + p)$, and of the form $\sigma_3 = (0;+; [-]; \{(-)^r, (2^{2p_1} , \underline{2}, \underline{2}, 2^{2p_2} , \underline{2}, \underline{2} , 2^{2p_3},\underline{4})\})$, where $4 = 4(2r + p + 1)$. The respective solutions of these two equations are $(r,p)= (0,1)$ and $(r,p)= (0,0)$, and from this we obtain two and one non-equivalent marked signatures respectively. So in $\mcg[4]{\rp}$, there are three conjugacy classes of subgroups isomorphic to $\dih{8}$.

\item For $\tonestar$ and $\oonestar$, the calculations are similar to those for $B_{3}(\rp)$, where we use the marked signatures $\sigma_3 = (0;+;[-];\{(-)^{r+1},(2^{2p}, \underline{2})\})$ for $\an[4]$, where $4 = 12r + 6p + 4$, and $\sigma_3 =(0;+;[-];\{(-)^r,(2^{2p_1},\underline{2},2^{2p_2},\underline{2}, \underline{2} ,2^{2p_3},\underline{4})\}) $ for $\sn[4]$, where $4=24r+12p+4$. Once more, we see that in $\mcg[4]{\rp}$, there is a single conjugacy class of subgroups isomorphic to $\an[4]$ and to $\sn[4]$.

\end{enumerate}

\item Let $n=5$. We consider the finite subgroups of $B_{5}(\rp)$ that are isomorphic to $\quat[16]$ or $\dic{20}$.  They correspond to the subgroups of $\mcg[5]{\rp}$ that are isomorphic to $\dih{8}$ and $\dih{10}$ respectively.

\begin{enumerate}[(i)]
\item For $\quat[16]$, the calculation is similar to that for $B_{4}(\rp)$, where we use the marked signatures $\sigma_4 = (0;+; [-]; \{(-)^r, (2^{2p_1} ,\underline{2}, 2^{2p_2} , \underline{2}, 2^{2p_3} , \underline{2}, \underline{2})\})$, where $5 = 4(2r + p) + 1$, and $\sigma_6 = (0;+; [-]; \{(-)^r, (2^{2p_1},\underline{2}, \underline{2}, 2^{2p_2} ,\underline{2}, \underline{2}, 2^{2p_3} ,\underline{2}, \underline{2})\}) $, where $5= 4(2r + p + 1) + 1$. Once again, we see that in $\mcg[5]{\rp}$, there are three conjugacy classes of subgroups isomorphic to $\quat[16]$.

\item The group $\dic{20}$ corresponds to $\dih{10}$ in $\mcg[5]{\rp}$. The calculation is similar to that for $\dih{6}$ in $\mcg[3]{\rp}$, and in the same manner, we see that there are two conjugacy classes of subgroups isomorphic to $\dih{10}$.\qedhere
\end{enumerate}
\end{enumerate}
\end{proof}

\begin{rem}
The groups $\oonestar$ and $\istar$ are maximal in $B_n(\rp)$ for all values of $n$ for which they exist. For $n\leq 5$, there is at most  
one conjugacy class of $\oonestar$. This is not case in general. By~\cite[pp.~269-70]{BCG} (resp.\ by~\cite[p.~270]{BCG}), the smallest integer for which $B_n(\rp)$ admits a copy of $\oonestar$ (resp.\ of $\istar$) is $n=3$ (resp.\ $n=6$), and there is a single conjugacy class. 
The smallest integer for which there exists more than one conjugacy class of $\oonestar$ (resp.\ of $\istar$) in $B_n(\rp)$ is $n=12$ (resp.\ $n=30$), and there are three conjugacy classes in each case.
\end{rem}

\end{document}